\documentclass[12pt]{amsart}
\usepackage{amssymb,latexsym,amsmath,amsfonts}
\usepackage{supertabular}
\usepackage{enumerate}
\usepackage[mathscr]{eucal}
\usepackage{enumitem}
\usepackage{hyperref}

\usepackage[nameinlink]{cleveref}

\allowdisplaybreaks

\numberwithin{equation}{section}
\theoremstyle{definition}
\newtheorem{definition}{Definition}[section]
\newtheorem{example}[definition]{Example}
\newtheorem{question}[definition]{Question}

\theoremstyle{remark}
\newtheorem{remark}[definition]{Remark}

\theoremstyle{plain}
\newtheorem{proposition}[definition]{Proposition}
\newtheorem{theorem}[definition]{Theorem}
\newtheorem{lemma}[definition]{Lemma}
\newtheorem{result}[definition]{Result}
\newtheorem{corollary}[definition]{Corollary}

\usepackage[usenames]{color}

\definecolor{Red}{rgb}{1,0,0}
\definecolor{Magenta}{rgb}{0.69,0,0.83}
\definecolor{Green}{rgb}{0.117,0.706,0.314}
\definecolor{Grey}{rgb}{0.8,0.8,0.8}
\definecolor{Blue}{rgb}{0.1,0.1,1}

\renewcommand{\l}{\lambda}
\renewcommand{\k}{\kappa}

\newcommand{\lk}{(\lambda, \kappa)}

\newcommand{\rprt}{\mathsf{Re}}
\newcommand{\iprt}{\mathsf{Im}}
\newcommand{\distance}{\mathrm{dist}}

\newcommand{\clos}[1]{\overline{#1}}

\newcommand{\wt}{\widetilde}
\newcommand{\wh}{\widehat}

\newcommand{\diam}{\mathsf{diam}}

\newcommand*{\defeq}{\mathrel{\vcenter{\baselineskip0.5ex \lineskiplimit0pt \hbox{\scriptsize.}\hbox{\scriptsize.}}}=}

\newcommand{\norm}[1]{\lVert #1 \rVert}
\newcommand*{\rom}[1]{\expandafter\@slowromancap\romannumeral #1@}

\newcommand{\e}{\acute{e}}

\newcommand{\bdy}{\partial}
\newcommand{\OM}{\Omega}

\newcommand{\C}{\mathbb{C}} 
\newcommand{\R}{\mathbb{R}}
\newcommand{\N}{\mathbb{N}}
\newcommand{\D}{\mathbb{D}}
\newcommand{\Z}{\mathbb{Z}}
\newcommand{\K}{\mathbb{K}}

\usepackage{tikz}
\usepackage{graphicx}

\tikzset{
  region diagram/.style={
    x=0.0085cm,
    y=-0.0085cm,
    line width=0.9pt,
    line cap=round,
    line join=round,
    every node/.style={font=\Large,inner sep=2pt}
  }
}

\newcommand{\domainoutline}{%
  \coordinate (upper-left) at (47,379);
  \coordinate (middle-left) at (63,485);
  \coordinate (lower-left) at (102,651);
  \coordinate (inner-left) at (140,765);
  \coordinate (inner-right) at (622,840);
  \coordinate (lower-right) at (754,725);
  \coordinate (global-right) at (777,704);
  \coordinate (middle-right) at (870,628);
  \coordinate (upper-right) at (961,576);
  \draw (34,30)
    .. controls (29,147) and (31,254) .. (upper-left)
    .. controls (51,416) and (56,451) .. (middle-left)
    .. controls (74,545) and (86,598) .. (lower-left)
    .. controls (113,693) and (126,729) .. (inner-left)
    .. controls (185,886) and (231,974) .. (350,974)
    .. controls (443,974) and (538,909) .. (inner-right)
    .. controls (666,803) and (710,763) .. (lower-right)
    .. controls (762,718) and (769,711) .. (global-right)
    .. controls (808,677) and (839,651) .. (middle-right)
    .. controls (900,608) and (930,590) .. (upper-right)
    .. controls (1048,535) and (1136,514) .. (1221,502);
}

\newcommand{\commonregions}{%
  \draw (upper-left)
    .. controls (130,301) and (210,247) .. (288,280)
    .. controls (348,302) and (393,367) .. (467,378)
    .. controls (535,389) and (612,359) .. (677,363)
    .. controls (819,364) and (918,460) .. (upper-right);
  \draw (inner-left)
    .. controls (201,704) and (242,752) .. (288,771)
    .. controls (351,799) and (408,753) .. (472,765)
    .. controls (532,774) and (579,800) .. (inner-right);
  \node at (785,447) {$U$};
  \node at (502,818) {$W$};
}

\usepackage{subcaption}

\definecolor{combLColor}{HTML}{8E44AD}
\definecolor{visred}{HTML}{E74747}
\definecolor{visblue}{HTML}{58A8EF}
\definecolor{viscyan}{HTML}{25A9BD}

\newcommand{\CombA}{0.72}
\newcommand{\CombB}{0.65}
\newcommand{\CombC}{0.58}

\newcommand{\CombDomainPanel}[2]{%
  \begin{tikzpicture}[
    x=.01cm,y=-.01cm,
    line width=1pt,line cap=round,line join=miter,
    every node/.style={font=\large,inner sep=0pt}
  ]
    \path[use as bounding box] (-4,-4) rectangle (534,534);
    \draw (24,22) rectangle (509,507);
    \foreach \tooth in {267,159,114,80,53,44,36,29}{
      \draw (\tooth,507) -- (\tooth,270);
    }
    \ifnum#2>0\relax
      \draw[combLColor] (24,270) -- (24,507);
      \node[text=combLColor,font=\normalsize,anchor=east]
        at (18,388.5) {$L$};
    \fi
    \ifnum#1>1\relax
      \pgfmathsetmacro{\CombAX}{24+485*\CombA}
      \pgfmathsetmacro{\CombAY}{507-485*\CombA}
      \pgfmathsetmacro{\CombUX}{24+485*0.52}
      \pgfmathsetmacro{\CombUY}{507-485*(\CombA+\CombB)/2}
      \draw[visred]
        (24,\CombAY) -- (\CombAX,\CombAY) -- (\CombAX,507);
      \node[text=visred,font=\normalsize] at (\CombUX,\CombUY) {$U$};
    \fi
    \ifnum#1>2\relax
      \pgfmathsetmacro{\CombBX}{24+485*\CombB}
      \pgfmathsetmacro{\CombBY}{507-485*\CombB}
      \pgfmathsetmacro{\CombWZeroX}{24+485*0.46}
      \pgfmathsetmacro{\CombWZeroY}{507-485*(\CombB+\CombC)/2}
      \draw[visblue]
        (24,\CombBY) -- (\CombBX,\CombBY) -- (\CombBX,507);
      \node[text=visblue,font=\scriptsize]
        at (\CombWZeroX,\CombWZeroY) {$W_1$};
    \fi
    \ifnum#1>3\relax
      \pgfmathsetmacro{\CombCX}{24+485*\CombC}
      \pgfmathsetmacro{\CombCY}{507-485*\CombC}
      \pgfmathsetmacro{\CombWX}{24+485*(0.5+\CombC)/2}
      \pgfmathsetmacro{\CombWY}{507-485*(0.5+\CombC)/2}
      \draw[viscyan]
        (24,\CombCY) -- (\CombCX,\CombCY) -- (\CombCX,507);
      \node[text=viscyan,font=\scriptsize] at (\CombWX,\CombWY) {$W$};
    \fi
    \node at (418,70) {$\Omega$};
  \end{tikzpicture}%
}

\begin{document}
\title[Intrinsic visibility and Denjoy-Wolff domains]{An intrinsic notion of visibility and classification of planar Denjoy-Wolff domains}
\author{Nirmal Rawat and Amar Deep Sarkar}

\address{NR: Department of mathematics, Indian Institute of Technology Bhubaneswar, Argul 752051, India}
\email{a24ma09005@iitbbs.ac.in}
    
\address{ADS: Department of Mathematics, Indian Institute of Technology Bhubaneswar, Argul 752051, India}
\email{amar@iitbbs.ac.in}

\keywords{Kobayashi distance, intrinsic visibility, localization result, Busemann horosphere, 
Denjoy-Wolff property, H-limit, Condition~(M$\star$), Gromov hyperbolic space, Condition $(\star\star)$.}

\subjclass[2020]{Primary: 32F45, 30D05, 30D40; Secondary: 32Q45, 53C23.}

\begin{abstract}
In this article, we introduce an intrinsic notion of visibility with respect to the Kobayashi distance that does not require any boundary. We establish localization results for the Kobayashi distance on Kobayashi hyperbolic domains in $\C^d$ that satisfy this intrinsic visibility property. As an application of these localization results, we compare the Busemann horospheres of a domain with those of a subdomain.

Recently, Benini and Bracci characterized simply connected planar domains with the Denjoy-Wolff property in terms of Busemann horospheres. More precisely, they proved that a proper simply connected planar domain has the Denjoy-Wolff property if and only if the closure of every Busemann horosphere meets the boundary of the domain in exactly one point. We introduce Condition~(M$\star$), a condition of attaching a proper disk to the domain at a boundary point. Using our localization results together with this comparison of Busemann horospheres, we extend the characterization of Benini and Bracci to domains satisfying Condition~(M$\star$) at every boundary point, and to domains satisfying it outside a totally disconnected subset of the boundary. In particular, we obtain a characterization for finitely connected domains and for a large class of infinitely connected domains having the Denjoy-Wolff property. These characterizations yield a variety of domains with the Denjoy-Wolff property, which we illustrate through example.
\end{abstract}

\maketitle
\section{Introduction and statement of results}

In this article, we introduce an intrinsic notion of visibility with respect to the Kobayashi distance that does not require any boundary. We show that this intrinsic notion of (weak) visibility is more general than the standard (weak) visibility, which is defined in terms of a boundary and has been studied in \cite{BZ2017, Bharali_Maitra, BNT2022, CMS2024, BZ2023, BGNT2024, Masanta, CGMS2024}. Using intrinsic visibility, we obtain several results that hold even for domains without standard visibility: localization of the Kobayashi distance, a local-global visibility equivalence, and a classification of the planar Denjoy-Wolff domains (see Definition~\ref{Def:DWP} for the Denjoy-Wolff property). As an application, we show that certain infinitely connected domains do not have visibility in the standard sense but still have the Denjoy-Wolff property, and that this follows from intrinsic visibility. This phenomenon was first observed by Benini and Bracci for simply connected domains, as we now recall.

Benini and Bracci showed that standard visibility does not capture the Denjoy-Wolff property; that is, standard visibility is not a necessary condition for a domain to have the Denjoy-Wolff property \cite[Corollary~1.6]{Benini-Bracci}. They constructed a non-convex simply connected planar domain that has the Denjoy-Wolff property but whose boundary is not locally connected, so it cannot have the standard visibility property with respect to the Euclidean boundary. Moreover, they classified all hyperbolic simply connected domains in the complex plane with the Denjoy-Wolff property. Their approach is to transfer the problem to the unit disk $\mathbb{D}$ via the Riemann map. They showed that such a domain has the Denjoy-Wolff property if and only if the Riemann map has an H-limit (see Definition~\ref{Def:H-limit}) at every point of $\partial\mathbb{D}$, and if and only if, for every geodesic ray, the closure in the Riemann sphere of the associated Busemann horosphere meets the boundary of the domain in exactly one point. Since the boundary of such a domain can be very irregular and need not be locally connected, the domain need not have the standard visibility property. Thus the standard notion of visibility is inadequate for classifying Denjoy-Wolff domains.

In this article, we show that intrinsic visibility allows us to classify a large class of Denjoy-Wolff domains, including finitely connected domains. A key observation is that, in the presence of intrinsic visibility, the localization of the Kobayashi distance yields information about Busemann horospheres. More precisely, localization allows us to compare the Busemann horospheres of the global domain with those of a local domain, even when the domain does not have the standard visibility property (see Theorem~\ref{Thm:Geo_horosphere_of_small_domain_is_subset_of_geo_horosphere_of_big_domain} and Theorem~\ref{Thm:Geo_horosphere_of_big_domain_is_subset_of_geod_of_small_domain_upto_some_compact_set_and-constant}).

To produce a large class of domains that have intrinsic visibility but need not have standard visibility, we introduce a condition of attached discs at a given boundary point of the domain, which we call Condition~(M$\star$) (see Definiton~\ref{Def:Condtion(M*)}). We show that Condition~(M$\star$) is equivalent to another condition, Condition~($\star$) (see Definition~\ref{Def:Condtion(*)}). We show that Condition~(M$\star$) is satisfied by simply connected domains, finitely connected domains, and infinitely connected domains whose boundary components are bounded, uniformly separated and non-degenerate. For domains satisfying Condition~(M$\star$), we give a complete characterization of the Denjoy-Wolff property in terms of Busemann horospheres and H-limits (see Theorem~\ref{Thm:Characterization_of_Condition(*)}), analogous to the one obtained by Benini and Bracci for simply connected domains.

We now recall the notion of standard visibility property with respect to the Kobayashi distance introduced by Bharali and Zimmer \cite{BZ2017,BZ2023}. Roughly speaking, the visibility property asserts that geodesics whose endpoints lie sufficiently close to two distinct boundary points of $\OM$ must bend uniformly into $\OM$. Bharali and Zimmer called it the visibility property with respect to the Kobayashi distance $K_\OM$, by analogy with the notion of visibility introduced by Eberlein and O'Neill \cite{EON}. Since $(\OM,K_\OM)$ need not be a geodesic space in general, Bharali and Zimmer introduced a class of curves with respect to the Kobayashi distance, which they called $\lk$-almost-geodesics. We refer the reader to Section~\ref{Subsec:Kob_dist} for the definition of $\lk$-almost-geodesics and to \cite{BZ2023, CGMS2024} for details on the end compactification.
\begin{definition}
Let $\OM\subset\C^d$ be a Kobayashi hyperbolic domain.
\begin{itemize}
    \item[(i)] Let $\xi,\eta\in\bdy\clos{\OM}^\text{End}$ with $\xi\neq \eta$. The pair $\xi, \eta$ has {\em visibility property} with respect to $K_{\OM}$ if, given any $\l\geq 1$ and $\k\geq 0$, there exist neighborhoods $U_{\xi}, U_{\eta}$ of $\xi , \eta$ respectively in $\clos{\OM}^{\text{End}}$ with $\clos{U}^{\text{End}}_{\xi}\cap \clos{U}^{\text{End}}_{\eta} = \emptyset$ and a compact set $K\subset\OM$ such that for every $\lk$-almost-geodesic $\gamma:[a,b] \to \OM$ with $\gamma(a) \in U_{\xi}$ and $\gamma(b) \in U_{\eta}$, $\gamma\cap K \neq \emptyset$.
    \item[(ii)] We say $\OM$ is a {\em visibility domain} (or satisfies {\em visibility property}) with respect to $K_{\OM}$ if every distinct pair $\xi , \eta \in \bdy\clos{\OM}^{\text{End}}$ has visibility property with respect to $K_{\OM}$.
    \item[(iii)] Let $\xi,\eta\in\bdy\clos{\OM}^\text{End}$ with $\xi\neq \eta$. The pair $\{\xi,\eta\}$ has {\em weak visibility property} with respect to $K_\OM$ if the condition in (\text{i}) is satisfied only for $\l =1$.
    \item[(iv)] We say $\OM$ is a {\em weak visibility domain} (or satisfies {\em weak visibility property}) with respect to $K_{\OM}$ if every distinct pair $\xi , \eta \in \bdy\clos{\OM}^{\text{End}}$ has weak visibility property with respect to $K_{\OM}$.
\end{itemize}
\end{definition}

\begin{remark}\label{Rem:Comment_on_End-space}
If $\OM$ is weak visibility domain, then $\clos{\OM}^\text{End}$ is sequentially compact (see for instance \cite[Theorem~3.3]{CGMS2024}). Moreover, if every distinct pair of Euclidean boundary points of domain $\OM$ has (weak) visibility property, then every distinct pair of ends of $\clos{\OM}$ have (weak) visibility property (see \cite{Masanta, CGMS2024}).
\end{remark}

The (weak) visibility property with respect to the Kobayashi distance has many applications in several complex variables, such as extension and boundary behavior of certain holomorphic mappings and as well as quasi-isometries with respect to Kobayashi distance (see \cite{BZ2017,Bharali_Maitra,BNT2022,BZ2023,CMS2024,CGMS2024}).
\smallskip

One may ask the following:
\begin{question}\label{Q:intrinsic_visibility}
Is there a notion of visibility with respect to the Kobayashi distance that is
intrinsic to the domain, i.e., independent of its boundary, and that allows one
to classify a natural class of Denjoy-Wolff domains in $\C$?
\end{question}
We define a notion of visibility property with respect to Kobayashi distance which is independent of the boundary points. We refer to this notion of visibility as intrinsic visibility.
\begin{definition}
 Let $\OM \subset \C^d$ be a Kobayashi hyperbolic domain and let $A, B \subset \OM$.
 \begin{itemize}
\item[(i)] We say {\em $A$ and $B$ have visibility property} with respect to $(\OM , K_\OM)$ if, given any $\l\geq 1$ and $\k\geq 0$, there exists a compact set $K \subset \OM$ such that every $\lk$-almost-geodesic $\gamma : [a,b] \rightarrow \OM$  with respect to $(\OM, K_\OM)$ with $\gamma(a) \in A$ and $\gamma(b) \in B, \gamma([a,b]) \cap K \neq \emptyset$.
\item[(ii)] We say {\em $A$ and $B$ have weak visibility property} with respect to $(\OM,K_\OM)$ if the condition in (\text{i}) is satisfied only for $\l=1$.
 \end{itemize}   
\end{definition}

In Section~\ref{Sec:Localization_results}, Lemma~\ref{Lem:Vis_of_BZ_same_to_our_vis} provides a relation between the intrinsic visibility property and standard visibility property. 
\subsection{The localization of Kobayashi distance}
The second named author in \cite{ADS2023} used the standard visibility to localize the Kobayashi distance, such results are commonly referred to as localization results. A notable feature of his results is that they require only local weak visibility assumption, rather than additional stronger geometric hypotheses such as local convexity or local convexifiability. We now recall the localization results established by the second named author.
\begin{result}[{\cite[Theorem~1.1]{ADS2023}}]\label{Res:Loc_loclization_res}
Suppose $\OM \subset  \C^d$ is a Kobayashi hyperbolic domain and $U$ is a subdomain of $\OM$ such that $\bdy U \cap \partial \Omega \neq \emptyset$. Suppose $(U, K_U)$ has weak visibility property for every pair of distinct points in $ \bdy U \cap \partial \Omega$. Then for every open set $W,W_1$ with $W \subset \subset W_1\subset\subset U$ and $K_\OM(W_1, \OM\setminus U)>0$, there exists a constant $C> 0$ which depends only on $U, W_1, W$ such that for every $z, w \in W$,
\[
K_U (z, w) \leq K_{\Omega} (z, w) + C.
\]
\end{result}

\begin{result}[{\cite[Theorem~1.3]{ADS2023}}]\label{Res:Global_localization_result}
Suppose $\OM \subset  \C^d$ is a Kobayashi hyperbolic domain and $U$ is a subdomain of $\OM$ such that $\bdy U \cap \partial \Omega \neq \emptyset$. Suppose $(\OM, K_\OM)$ has weak visibility property for every pair of distinct points in $ \bdy U \cap \partial \Omega$. Then for every open set $W$ such that $W \subset \subset U$, there exists a constant $C> 0$ which depends only on $U, W$ such that for every $z, w \in W$,
\[
K_U (z, w) \leq K_{\Omega} (z, w) + C.
\]
\end{result}

One of the main objectives of this article is to extend results quoted above to Kobayashi hyperbolic domains where standard visibility may fail to hold. The following theorem requires intrinsic visibility with respect to subdomain and, in additional, separation property which are analogues to the assumptions of Result~\ref{Res:Loc_loclization_res} in our setting.

\begin{theorem}\label[theorem]{Thm:Local_Loc_Res}
Suppose $\OM \subset \C^d$ is a Kobayashi hyperbolic domain and $U \subset \OM$ is a sub-domain. Then for every set $W,W_1,W_2$ with $W_1$ and $W_2$ are open sets such that
\begin{itemize}
\item[(i)] $W \subset W_1 \subset W_2 \subset U$ and $K_\OM(W_2, \OM \setminus U) > 0$, and
\item[(ii)] $W $ and $U \setminus W_1$, and $W_1$ and $U \setminus W_2$ have weak visibility property with respect to $(U, K_U)$,
\end{itemize}
there exists a constant $C > 0$ which depends only on $U,W_2,W_1$ and $W$ such that for all $z , w \in W$, 
\[
K_U(z,w) \leq K_\OM(z,w ) + C.
\]
\end{theorem}

\begin{figure}[htbp]
\centering
\begin{minipage}[t]{0.48\textwidth}
  \centering
  \resizebox{0.70\linewidth}{!}{%
  \begin{tikzpicture}[region diagram]
    \domainoutline
    \commonregions
    \draw (middle-left)
      .. controls (371,404) and (641,472) .. (middle-right);
    \draw (lower-left)
      .. controls (143,560) and (204,585) .. (248,615)
      .. controls (310,663) and (373,623) .. (428,609)
      .. controls (547,578) and (666,647) .. (lower-right);
    \node at (876,318) {$\Omega$};
    \node at (643,553) {$W_2$};
    \node at (578,672) {$W_1$};
  \end{tikzpicture}%
  }
\end{minipage}\hfill
\begin{minipage}[t]{0.48\textwidth}
  \centering
  \resizebox{0.68\linewidth}{!}{%
  \begin{tikzpicture}[region diagram]
    \domainoutline
    \commonregions
    \coordinate (global-left) at (82,577);
    \draw (global-left)
      .. controls (145,494) and (207,507) .. (253,541)
      .. controls (292,570) and (310,609) .. (375,620)
      .. controls (464,638) and (548,591) .. (624,598)
      .. controls (697,603) and (747,644) .. (global-right);
    \node at (876,318) {$\Omega$};
    \node at (605,643) {$W_0$};
  \end{tikzpicture}%
  }
\end{minipage}
\caption{The first and second figure demonstrate the number of auxiliary sets required in Theorem~\ref{Thm:Local_Loc_Res} and Theorem~\ref{Thm:Global_Loc_Res} respectively.}
\end{figure}
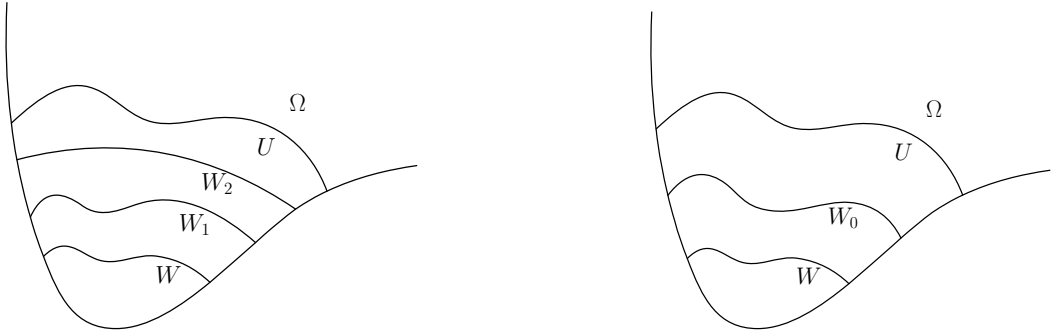
The next theorem assumes only intrinsic visibility with respect to ambient domain which is analogous to the assumptions of Result~\ref{Res:Global_localization_result} in our setting.
\begin{theorem}\label[theorem]{Thm:Global_Loc_Res}
Suppose $\OM \subset \C^d$ is a Kobayashi hyperbolic domain and $U \subset \OM$ is a sub-domain. Then for every set $W, W_1$ with $W_1$ is an open set such that
\begin{itemize}
\item[(i)] $W \subset W_1 \subset U$ and $\overline{W_1} \cap \OM \subset U$ and
\item[(ii)] $W$ and $\OM\setminus W_1$, and $W_1$ and $\OM\setminus U$ have weak visibility property with respect to $(\OM, K_\OM)$,
\end{itemize} 
there exists a constant $C>0$ which depends only on $U,W_1$ and $W$ such that for every $z, w \in W$,
\[
K_{U}(z,w) \leq K_\OM(z,w) + C.
\]
\end{theorem}
\begin{remark}
In Theorem~\ref{Thm:Local_Loc_Res} and Theorem~\ref{Thm:Global_Loc_Res}, intrinsic visibility allows that boundary of the subdomain need not intersect the boundary of the ambient domain. This is not the case in Result~\ref{Res:Loc_loclization_res} and Result~\ref{Res:Global_localization_result}.
\end{remark}

We provide an example of planar domain where Result~\ref{Res:Global_localization_result} is not applicable while Theorem~\ref{Thm:Global_Loc_Res} is applicable see Example~\ref{Ex:Our_loc_holds_while_previous_one_doesnot} in Section~\ref{Sec:Examples}. For the proofs of these results we refer to Section~\ref{Sec:Localization_results}. However, the proofs of both theorems are almost similar to the results quoted above i.e., localization results under the standard visibility.

\subsection{Comparison between Busemann horospheres}
\begin{definition}
Let $\OM\subset\C^d$ be a complete kobayashi hyperbolic domain. Given $R>0$ and a geodesic ray $\gamma: [0, \infty) \to \OM$ of $K_{\OM}$, consider the {\em Busemann horosphere} of $\gamma$ as
\[
H_{\OM}(\gamma,R)\defeq \Big\{z \in \OM : \lim_{t \rightarrow \infty}[K_{\OM}(z, \gamma(t)) - K_{\OM}(\gamma(0), \gamma(t))] < \frac{1}{2}\log R\Big\}.
\]
\end{definition}

In this article, we study the relation between Busemann horospheres of two geodesics if one is geodesic of subdomain and another is geodesic of the ambient Kobayashi hyperbolic domain. Their proofs are given in Section~\ref{Sec:Busemann_horo}. We now state such results.

\begin{theorem}\label[theorem]{Thm:Geo_horosphere_of_small_domain_is_subset_of_geo_horosphere_of_big_domain}
Suppose $\OM \subset \C^d$ is complete hyperbolic domain and $U \subset \OM$ is a subdomain such that $(U,K_{U})$ is complete. Suppose 
\begin{enumerate}
\item[(i)] $W, W_1 \subset \OM$ such that $W\subset W_1 \subset U, W_1$ is an open set and $\overline{W_1} \cap \OM \subset U$, and
\item[(ii)] $W$ and $\OM \setminus W_1$, and $W_1$ and $\OM\setminus U$ have weak visibility property with respect to $(\OM, K_\OM)$.
\end{enumerate}
Suppose $\gamma : [0, \infty) \rightarrow \OM$ and $\wt{\gamma}:[0,\infty) \to U$ are geodesic rays of $(\OM, K_\OM)$ and $(U,K_{U})$, respectively, with $\gamma, \wt{\gamma} \subset W$ and there exists a constant $C^{\prime}>0$ such that for all $t \in[0,\infty)$
\begin{equation}\label{ineq:Asymp_Condition_for_quasi_geo_1}
        K_\OM(\gamma(t), \wt{\gamma}(t)) \leq C^{\prime}.
\end{equation}
Then there exists a constant $C> 0$ such that for all $R>0$
\[
H_{U}(\wt{\gamma}, R) \subseteq H_{\OM}(\gamma, CR).
\]
\end{theorem}

\begin{theorem}\label[theorem]{Thm:Geo_horosphere_of_big_domain_is_subset_of_geod_of_small_domain_upto_some_compact_set_and-constant}
Suppose $\OM \subset \C^d$ is complete hyperbolic domain and $U \subset \OM$ is a sub-domain such that $(U,K_{U})$ is complete. Suppose 
\begin{enumerate}
\item[(i)] $W, W_1, W_2 \subset \OM$ such that $W \subset W_1\subset W_2 \subset U, W_2$ is an open set and $\overline{W_2} \cap \OM \subset U$, and
\item[(ii)] $W$ and $\OM \setminus W_1$, $W_1$ and $\OM \setminus W_2$, and $W_2$ and $\OM\setminus U$ have weak visibility property with respect to $(\OM, K_\OM)$.
\end{enumerate}
Suppose $\gamma : [0, \infty) \rightarrow \OM$ and $\wt{\gamma}:[0,\infty) \to U$ are geodesic rays of $(\OM, K_\OM)$ and $(U,K_{U})$, respectively, with $\gamma, \wt{\gamma} \subset W$ and there exists a constant $C^{\prime}>0$ such that for all $t \in[0,\infty)$
\begin{equation}\label{ineq: Asymp_Condition_for_quasi_geo_2}
K_\OM(\gamma(t), \wt{\gamma}(t)) \leq C^{\prime}.
\end{equation}
Then there exist constant $C > 0$ and compact set $K \subset \OM$ such that
\[
H_{\OM}(\gamma, R) \setminus K \subseteq H_{U}(\wt{\gamma}, CR)
 \]
where $K$ depend on $R$, geodesic ray $\gamma$ and $W_1$.
\end{theorem}

In Section~\ref{Sec:Busemann_horo}, we also study the relation between Busemann horospheres of two geodesics of the same complete hyperbolic domain see Lemma~\ref{Lem:Asym_compares_two_geo_horo} and Lemma~\ref{Lem:Strong_asym_imp_geo_horo_are_same}.

\subsection{The Denjoy-Wolff theorem}
The sequence of iterates of holomorphic self maps of Koabyashi hyperbolic domains usually have very simple behavior. In the case of the unit disc, it is the celebrated theorem due, independently, to Denjoy and Wolff which states that 
\begin{result}[{\cite{Denjoy,Wolff1,Wolff2}}]\label{Res:DWinD}
Let $f: \D \rightarrow \D$ be a holomorphic map. Then either
\begin{enumerate}
\item $f$ has a unique fixed in $\D$, or
\item there exists a point $ \xi \in \partial{\D}$ such that the sequence of iterates $(f^{n})_{n \geq 1}$ converges to $\xi$ uniformly on compact subsets of $\D$.
\end{enumerate}
\end{result}

If the domain is not contractible (e.g., annulus in the plane), then in general the dichotomy in Result~\ref{Res:DWinD} fails. However, it was Heins who proved the Denjoy-Wolff type theorem for hyperbolic Riemann surface (see for instance \cite{AbateNewbook}). The Denjoy-Wollf theorem was proved in the case of the unit ball in $\C^{d}$ by Herv\'e \cite{Herve}. Abate generalized the result of Herv\'e to strongly convex domain with $C^{2}$-smooth boundaries \cite{Abate:Horo_and_itera,AbateBook}. Abate introduced an appropriate dichotomy which also works for domains in higher dimensions and generalized Result~\ref{Res:DWinD} for strongly pseuodoconvex domain with $C^{2}$-smooth boundaries \cite{Abate:Ite_theory_and_comp_div}. The latter result was also extended to strictly convex domains with out the boundary regularity assumptions (see \cite{AR2014} and the references therein). The Denjoy-Wolff theorem has also been extended for certain type of metric spaces where a boundary at inifinity replaces the topological boundary (see for instance \cite{Beardon1,Beardon2,Karlsson}), and even to certain holomorphic self maps of open unit balls of reflexive Banach spaces that are reasonably ``nice'' (see \cite{BKR} and the references therein). Moreover, Bracci-Gaussier-Zimmer generalized the Denjoy-Wolff theorem for bounded convex domains which are Gromov hyperbolic with respect to the Kobayashi distance \cite{Bracci_Gaussier_Zimmer}.

One important application of the standard visibility property is its role in extending the Denjoy–Wolff theorem to broader classes of domains. Bharali and Maitra extended the Denjoy-Wolff theorem for bounded domain $\OM\subset \C^d$ having the standard visibility property under the additional mild assumption that $\OM$ is taut \cite{Bharali_Maitra}. The latter result was subsequently extended by Chandel-Maitra-Sarkar to bounded domain $\OM\subset \C^d$ satisfying weak visibility property \cite{CMS2024}. Futhermore, Bharali and Zimmer extended the Denjoy-Wolff theorem for unbounded weak visibility domain \cite{BZ2023}. Recently, Masanta studied the standard visibility property for domains in complex manifolds and established a corresponding Denjoy–Wolff theorem in this setting \cite{Masanta}.

\begin{remark}
In general, for Kobayashi hyperbolic domain $\OM \subset \C^{d}, d\geq 2$, it is not known whether the Denjoy-Wolff property is weaker than the standard visibility property. In fact, Bracci and {\"O}kten conjectured that for bounded convex domain the standard visibility property is equivalent to the Denjoy-Wolff property \cite{BO_Failure,BB1}. This is an open problem, even for domains with smooth boundaries.
\end{remark}

Recently, Benini and Bracci characterized simply connected domain $\OM\subsetneq \C$ that has Denjoy-Wolff property in terms of Busemann horosphere.
\begin{definition}\label{Def:DWP}
Let $\OM\subset \C$ be a hyperbolic domain. We say that $\OM$ has the {\em Denjoy-Wolff} property if for every $f \in \hbox{Hol}(\OM, \OM)$ either for any $z\in\OM$ the orbit $\{f^{n}(z):n\geq1\}$ is relatively compact in $\OM$ or there exists $\xi\in\bdy_{\infty}\OM$ such that the sequence of iterates $(f^n)_{n\geq 1}$ converges uniformly on compact subset of $\OM$ to $\xi$.
\end{definition}

\begin{definition}\label{Def:H-limit}
Let $f: \D \rightarrow \C$ be a function. We say that $L \in \C_{\infty}$ is the {\em Horosphere limit} (in short {\em H-limit}) of $f$ at $\xi \in \partial{\D}$ if $f(z) \to L$ as $z$ tends to $\xi$ within $H_\D(\sigma, R)$ for all $R>0$, where $\sigma$ is a radial line that joins $0$ and $\xi$.
\end{definition}

With these definitions in place, we state Benini and Bracci's result below. We refer the reader to Section~\ref{Sec:Notations} for the notations used in the statement.
\begin{result}[{Paraphrasing Theorem~3.1 in \cite{Benini-Bracci}}]\label{Thm:Main_Thm_of_BB}
Let $\Omega\subsetneq \C$ be a simply connected domain. The following are equivalent:
\begin{enumerate}
\item $\Omega$ has the Denjoy-Wolff Property;
\item for every geodesic ray $\gamma:[0,\infty)\to\OM$ of $(\OM, K_{\OM})$, there exists $x\in\bdy_{\infty}\OM$ such that $\overline{H_\Omega(\gamma, R)}^{\infty}\cap \bdy_{\infty}\OM=\{x\}$ for every $R>0$;
\item there exist $p\in\Omega$ such that for every geodesic ray $\gamma:[0,\infty)\to\OM$ of $(\OM, K_{\OM})$ with $\gamma(0)=p$, there exists $x\in\bdy_{\infty}\OM$ and $R>0$ such that $\overline{H_\Omega(\gamma, R)}^{\infty}\cap \bdy_{\infty}\OM=\{x\}$;
\item for every Riemann map $h:\D\to \Omega$, $h$ has H-limit at every $\xi\in\partial\D$.
\end{enumerate}
\end{result}

Before presenting the definition of Condition ($\star$) and Condition (M$\star$), we note the following definition.

\begin{definition}
Let $U\subsetneq \D$ a subdomain such that $\bdy U\cap\bdy\D$ contains $1$ as an interior point. We say $U$ is a {\em bi-geodesic region} in $\D$ if $\bdy U\cap\D$ is a bi-geodesic ray in $\D$ for the Poincar$\e$ distance on the unit disc $\D$.
\end{definition}

\begin{definition}\label{Def:Condtion(*)}
Let $\OM \subset \C$ be a hyperbolic domain and let $x\in\bdy\OM$. We say $\OM$ satisfies Condition~($\star$) at $x$ if for every sequence $(z_{n})_{n \geq 1} \subset \OM$ with $z_{n} \rightarrow x$, there exist an injective holomorphic map $\varphi : \D \rightarrow \OM$, a bi-geodesic region $U \subset \D$ and a subsequence $(z_{n_{k}})_{k \geq 1}$ such that 
\begin{enumerate}
\item $x \in \partial{\varphi(\D)}\cap\partial{\OM}$; 
\item $(z_{n_{k}})_{k \geq 1} \subset \varphi(U), \varphi^{-1}(z_{n_{k}}) \rightarrow 1$ as $k \rightarrow \infty$; 
\item for every sequence $(\zeta_{n})_{n \geq 1} \subset U$ with $\delta_{\D}(\zeta_{n}) \rightarrow 0$ as $n \rightarrow \infty$, either $\delta_{\OM}(\varphi(\zeta_{n})) \rightarrow 0$ or $|\varphi(\zeta_{n})| \rightarrow \infty$ as $n \rightarrow \infty$;
\item for every bi-geodesic regions $U_{1}$ and $U_{2}$ with $U_{1} \subsetneq U_{2} \subseteq U$, we have $K_\OM(\varphi(U_{1}), \varphi(\D) \setminus \varphi(U_2)) > 0$.
\end{enumerate}
We say that $\OM$ satisfies {\em Condition~($\star$)} if it satisfies Condition~($\star$) at every boundary point of $\OM$.
\end{definition}
\begin{remark}
The condition $(3)$ in the definition above corresponds to the properness of map $\varphi$ on the bi-geodesic region $U$.
\end{remark}

We now introduce the mild version of the Condition~($\star$) as follows.
\begin{definition}\label{Def:Condtion(M*)}
Let $\OM \subset \C$ be a hyperbolic domain and let $x\in\bdy\OM$. We say $\OM$ satisfies Condition~(M$\star$) at $x$ if for every sequence $(z_{n})_{n \geq 1} \subset \OM$ with $z_{n} \rightarrow x$, there exist an injective holomorphic map $\varphi : \D \rightarrow \OM$, a bi-geodesic region $U \subset \D$ and a subsequence $(z_{n_{k}})_{k \geq 1}$ such that 
\begin{enumerate}
\item $x \in \partial{\varphi(\D)}\cap\partial{\OM}$;
\item for all $k \in \N, z_{n_k} \in \varphi(U)$ and $\varphi^{-1}(z_{n_{k}}) \rightarrow 1$ as $k \rightarrow \infty$;
\item for every sequence $(\zeta_{n})_{n \geq 1} \subset U$  with $\delta_{\D}(\zeta_{n}) \rightarrow 0$ as $n \rightarrow \infty$, either $\delta_{\OM}(\varphi(\zeta_{n})) \rightarrow 0$ or $|\varphi(\zeta_{n})| \rightarrow \infty$ as $n \rightarrow \infty$.
\end{enumerate}
We say that $\OM$ satisfies {\em Condition~(M$\star$)} if it satisfies Condition~(M$\star$) at every boundary point of $\OM$.
\end{definition}
Condition~($\star$) consists of Condition~(M$\star$) together with some additional requirement i.e., the condition $(4)$ in the definition of Condition ($\star$). Therefore, it is clear that Condition~($\star$) implies Condition~(M$\star$). In fact, the converse also holds; see Proposition~\ref{Pro:Cond(M*)_implies_Cond(*)}.

We now state one of the main results of this article which is a generalization of Result~\ref{Thm:Main_Thm_of_BB}.
\begin{theorem}\label[theorem]{Thm:Characterization_of_Condition(*)}
Suppose $\OM \subset \C$ is a hyperbolic domain that satisfies Condition~(M$\star$). The following are equivalent:
\begin{enumerate}
\item $\OM$ has the Denjoy-Wolff Property;
\item for every geodesic ray $\gamma :[0,\infty)\to \OM$ of $(\Omega,K_\OM)$ either $\overline{H_\Omega(\gamma, R)}\cap \partial \Omega=\emptyset$ for every $R>0$ or there exists $x\in\partial\Omega$  such that $\overline{H_\Omega(\gamma, R)}\cap \partial \Omega=\{x\}$ for  every $R > 0$;
\item there exist $z_0\in\Omega$ such that for every geodesic ray $\gamma :[0,\infty)\to \OM$ of $(\Omega,K_\OM)$ with $\gamma(0)=z_0$ either there exists $R>0$ such that $\overline{H_\Omega(\gamma, R)}\cap \partial \Omega=\emptyset$ or there exist $x\in\partial\Omega$ and $R>0$ such that $\overline{H_\Omega(\gamma, R)}\cap \partial\Omega=\{x\}$;
\item for every sequence $(z_{n})_{n \geq 1}\subset \OM$ that converges to some boundary point $x\in \bdy\OM$ such that if $\varphi :\D\to \OM$ is an injective holomorphic map obtained by $x$ and $(z_{n})_{n \geq 1}$ as in Condition~(M$\star$), then either {\sl $H$-limit} of $\varphi$ at $1$ exists and it is finite or it is infinity.
\end{enumerate}
\end{theorem}

The definition of Condition~(M$\star$) may seem very technical. However, it is not technical as it first appears. There is an alternate way to describe Condition~(M$\star$) which we call Condition $(\star\star)$. One of the advantages to have an alternate description of Condition~(M$\star$) is to produce examples of such domains. Before introducing such definition, let us set $A(0,r,1)\defeq\{z\in \C:r<|z|<1\}$. We now introduce Condition $(\star\star)$.
\begin{definition}
Let $\OM \subset \C$ be a hyperbolic domain. We say that $\OM$ satisfies {\em Condition $(\star\star)$} if for every connected component $\Gamma$ of $\bdy\OM$, there exists an injective holomorphic map $f:A(0,r,1)\to\OM$ such that for every sequence $(\zeta_n)_{n\geq 1}\in A(0,r,1)$ with $|\zeta_n|\to 1$ as $n\to\infty$, we have either $\delta_\OM(f(\zeta_n))\to 0$ as $n\to\infty$ or  $|f(\zeta_n)|\to \infty$ as $n\to\infty$.
\end{definition}

Geometrically, Condition $(\star\star)$ can be thought of domains with attached annulus at its every boundary component. We prove that Condition $(\star\star)$ implies Condition~(M$\star$); see Lemma~\ref{Lem:Contions_double_star_implies_mild_star}. Therefore, following two corollaries are direct consequences of Theorem~\ref{Thm:Characterization_of_Condition(*)}. We refer Section~\ref{Sec:Prop_of_Con(*)_and_relation_bw_Cond(*)_and_Cond(**)} for more details about the relation between Condition~(M$\star$) and Condition $(\star\star)$.

\begin{corollary}\label[corollary]{Cor:Finitely_conn_dom_charac_for_DWP}
Suppose $\OM\subset \C$ is a finitely connected domain. Then the conditions $(1)-(4)$ in Theorem~\ref{Thm:Characterization_of_Condition(*)} are equivalent.
\end{corollary}

\begin{corollary}\label[corollary]{Cor:Well_seperated_domains_charac}
Suppose $\OM \subset \C$ is an (possibly infinitely connected) domain such that every connetced component of $\partial \OM$ is bounded, non-degenerate and there exists a positive constant that is a lower bound for the Euclidean distance between any two of connected components of $\partial\OM$. Then the conditions $(1)-(4)$ in Theorem~\ref{Thm:Characterization_of_Condition(*)} are equivalent.
\end{corollary}

In fact, we also extend Theorem~\ref{Thm:Characterization_of_Condition(*)} for planar domain which satisfies Condition~(M$\star$) outside the totally disconnected set. We state it in the following way.
\begin{theorem}\label[theorem]{Thm:Main_thm_outside_totally_disc_set}
Let $\OM \subset \C$ be a hyperbolic domain. Suppose that there exists a totally disconnected set $S \subset \bdy\OM$ such that $\OM$ satisfies Condition $(M \star)$ at every $\xi \in \partial\OM\setminus S$. Then the following are equivalent:
\begin{itemize}
\item[(1)] $\OM$ has the Denjoy-Wolff Property;
\item[(2)] for every geodesic ray $\gamma:[0, \infty) \to \OM$ of $(\OM, K_\OM)$, either $\clos{H_{\OM}(\gamma, R)}\cap \bdy\OM =\emptyset$ for all $R>0$ or there exists $x \in \bdy \OM$ such that for all $R>0$
\[
\clos{H_{\OM}(\gamma, R)}\cap \bdy\OM = \{x\};
\]
\item[(3)] for every sequence $(z_{n})_{n \geq 1}$ in $\OM$ that converges to some boundary point $y \in \bdy\OM\setminus S$ such that if $(\gamma_{k})_{k\geq 1}$ is a sequence of geodesic $(\OM,K_\OM)$ that joins some $p \in \OM$ and $z_{n_{k}}$ such that $(\gamma_{k})_{k\geq 1}$, up to a subsequence, converges uniformly on compact set to a geodesic ray $\gamma$ of $(\OM, K_\OM)$; where $\varphi:\D \to \OM$ is an injective holomorphic map and $(z_{n_k})_{k \geq 1}$ is a subsequence of $(z_{n})_{n \geq 1}$ as in Condition $(M \star)$ obtained by $(z_{n})_{n \geq 1}$ and $y$, then either $\clos{H_{\OM}(\gamma, R)}\cap \bdy\OM =\emptyset$ for all $R>0$ or there exists $x \in \bdy\OM$ such that for all $R>0$
\[
\clos{H_{\OM}(\gamma, R)}\cap \bdy\OM = \{x\}.
\]
\end{itemize}
\end{theorem}

\noindent The proof of the Theorem~\ref{Thm:Characterization_of_Condition(*)} and Theorem~\ref{Thm:Main_thm_outside_totally_disc_set} are presented in Section~\ref{Sec:Proofs_of_main_theorems} while Section~\ref{Sub:Embedded_disc_preserve_vis} and Section~\ref{Sub:Prop_of_M*} contain some auxiliary lemmas which are used their proofs.

\subsection{Notations}\label{Sec:Notations} We now fix some notations for this article.
\begin{itemize}
\item Let $\N$ denote the set of all natural numbers.
\item Let $\D\defeq\{z\in\C:|z|<1\}$ denote the open unit disc in the complex plane and $\rho_\D$ denote the Poincar$\e$ metric on $\D$.
\item Let $\K\defeq\{z\in \C:\rprt z>0\}$ denote the right half plane.
\item For $z \in \C^d$, let $\norm{z}$ denote the standard Euclidean norm of $z$.
\item For a point $p \in \C^d$ and $r>0$, $B(z,r)\defeq\{z\in\C^d:\norm{z-p}<r\}$.
\item Given a domain $\OM \subset \C^d$ and $z \in \OM$ let
 \[
\delta_\OM(z)\defeq\inf\{\norm{w-z}:w\in\bdy\OM\}.
\]
\item Given two domains $\OM_1\subseteq\C^{d_1}$ and $\OM_2\subseteq\C^{d_2}$, let $\hbox{Hol}(\OM_1,\OM_2)$ be the space of holomorphic maps from $\OM_1$ to $\OM_2$.
\item For a set $A\subset \C^d$ and $\epsilon>0$,
\[
N_{\C^d}(A,\epsilon)\defeq\{z\in\C^d:\inf_{w\in A}\|z-w\|<\epsilon\}.
\]
\item Let $\OM\subseteq\C^d$ be a domain. For a curve $\gamma:I\to\OM$ of an interval $I\subset\R$, by abuse of notation, we denote image of the curve $\gamma$ by $\gamma$ itself.
\item Let $\C_{\infty}$ be the one point compactification of $\C$, i.e. $\C_\infty\defeq\C\cup\{\infty\}$.
\item Given a domain $\OM\subseteq \C$, let $\clos{\OM}^{\infty}$ be the closure of $\OM$ in $\C_\infty$.
\item Given a domain $\OM\subseteq \C$, let $\bdy_{\infty}\OM$ denotes its boundary in $\C_\infty$. 

\end{itemize}

\section{Localization of the Kobayashi distance}\label{Sec:Localization_results}
We start this section by introducing Kobayashi metric and Kobayashi distance. The main aim of this section is to present proof of Theorem~\ref{Thm:Local_Loc_Res}, Theorem~\ref{Thm:Global_Loc_Res} and their consequences.

\subsection{The Kobayashi pseudometric and pseudodistance}\label{Subsec:Kob_dist} Given a domain $\OM\subseteq\C^d$ the Kobayashi pseudometric is defined by
\[
\k_\OM(z;v)\defeq \inf\{\rho_\D(\zeta; \xi): \exists \phi \in \hbox{Hol}(\D, \OM) \text{ such that } \phi(\zeta) = z \text{ and } \xi\phi^{\prime}(\zeta) = v\}
\]
where $z\in\OM$ and $v\in\C^d$. It follows from the result of Royden \cite[Proposition~3]{RoyKobDist} that the Kobayashi metric is upper semi-continuous. Thus, one can define the length of piecewise $C^{1}$ curve $\gamma:[a,b]\to\OM$ to be
\[
l_\OM^{\k}(\gamma)\defeq \int_{a}^{b}\k_\OM(\gamma(t);\gamma^\prime(t))dt.
\]
Now, one can define the Kobayashi pseudodistance to be
\[
K_\OM(z,w)\defeq\inf\{l_\OM^{\k}(\gamma): \gamma:[a,b] \to \OM \text{ is piecewise } C^{1} \text{ curve with } \gamma(a)=z \text{ and }\gamma(a)=w\}.
\]
A domain $\OM\subset\C^d$ is said to be {\em Kobayashi hyperbolic} if the Kobayashi pseudodistance $K_\OM$ is distance. Moreover, a domain $\OM$ is said to be {\em complete hyperbolic} if the space $(\OM,K_\OM)$ is Cauchy complete. A planar Kobayashi hyperbolic domain is commonly referred to as a {\em hyperbolic} domain. For a detailed discussion of the Kobayashi pseudometric and Kobayashi pseudodistance, we refer the reader to \cite{JarPfl}.

\begin{definition}
Let $\Omega \subset \C^{d}$ be a Kobayashi hyperbolic domain and let $\l\geq 1$ and $\k\geq 0$. A map $\gamma: I\to\Omega$, where $I\subset\R$ is an interval, is 
said to be a {\em $(\lambda,\kappa)$-almost-geodesic} if the following are satisfied:
\begin{itemize}
\item[(i)]for all $\,s,t\in I$
\[
\frac{1}{\lambda}|s-t|-\kappa\leq K_\OM(\gamma(s),\gamma(t))\leq \lambda |s-t|+\kappa \text{ and }
\]
\item[(ii)] $\gamma$ is absolutely continuous such that for 
almost every $t\in I$
\[
\k_\Omega(\gamma(t);\gamma'(t))\leq\lambda.
\] 
\end{itemize}
\end{definition}

The curve $\lk$-almost-geodesic plays an important role in the standard visibility property due to the following result.
\begin{result}[{\cite[Proposition~5.3]{BZ2023}}]
Let $\OM\subset\C^d$ be a Koabayshi hyperbolic domain. Then for any $\k>0$ and $z,w\in\OM$, there exists a $(1,\k)$-almost-geodesic $\gamma:[a,b]\to\OM$ such that $\gamma(a)=z$ and $\gamma(b)=w$.    
\end{result}

\begin{definition}
Let $\OM\subset\C^d$ be a Kobayashi hyperbolic domain. Let $\lambda\geq 1$, $\kappa>0$ and $I\subset\R$ be an interval. An absolutely continuous curve $\gamma:I\to\OM$ is called a $\lk${\em-geodesic} for $(\OM,K_\OM)$ if for all $s\leq t \in I$
\[
l_\OM^{\k}(\left.\gamma\right|_{[s,t]}) \leq \lambda K_\OM(\gamma(s), \gamma(t)) +\kappa.
\]
\end{definition}

\begin{remark}
Every $\lk$-almost-geodesic of $(\OM,K_\OM)$ is $(\lambda^2, \lambda^2\kappa)$-geodesic for $(\OM,K_\OM)$.
\end{remark}

We now would like to record some technical results about almost-geodesic, geodesic and about their reparameterization.

\begin{result}[{\cite[Proposition~2.3]{CGMS2024}}]\label{Res:Unit_speed_rep}
Suppose $\OM \subset \C^d$ is a Kobayashi hyperbolic domain and $\gamma:[a,b]\to\OM$ is an absolutely continuous map such that $\gamma^\prime$ is almost everywhere non-vanishing. Then there exists an absolutely continuous curve $\wt\gamma:[0,L]\to\OM$ such that it is a reparameterization of $\gamma$ and such that $\kappa_\OM(\gamma(t),\gamma^\prime(t)) = 1$ for a.e. $t\in [0,L]$.
\end{result}

\begin{result}[{\cite[Proposition~2.6]{CGMS2024}}]\label{Res:(lk)-geo_is_(lk)-al-geo}
Suppose $\OM \subset \C^d$ is a Kobayashi hyperbolic domain. Let $\lambda\geq 1$, $\kappa>0$ and let $\gamma:[a,b]\to \OM$ be a $\lk$-geodesic for $(\OM, K_\OM)$ such that $\gamma^\prime$ is almost everywhere non-vanishing. Let $\wt\gamma$ be the curve as in the Result~\ref{Res:Unit_speed_rep} above. Then $\wt\gamma$ is $\lk$-almost-geodesic for $(\OM,K_\OM)$.
\end{result}

\begin{result}[{\cite[Proposition~2.7]{CGMS2024}}]\label{Res:Fattenig_of_(lk)-geo}
Suppose $\OM \subset \C^d$ is a Kobayashi hyperbolic domain. Given $\lambda\geq 1$ and $\kappa>0$, let $\gamma:[a,b]\to\OM$ be a $\lk$-geodesic for $(\OM,K_\OM)$ and let $\epsilon>0$ be such that $N_{\C^d}(\gamma,\epsilon)\subset\OM$. Then there exists a absolutely continuous map $\wt\gamma:[a,b]\to\OM$ with $\wt\gamma(a) = \gamma(a)$ and $\wt\gamma(b) = \gamma(b)$ such that
\begin{itemize}
\item[(i)] $\wt\gamma$ is $(\lambda,\lambda\epsilon+\kappa)$-geodesic for $(\OM,K_\OM)$,
\item[(ii)] for every $t\in [a,b]$, $\|\gamma(t)-\wt\gamma(t)\|\leq \epsilon/2$ and
\item[(iii)] $\wt{\gamma}^\prime$ is almost everywhere non-vanishing.
\end{itemize}
\end{result}

We end this subsection by proving the following lemma which establishes the relation between standard visibility property and intrinsic visibility property.
\begin{lemma}\label{Lem:Vis_of_BZ_same_to_our_vis}
Suppose $\OM \subset \C^d$ is Kobayashi hyperbolic domain. Then the following are equivalent:
\begin{itemize}
\item[(1)] $\OM$ is (weak) visibility domain;
\item[(2)] For any two sets $A, B \subset \OM$ such that $\clos{A}^{End} \cap \clos{B}^{End} = \emptyset$, $A$ and $B$ have (weak) visibility property with respect to $(\OM, K_\OM)$.
\end{itemize}
\end{lemma}

\begin{proof}
First, we show that $(2) \implies (1)$. Let $p, q \in \bdy{\OM}$ with $p \neq q$. There exists $r>0$ such that $\clos{B(p,r)} \cap \clos{B(q,r)} = \emptyset$. We take $A := B(p,r) \cap \OM$ and $B := B(q,r) \cap \OM$. Since $A$ and $B$ are bounded sets, so we have $\clos{A}^{End} \cap \clos{B}^{End} = \emptyset$. Since $A$ and $B$ have (weak) visibility property with respect to $(\OM, K_\OM)$, the pair $\{p, q\}$ is visible pair with respect to $(\OM, K_\OM)$. Since $p,q \in \bdy{\OM}$ are arbitrary points, by Remark~\ref{Rem:Comment_on_End-space}, $\OM$ is (weak) visibility domain with respect to $(\OM, K_\OM)$.

We will show that $(1) \implies (2)$ by the method of contradiction. To see this assume that there exist $A, B \subset \OM$ such that $\clos{A}^{End} \cap \clos{B}^{End} = \emptyset$ and $A$ and $B$ does not have visibility property with respect to $(\OM, K_\OM)$. There exist $\kappa>0$ and a sequence $\gamma_n : [a_n,b_n] \rightarrow \OM$ of $\lk$-almost-geodesic of $(\OM, K_\OM)$ with $\gamma_n(a_n)\in A$ and $\gamma_n(b_n)\in B$ such that $(\gamma_n)_{n \geq 1}$ eventually avoids every compact set in $\OM$. Since $\clos{\OM}^{End}$ is sequentially compact, passing to a subsequence and relabeling, $\gamma_n(a_n) \rightarrow \xi_1 \in \bdy{\clos{\OM}}^{End}$, $\gamma_n(b_n) \rightarrow \xi_2 \in \bdy{\clos{\OM}}^{End}$ as $n \rightarrow \infty$ where $\xi_1 \neq \xi_2$. Since $\OM$ is (weak) visibility domain, there exist a neighborhoods $U_{\xi_1}$ and $U_{\xi_2}$ of $\xi_1$ and $\xi_2$, respectively, with $\clos{U}_{\xi_1}^\text{End}\cap \clos{U}_{\xi_2}^\text{End} = \emptyset$, a compact set $K\subset \OM$ and $N\in \N$ such that $\gamma_n\cap K\neq \emptyset$ for all $n \geq N$. This gives a contradiction to our starting assumption. This completes the proof.
\end{proof}

\subsection{Local localization result}
We start this subsection by stating the Royden's localization result of the Kobayashi metric and the lemma which is its consequence.
\begin{result}[{\cite[Lemma~2]{RoyKobDist}}]\label{Res:Royden_loc_result}
Suppose $\OM\subset \C^d$ is a Kobayashi hyperbolic domain and $U$ is a sub-domain. Then for all $z\in U$ and $v \in \C^d$,
\[
\k_U(z;v)\leq\coth(K_\OM(z,\OM\setminus U))\k_\OM(z;v).
\]
\end{result}
\begin{lemma}\label{Lem:Loc_of_metric}
Suppose $\OM \subset \C^d$ is a Kobayashi hyperbolic domain and $U \subset \OM$ is a sub-domain. Then for every set $W \subset U$ such that $K_\OM(W, \OM \setminus U) > 0$, there exists $L > 0$ such that for all $z \in W$ and $v \in \C^d$,
\[
\k_{U}(z;v) \leq \left(1 + Le^{-K_\OM(z, \OM \setminus U)}\right)\k_{\OM}(z;v).
\]
\end{lemma}

\begin{proof}
The proof follows exactly same as the proof of Lemma 3.1 in \cite{ADS2023}.
\end{proof}
The next two lemmas are analogous to Lemma~3.4 and Lemma~3.5 in \cite{ADS2023} which we prove in our setting.

\begin{lemma}\label{Fixed Point local Localization}
Suppose $\OM \subset \C^d$ is a Kobayashi hyperbolic domain and $U \subset \OM$ is a sub-domain. Then for every two sets $W, W_1$ with $W \subseteq W_1 \subset U$ such that $W_1$ is an open set and such that $W$ and $U\setminus W_1$ have weak visibility property with respect to $(U, K_U)$ and $o\in U$, there exists $L > 0$ such that for all $z \in W$
\[
K_U(z,o) \leq K_U(z, U \setminus W_1) + L.
\]
\end{lemma}

\begin{proof}
Let $z \in W$ and $w \in U \setminus W_1$ and $\kappa > 0$. Let $\gamma: [0,a] \rightarrow U$ be $(1,\kappa)$-almost-geodesic of $(U, K_U)$ such that $\gamma(0)=z$ and $\gamma
(a)=w$. Since $W$ and $U \setminus W_1$ have weak visibility property with respect to $(U, K_U)$, there exists a compact set $K \subset U$ such that $\gamma \cap K \neq \emptyset$. Let $\tilde{o} \in \gamma \cap K$ with $\gamma(\tilde{t}) = \tilde{o}$ for some $\tilde{t} \in (0,a)$. Then, using triangle inequality and the fact that $\gamma$ is $(1,\kappa)$-almost-geodeisc of $(U, K_U)$, we get following inequalities
\begin{align*}
K_U(z,o)&\leq K_U(z, \tilde{o}) + K_U(\tilde{o}, o)\\
        &\leq K_U(z, \tilde{o}) + \sup_{y \in K}K_U(y,o)\\
        &\leq \tilde{t} + \kappa + \sup_{y \in K}K_U(y,o)\\
        &\leq K_U(z, w) + 2\kappa + \sup_{y \in K}K_U(y,o).
\end{align*}
Since $w \in U \setminus W_1$ is arbitrary and the Kobayashi distance $K_\OM(\cdot, o)$ is continuous on compact set $K$, it follows that $K_{U}(z,o) \leq K_{U}(z, U \setminus W_0) + L$ where $L\defeq 2\kappa + \displaystyle \sup_{y \in K}K_\OM(y,o) < \infty$. This completes the proof of lemma.
\end{proof}

The proof of the following lemma is essentially the same as that of Lemma~3.5 in \cite{ADS2023}. Nevertheless, for completeness, we provide a proof in our setting, omitting some routine calculations.
\begin{lemma}\label{AG wrt U for Local}
Suppose $\OM \subset \C^d$ is a Kobayashi hyperbolic domain and $U \subset \OM$ is a sub-domain. Let $W, W_1$ be two sets such that $W \subseteq W_1 \subset U$, $W_1$ is an open, $K_\OM(W_1, \OM \setminus U) > 0$, and $W$ and $U \setminus W_1$ have weak visibility property with respect to $(U, K_U)$. Let $z, w \in W$ and $\gamma:[0,a]\to\OM$ be a $(1, \kappa)$-almost-geodesic of $(\OM, K_\OM)$ with $\gamma(0)=z$ and $\gamma(a)=w$ such that $\gamma \subset W$. Then there exists $\kappa_0 > 0$ such that $\gamma$ is $(1, \kappa_0)$-quasi-geodesic of $(U, K_U)$. 
\end{lemma}

\begin{proof}
Note that $K_\OM(z, \OM \setminus W_1) \leq K_\OM(z, \OM \setminus U)$ for all $z \in W$. Therefore, by Lemma \ref{Lem:Loc_of_metric}, there exists $L>0$ such that for all $z \in W_1$ and $v \in \C^{d}$, we have
\begin{equation}\label{Ineq:Loc_of_metric_under_separation}
\kappa_{U}(z;v) \leq \left(1 + Le^{-K_\OM(z, \OM \setminus W_1)}\right)\kappa_{\OM}(z;v).
\end{equation}
Note that $\gamma \subset W$. Using the fact that the Kobayashi distance is the integrated form of the Kobayashi metric, by~\eqref{Ineq:Loc_of_metric_under_separation} and the fact that $\gamma$ is $(1, \kappa)$-almost-geodesic of $(\OM, K_\OM)$,  we get that for every $s_1 \leq s_2 \in [0,a]$
\begin{align}
K_{U}(\gamma(s_1), \gamma(s_2)) &\leq \int_{s_1}^{s_2} \kappa_{U}(\gamma(t);\gamma^{\prime}(t))dt \notag\\
                                &\leq K_\OM(\gamma(s_1),\gamma(s_2)) + \k + L\int_{s_1}^{s_2} e^{-K_\OM(\gamma(t), \OM \setminus W_1)}dt \label{ineq: bound of Kob dist in terms of integral}.
\end{align}
\noindent {\bf Claim.} There exist $C_0, C_1 > 0$ such that for all $z \in W$
\begin{align*}
K_{U}(z, U \setminus W_1) \leq C_0 K_\OM(z, \OM \setminus W_1) + C_1.
\end{align*}
\noindent {\itshape Proof of Claim.} Its proof is exactly the same as that of the claim in the proof of Lemma~3.5 in \cite{ADS2023}. Therefore, we omit the proof.
\hfill $\blacktriangleleft$
\smallskip
    
Now, for $o \in U$, by Lemma \ref{Fixed Point local Localization}, there exists $L_0 > 0$ such that for all $z \in W$, we have 
\begin{align}
K_{U}(z, o) \leq K_{U}(z, U\setminus W_1) + L_0 \label{ineq: use of previous lem}. 
\end{align} 
Since $[0,a]$ is compact  and $K_{U}(\gamma(t), U \setminus W_1)$ is continuous function in the variable $t$, so it follows that there exists $t_0 \in [0, a]$ such that for all $t \in [0, a]$
\begin{align}
K_{U}(\gamma(t_0), U \setminus W_1) \leq K_{U}(\gamma(t), U \setminus W_1) \label{ineq: min attain for kobayashi distance}.
\end{align}
Now, by triangle inequality and using the fact that $\gamma$ is $(1,\kappa)$-almost-geodesic of $(\OM,K_\OM)$, and by \eqref{ineq: use of previous lem} and \eqref{ineq: min attain for kobayashi distance}, and by the previous claim, we get that for all $t \in [0,a]$
\begin{align*}
|t - t_0| - \kappa&\leq K_\OM(\gamma(t), \gamma(t_0))\\
                & \leq K_{U}(\gamma(t), o) + K_{U}(o, \gamma(t_0))\\
                & \leq K_{U}(\gamma(t), U\setminus W_1) + K_{U}(\gamma(t_0), U\setminus W_1) + 2L_0 \\
                & \leq 2K_{U}(\gamma(t), U\setminus W_1) + 2L_0\\
                & \leq 2C_0 K_\OM(\gamma(t), \OM \setminus W_1) + 2L_0 + 2C_1.
\end{align*}
Thus, for all $t \in [0,a]$
\begin{equation}\label{ineq: U-bound for the kob dist-1}
- K_\OM(\gamma(t), \OM \setminus W_1)  \leq \frac{-|t - t_0|}{2C_0} + \frac{2L_0 + 2C_1 + \kappa}{2C_0} 
\end{equation}
Using \eqref{ineq: bound of Kob dist in terms of integral} and \eqref{ineq: U-bound for the kob dist-1}, we get for all $s_1, s_2 \in [0,a]$,
\begin{align*}
K_{U}(\gamma(s_1), \gamma(s_2)) & \leq K_\OM(\gamma(s_1),\gamma(s_2)) + \k + L e^{\frac{2L_0 + 2C_1 + \kappa}{2C_0}} \int_{0}^{+ \infty} e^{\frac{-|t - t_0|}{2C_0}} dt\\
                                & \leq K_\OM(\gamma(s_1),\gamma(s_2)) + \k + 4LC_{2} C_{0}
\end{align*}
where $C_2 \defeq e^{\frac{2L_0 + 2C_1 + \kappa}{2C_0}}$ and the last inequality follows due to the following upper bound:
\begin{align}
\int_{0}^{+ \infty} e^{\frac{-|t - t_0|}{2C_0}} dt \leq 4C_0.
\end{align}
Since for all $s_1, s_2 \in [0,a]$, we have $K_\OM(\gamma(s_1), \gamma(s_2)) \leq K_{U}(\gamma(s_1), \gamma(s_2))$, and by taking $\kappa_0 \defeq \kappa + 4LC_{0} C_{2}$, $\gamma$ is $(1, \kappa_0)$-quasi-geodesic of $(U, K_{U})$. This completes the proof.
\end{proof}
We are now ready to present the proof of Theorem~\ref{Thm:Local_Loc_Res}. Although the proof is essentially the same as that of Theorem~1.1 in \cite{ADS2023}, we present it explicitly in our setting. We note that the proof in our setting is slightly simpler, as there is no need to care the boundary points.
\begin{proof}[Proof of {\Cref{Thm:Local_Loc_Res}}]
Let $z, w \in W$ and let $\kappa > 0$. Suppose $\gamma : [0, a] \rightarrow \OM$ is $(1,\kappa)$-almost-geodesic of $(\OM, K_\OM)$ with $\gamma(0)=z$ and $\gamma(a)=w$.
\smallskip
       
\noindent {\bf Case 1.} $\gamma \subset W_1$.
\smallskip

It is given that $K_\OM(W_2, \OM\setminus U) > 0$, $W_1$ and $U \setminus W_2$ have weak visibility property with respect to $(U, K_U)$, and $\gamma \subset W_1$. By replacing $W$ and $W_1$ with $W_1$ and $W_2$, respectively, in Lemma \ref{AG wrt U for Local}, there exists $\kappa_0 > 0$ such that $\gamma$ is a $(1,\kappa_0)$-quasi-geodesic of $(U, K_{U})$. Then, using the fact that $\gamma$ is $(1,\kappa)$-almost-geodesic of $(\OM, K_\OM)$, we get
\begin{align*}
K_U(z,w) &\leq a + \kappa_0\\
        &\leq K_\OM(z,w) + \kappa + \kappa_0.
\end{align*}
Take $C_1 \defeq \kappa + \kappa_0$. Therefore, in this case, we are done.
\smallskip

\noindent {\bf Case 2.} $\gamma \not \subset W_1$.
\smallskip

Fix
\begin{align*}
s_0 &\defeq \sup \{t : t \in [0,a] \,\text{such that}\, \gamma(s) \in W_{1} \forall s \in [0, t]\} \text{ and}\\
t_0 &\defeq \inf \{t : t \in [0,a] \,\text{such that}\, \gamma(s) \in W_{1} \forall s \in [t, a]\}.
\end{align*}

By construction, $\left.\gamma\right|_{[0,s_0)} \subset W_1$ and $\left.\gamma\right|_{(t_0, a]} \subset W_1$. Then, by the Case~1, $\left.\gamma\right|_{[0,s_0)}$ and $\left.\gamma\right|_{(t_0,a]}$ are $(1, \kappa_0)$-quasi-geodesic of $(U,K_U)$. Since $\k_0$ depends only on $W$, $W_1$, $W_2$, $U$ and $\kappa$, by using sequential argument, we get that $\left.\gamma\right|_{[0,s_0]}$ and $\left.\gamma\right|_{[t_0,a]}$ are $(1, \kappa_0)$-quasi-geodesic of $(U,K_U)$. Let $\sigma_{1}: [0,a_1] \rightarrow U$ and $\sigma_{2}: [0,a_2] \rightarrow U$ be $(1,\kappa_0)$-almost-geodesic of $(U,K_{U})$ with $\sigma_1(0)=z$ and $\sigma_1(a_1)=\gamma(s_0)$, and $\sigma_2(0)=w$ and $\sigma_2(a_2)=\gamma(t_0)$. Since $W_1$ is an open set, so $\gamma(s_0), \gamma(t_0) \in \OM \setminus W_1$. Since $W$ and $U \setminus W_1$ have weak visibility property with respect to $(U, K_{U})$, so it follows that there exists compact set $K \subset U$ such that $\sigma_1 \cap K \neq \emptyset$ and $\sigma_2 \cap K \neq \emptyset$.
\smallskip
       
Let $o_1 \in \sigma_1 \cap K$ and $o_2 \in \sigma_2 \cap K$, and $o \in U$. Now, using the fact that $\sigma_1$ and $\sigma_2$ are $(1,\kappa_0)$-almost-geodesic of $(U,K_{U})$ and the triangle inequality, we have the following:
\begin{align}
K_{U}(z, \gamma(s_0))&\geq K_{U}(z, o_1) + K_{U}(o_1, \gamma(s_0)) - 3\kappa_0 \notag\\
                    &\geq K_{U}(z, o) - K_{U}(o, o_1) + K_{U}(o, \gamma(s_0)) - K_{U}(o, o_1) - 3\kappa_0 \notag\\
                    & \geq K_{U}(z, o) + K_{U}(o, \gamma(s_0)) - 2 \sup_{y \in K} K_{U}(o,y) - 3\kappa_0. \label{ineq: U-bound of Kob dist in LLR for sigma1}
\end{align}
and similarly,
       \begin{align}
           K_{U}(w, \gamma(t_0)) \geq K_{U}(w, o) + K_{U}(o, \gamma(t_0)) - 2 \sup_{y \in K} K_{U}(o,y) - 3\kappa_0. \label{ineq: U-bound of Kob dist in LLR for sigma2}
       \end{align}

Now, using the fact that $\gamma$ is $(1,\kappa)$-almost-geodesic of $(\OM, K_\OM)$, $\left.\gamma\right|_{[0,s_0]}$ and $\left.\gamma\right|_{[t_0,a]}$ are $(1, \kappa_0)$-quasi-geodesic of $(U,K_U)$, and by \eqref{ineq: U-bound of Kob dist in LLR for sigma1} and \eqref{ineq: U-bound of Kob dist in LLR for sigma2}, we get
\begin{align*}
K_\OM(z,w) & \geq a - \kappa\\
                        & \geq a - t_{0} + s_{0} - \kappa\\
                        & \geq K_{U}(w, \gamma(t_0)) - \kappa_{0} + K_{U}(\gamma(s_0), z) - \kappa_{0} - \kappa\\
                        & \geq K_{U}(z, o) + K_{U}(o, w) +  K_{U}(o, \gamma(t_0)) + K_{U}(o, \gamma(s_0)) - 4 \sup_{y \in K} K_{U}(o,y) - 8\k_0 -\k\\
                        & \geq  K_{U}(z, w) - 4 \sup_{y \in K} K_{U}(o,y)- 8\k_0-\k.
\end{align*}
Thus, we have
\begin{align*}
           K_{U}(z,w) \leq K_\OM(z,w) + C_{2}
\end{align*}
where $C_2 \defeq  4 \sup_{y \in K} K_{U}(o,y) + 8\k_0+\k < + \infty$.
\smallskip

Therefore, by taking $C \defeq \max(C_1, C_2)$, we get that $K_{U}(z,w) \leq K_\OM(z,w) + C$ for every $z, w \in W$. This completes the proof.       
\end{proof}

\subsection{Global localization result}
The following lemma is analogous to the part $(3)$ of Proposition~3.2 in \cite{ADS2023} which we prove in our setting.
\begin{lemma}\label{PDwrtGKD}
Suppose $\OM \subset \C^d$ is a Kobayashi hyperbolic domain and $U \subset \OM$ is a sub-domain. Let $W \subset U$ such that $\overline{W} \cap \OM \subset U$. If $W$ and $\OM \setminus U$ have weak visibility property with respect to $(\OM , K_\OM)$, then $K_\OM(W, \OM \setminus U) > 0$.
\end{lemma}

\begin{proof}
Assume that $K_\OM(W, \OM \setminus U) = 0$. Then there exist sequences of points $z_n \in W$ and $w_n \in \OM \setminus U$ such that $K_\OM(z_n, w_n) \rightarrow 0$ as $ n \rightarrow \infty$. Let, for each $n \in \N$, $\gamma_n : [0, a_n] \rightarrow \OM$ be a $(1, 1/n)$-almost-geodeisc of $(\OM, K_\OM)$ such that $\gamma_n(0)=z_n$ and $\gamma_n(a_n)=w_n$. Since $W$ and $\OM \setminus U$ have weak visibility property with respect to $(\OM , K_\OM)$, there exists a compact set $K \subset \OM$ and $N\in\N$ such that $\gamma_n \cap K \neq \emptyset$ for all $n\geq N$. Let $o_n \in \gamma_n \cap K$. Since $K$ is compact set, up to a subsequence, we may assume that $o_n \rightarrow o \in K$ as $n \rightarrow \infty$. Let $o_n = \gamma_{n}(t_n)$ for some $t_n \in [0,a_n]$. Since $\gamma_{n}$ is $(1, 1/ n)$-almost-geodesic of $(\OM, K_\OM)$ passing through $o_n$, we get
\begin{align*}
K_\OM(z_n, w_n)\geq K_\OM(z_n, o_n) + K_\OM(o_n, w_n) - 3/n.
\end{align*}
The inequality above implies that as $n \rightarrow \infty $, we get $K_\OM(z_n, o_n) \rightarrow 0$ and $K_\OM(o_n, w_n) \rightarrow 0$. Since $\OM$ is Kobayashi hyperbolic, so $z_n \rightarrow o$ and $w_n \rightarrow o$ as $n \rightarrow \infty$. Therefore, we get that $o \in \overline{W}$ and $ o \in \overline{\OM \setminus U}$. Since $\overline{W} \cap \OM \subset U$, so it follows that $o \in U$. Hence, we have $o \in U \cap \overline{\OM \setminus U}$. This is a contradiction because $U \cap \overline{\OM \setminus U} = \emptyset$. This completes the proof.
\end{proof}

\begin{lemma}\label{BDforanypointwrtGD}
Suppose $\OM \subset \C^d$ is a Kobayashi hyperbolic domain and $U \subset \OM$ is a sub-domain. Then for every set $W$ with $W \subset U$ such that $W$ and $\OM \setminus U$ have weak visibility property with respect to $(\OM , K_\OM)$ and $o \in \OM$, there exists $L > 0$ such that for all $z \in W$
\begin{align*}
K_\OM(z,o) \leq K_\OM(z,\OM \setminus U) + L.
\end{align*}
\end{lemma}
\begin{proof}
The proof follows exactly same as the proof of Lemma~\ref{Fixed Point local Localization}.
\end{proof}
The next lemma is analogous to Lemma~3.6 in \cite{ADS2023}. The proof follows the same idea as that of Lemma~\ref{AG wrt U for Local}, but it is in fact simpler. We therefore omit the routine calculations; see also Lemma~3.6 in \cite{ADS2023}.
\begin{lemma}\label{Lem:Quasi-geo_wrt_Small_Domain}
Suppose $\OM \subset \C^d$ is a Kobayashi hyperbolic domain and $U \subset \OM$ is a sub-domain. Let $W \subset U$ such that $\overline{W} \cap \OM \subset U$ and such that $W$ and $\OM \setminus U$ have weak visibility property with respect to $(\OM , K_\OM)$. Let $z, w \in W$ and let $\gamma$ be a $(1, \kappa)$-almost-geodesic of $(\OM, K_\OM)$ with $\gamma(0)=z$ and $\gamma(a)=w$ such that $\gamma \subset W$. Then there exists a constant $\k_0 > 0$ such that $\gamma$ is a $(1, \kappa_0)$-quasi-geodesic of $(U, K_U)$.
\end{lemma}

\begin{proof}
First we observe that, by Lemma \ref{PDwrtGKD}, we have $K_\OM(W, \OM \setminus U) > 0$. Now, by  Lemma \ref{Lem:Loc_of_metric}, there exists $L > 0$ such that for all $z \in W$ and $v \in \C^d$,
\begin{equation}\label{ineq:PFofQGwithL1}
\kappa_{U}(z;v) \leq \left(1 + Le^{-K_\OM(z, \OM \setminus U)}\right)\kappa_{\OM}(z;v).
\end{equation}
Note that $\gamma \subset W$. Using the fact that the Kobayashi distance is the integrated form of the Kobayashi metric, from~\eqref{ineq:PFofQGwithL1} and the fact that $\gamma$ is $(1, \kappa)$-almost-geodesic of $(\OM, K_\OM)$, we get for every $s_1 \leq s_2 \in [0,a]$
\[
K_{U}(\gamma(s_1), \gamma(s_2))\leq K_\OM(\gamma(s_1), \gamma(s_2)) + \kappa + L \int_{0}^{a}e^{-K_\OM(\gamma(t), \OM \setminus U)}dt.
\]
For $o \in \OM$, by Lemma \ref{BDforanypointwrtGD}, there exists $L_o > 0$ such that for every $z \in W$, we have 
\begin{equation}\label{ineq:PFofQGwithL2}
K_\OM(z,o) \leq K_\OM(z,\OM \setminus U) + L_o.
\end{equation}
By the continuity of Kobayashi distance, there exists $t_0 \in [0,a]$ such that for all $t \in [0,a]$
\begin{equation}\label{ineq:PFofQGwithL3}
K_\OM(\gamma(t_0), \OM \setminus U) \leq K_\OM(\gamma(t), \OM \setminus U).
\end{equation}
Using the fact that $\gamma$ is $(1,\kappa)$-almost-geodesic of $(\OM, K_\OM)$, by triangle inequality, by~\eqref{ineq:PFofQGwithL2} and \eqref{ineq:PFofQGwithL3}, we have 
\[
|t - t_0| - \kappa \leq 2K_\OM(\gamma(t), \OM \setminus U) + 2L_o.
\]
Thus, for all $t \in [0, a]$, we have
\begin{equation}\label{ineq:PFofQGwithL4}
-K_\OM(\gamma(t), \OM \setminus U) \leq \frac{2L_o + \kappa}{2} - \frac{|t - t_0|}{2}.
\end{equation}
Now, using~\eqref{ineq:PFofQGwithL4}, the fact that $\gamma$ is $(1, \kappa)$-almost-geodesic of $(\OM, K_\OM)$ and the estimate $\int_{0}^{+\infty} e^{\frac{-|t - t_0|}{2}}dt \leq 4$, for all $s_1, s_2 \in [0,a]$, we get
\[
K_{U}(\gamma(s_1), \gamma(s_2))\leq |s_1 - s_2| + 2\kappa + 4LC
\]
where $C\defeq e^{\frac{2L_o + \kappa}{2}}$. 
Therefore, $\gamma$ is $(1,\kappa_0)$-quasi-geodesic of $(U, K_U)$ for $\kappa_0 \defeq 2\kappa + 4LC$. 
\end{proof}

We now present the proof of Theorem~\ref{Thm:Global_Loc_Res}.
\begin{proof}[Proof of {\Cref{Thm:Global_Loc_Res}}]
Let $z, w \in W$ and let $\kappa > 0$. Let $\gamma : [0,a] \rightarrow \OM$ be a $(1,\kappa)$-almost-geodesic of $(\OM, K_\OM)$ joining $z$ and $w$.
\smallskip
        
\noindent {\bf Case~1.} $\gamma \subset W_1$.
\smallskip
        
Then, by replacing $W$ with $W_1$ in Lemma \ref{Lem:Quasi-geo_wrt_Small_Domain}, we get that there exists $\kappa_0 > 0$ such that $\gamma$ is $(1,\kappa_0)$-quasi-geodesic of $(U,K_U)$. Thus, using $\gamma$ is $(1,\kappa_0)$-quasi-geodesic of $(U,K_U)$ and $(1,\kappa)$-almost-geodesic of $(\OM, K_\OM)$, we get the following inequality,
\[
K_{U}(z,w) \leq K_\OM(z,w) + \kappa + \kappa_0.
\]
Take $C_1 \defeq \kappa + \kappa_0$. Therefore, in this case, we are done.
\smallskip
        
\noindent {\bf Case~2.} $\gamma \not \subset W_1$
\smallskip

Now, we set 
\begin{align*}
s_0 &\defeq \sup\{t : t \in [0,a] \,\text{such that} \, \gamma(s) \in W_1 \, \forall s \in [0,t]\} \text{ and}\\
t_0 &\defeq \inf\{t : t \in [0,a] \,\text{such that} \, \gamma(s) \in W_1 \, \forall s \in [t,a]\}.
\end{align*}
By construction, $\left.\gamma\right|_{[0,s_0)} \subset W_1$ and $\left.\gamma\right|_{(t_0, a]} \subset W_0$. Then, by the Case 1, $\left.\gamma\right|_{[0,s_0)}$ and $\left.\gamma\right|_{(t_0,a]}$ are $(1, \kappa_0)$-quasi-geodesic of $(U,K_U)$. Since $\k_0$ depends only on $W$, $W_1$, $U$ and $\kappa$, by using sequential argument, $\left.\gamma\right|_{[0,s_0]}$ and $\left.\gamma\right|_{[t_0,a]}$ are $(1, \kappa_0)$-quasi-geodesic of $(U,K_U)$. Since $W_1$ is an open set, so $\gamma(s_0), \gamma(t_0) \in \OM \setminus W_1$. Since $W$ and $\OM \setminus W_1$ have weak visibility property with respect to $(\OM, K_\OM)$, there exists compact set $K \subset \OM$ such that $\left.\gamma\right|_{[0,s_0]} \cap K \neq \emptyset$ and $\left.\gamma\right|_{[t_0, a]} \cap K \neq \emptyset$.
\smallskip
        
Let $o_1 \in \left.\gamma\right|_{[0,s_0]} \cap K$ and $o_2 \in \left.\gamma\right|_{[t_0, a]} \cap K$. Now, using the fact that $\gamma$ is $(1,\kappa)$-almost-geodesic of $(\OM, K_\OM)$, $\left.\gamma\right|_{[0,s_0]}$ and $\left.\gamma\right|_{[t_0,a]}$ are $(1, \kappa_0)$-quasi-geodesic of $(U,K_U)$, and by the triangle inequality, we get following inequalities 
\begin{align*}
K_\OM(z,w) &\geq a - \kappa\\
                         &\geq K_U(z,\gamma(s_0)) - \kappa_0 + K_U(w, \gamma(t_0)) - \kappa_0 - \kappa\\
                         &\geq K_U(z,o_1) + K_U(o_1, \gamma(s_0))  + K_U(w,o_2) + K_U(o_2, \gamma(t_0)) - 8\kappa_0 - \kappa\\
                         &\geq K_U(z,o_1) + K_U(w,o_2) - 8\kappa_0 - \kappa\\
                         &\geq K_U(z,o_1) + K_U(w,o_1) - K_U(o_1, o_2) - 8\kappa_0 - \kappa\\
                         &\geq K_U(z,w) - K_U(o_1, o_2) - 8\kappa_0 - \kappa.
\end{align*}
Since $K \subset \OM$ is compact, so $\overline{W_1} \cap K$ is compact. Since $\overline{W_1} \cap \OM \subset U$, so it follows that $\overline{W_1} \cap K \subset \overline{W_1} \cap \OM \subset U$. Since the Kobayashi distance is continuous, so it follows that $\sup_{x, y \in \overline{W_1} \cap K}K _U(x,y) < \infty$. Thus we have
\[
K_{U}(z,w) \leq K_\OM(z,w) + C_2
\]
where $C_2 \defeq \sup_{x,y \in \overline{W_1} \cap K}K_U(o_1, o_2) + 8\kappa_0 + \kappa$.
\smallskip

Therefore, by taking $C \defeq \max(C_1, C_2)$, we obtain that, for every $z, w \in W$, $K_{U}(z,w) \leq K_\OM(z,w) + C$. This completes the proof.
\end{proof}

\subsection{Some consequences of additive localization results}
We start this subsection by proving some auxiliary results which are required to prove the result regarding equivalence of local-global intrinsic visibility. We first prove the following lemma.
\begin{lemma}\label{Lem:Rep_lemma}
Suppose $\OM \subset \C^d$ is a Kobayashi hyperbolic domain and $U \subset \OM$ is a sub-domain. Suppose $W,W_1,W_2$ are sets such that $W_2$ is open set and such that
\begin{itemize}
\item[(i)] $W \subset W_1 \subset W_2 \subset U$ and $K_\OM(W_2, \OM \setminus U) > 0$, and
\item[(ii)] $W_1$ and $U \setminus W_2$ have weak visibility property with respect to $(U, K_U)$.
\end{itemize}
Let $z,w\in W_1$ and let $\gamma:[0,a]\to\OM$ be a $(1,\k)$-almost-geodesic of $(\OM,K_\OM)$ that joins $z$ ad $w$ such that $\gamma\subset W_1$. Then there exists $\wt\k>0$ such that $\gamma$, up to a reparameterization, is $(1,\wt\k)$-almost-geodesic for $(U,K_U)$. 
\end{lemma}
\begin{proof}
 We first note that $\gamma$ is $(1,\kappa)$-geodesic for $(\OM,K_\OM)$. Choose $\epsilon>0$ such that $\epsilon\in(0,1)$ and such that $N_{\C^d}(\gamma,\epsilon)\subset\subset\OM$. By Result~\ref{Res:Fattenig_of_(lk)-geo}, there exists an absolutely continuous map $\wt\gamma:[a,b]\to\OM$ with $\wt\gamma(a) = \gamma(a)$ and $\wt\gamma(b) = \gamma(b)$ such that
\begin{itemize}
\item[(i)] $\wt\gamma$ is $(1,\epsilon+\kappa)$-geodesic for $(\OM,K_\OM)$,
\item[(ii)] for every $t\in [a,b]$, $\|\gamma(t)-\wt\gamma(t)\|\leq \epsilon$ and
\item[(iii)] $\wt{\gamma}^\prime$ is almost everywhere non-vanishing.
\end{itemize}
Note that $\wt\gamma$ is $(1,\kappa+1)$-geodesic for $(\OM,K_\OM)$. Using Result~\ref{Res:(lk)-geo_is_(lk)-al-geo}, we can reparametrize $\wt\gamma$ to get $\k_\OM$-unit-speed. By the abuse of notation, we continue to call this reparametrization $\wt\gamma$. Since $K_\OM(W_2,\OM\setminus U)>0$, by Lemma~\ref{Lem:Loc_of_metric}, there exists $L > 0$ such that for all $z \in W_2$ and $v \in \C^d$
\begin{equation}\label{Ineq:Loc_of_metric}
\k_{U}(z;v) \leq \left(1 + Le^{-K_\OM(z, \OM \setminus U)}\right)\k_{\OM}(z;v).
\end{equation}
Since $K_\OM(W_2,\OM\setminus U)>0$, by following the same argument as in the proof of the Claim in Lemma~\ref{AG wrt U for Local}, we get that there exist $C_0,C_1>0$ such that for all $z \in W_1$
\begin{equation}\label{Ineq:Dist_of_sets_Comp}
K_U(z,U\setminus W_2)\leq C_0 K_\OM(z,\OM\setminus W_2) + C_1.
\end{equation}
Since $W_1$ and $U\setminus W_2$ have weak visibility property with respect to $(U,K_U)$, for $o\in U$, by Lemma~\ref{Fixed Point local Localization}, there exists $L_0>0$ such that for all $z\in W_1$
\begin{equation}\label{Ineq:Comp_of_Dist_of_points_with_set}
K_U(z,o)\leq K_U(z,U\setminus W_2) + L_0.
\end{equation}
We now make the following claim.
\smallskip

\noindent{\bf Claim.} There exists $\wt\k>0$ such that $\wt\gamma$ is $(1,\wt\k)$-geodesic for $(U,K_U)$.
\smallskip

\noindent{\itshape Proof of Claim.} Using the fact that the Kobayashi distance is the integrated form of the Kobayashi metric, by~\ref{Ineq:Loc_of_metric} and $\kappa_\OM(\wt\gamma(t);\wt\gamma^\prime(t))=1$ almost everywhere, for all $s_1\leq s_2$, we get
\[
 \int_{s_1}^{s_2} \kappa_{U}(\gamma(t);\gamma^{\prime}(t))dt \leq |s_2 - s_1| + L\int_{0}^{a} e^{-K_\OM(\gamma(t), \OM \setminus W)}dt.
\]
Now, using the fact that $\wt\gamma$ is $(1,\kappa+1)$-geodesic of $(\OM,K_\OM)$, by triangle inequality, and by~\eqref{Ineq:Comp_of_Dist_of_points_with_set} and~\eqref{Ineq:Dist_of_sets_Comp}, we get the following inequalities, for all $t \in [0,a]$
\[
|t - t_0| - \kappa-1 \leq 2C_0 K_\OM(\wt\gamma(t), \OM \setminus W) + 2L_0 + 2C_1.
\]
Therefore, for all $t \in [0,a]$
\begin{align}
-K_\OM(\gamma(t), \OM \setminus W)  \leq \frac{-|t - t_0|}{2C_0} + \frac{2L_0 + 2C_1 + \kappa+1}{2C_0} \label{ineq: U-bound for the kob dist}
\end{align}
Therefore, we get for all $s_1, s_2 \in [0,a]$,
\[
K_{U}(\wt\gamma(s_1), \wt\gamma(s_2))\leq |s_2 - s_1| + 4LC_{2} C_{0}
\]
where $C_2 \defeq e^{\frac{2L_0 + 2C_1 + \kappa+1}{2C_0}}$. Now, again using the fact that $\wt\gamma$ is $(1,\k+1)$-geodesic for $(\OM,K_\OM)$, we get
\[
K_{U}(\wt\gamma(s_1), \wt\gamma(s_2)) \leq K_U(\wt\gamma(s_1), \wt\gamma(s_2))+ \k+1+ 4LC_0C_2.
\]
Therefore, by taking $\wt\kappa \defeq4LC_{0} C_{2}+\k+1$, $\wt\gamma$ is $(1,\wt\kappa)$-geodesic for $(U,K_U)$.
\smallskip
\hfill$\blacktriangleleft$

Now, using Result~\ref{Res:(lk)-geo_is_(lk)-al-geo}, we can again reparametrize $\wt\gamma$ so that $\wt\gamma$ is $(1,\wt\k)$-almost-geodesic for $(U,K_U)$. This completes the proof of the lemma.
\end{proof}

\begin{theorem}\label{Thm:Local_to_global_vis}
Suppose $\OM \subset \C^d$ is a Kobayashi hyperbolic domain and $U \subset \OM$ is a sub-domain. Suppose $W, W_1, W_2$ are sets such that $W_1$ and $W_2$ are open and such that
\begin{itemize}
\item[(i)] $W\subset W_1\subset W_2\subset U$ and $K_\OM(W_2, \OM\setminus U)>0$
\item[(ii)] $W$ and $U\setminus W_1$, and $W_1$ and $U\setminus W_2$ have weak visibility property with respect to $(U,K_U)$
\end{itemize}
Let $A \subset W$ and $B \subset U$. If $A$ and $B$ have weak visibility property with respect to $(U,K_U)$, then $A$ and $B$ have weak visibility property with respect to $(\OM, K_\OM)$.
\end{theorem}
\begin{proof}
Assume, to get a contradictions, that $A$ and $B$ does not have visibility property with respect to $(\OM,K_\OM)$. Then there exist $\kappa> 0$ and a sequence $\gamma_n:[0,a_n]\to \OM$ of $(1,\k)$-almost-geodesic for $(\OM,K_\OM)$ such that $\gamma_n(0)\in A$, $\gamma_n(a_n)\in B$ and such that $(\gamma_n)_{n \geq 1}$ eventually avoids every compact subset of $\OM$.  We may have following two cases:
\smallskip
    
\noindent{\bf Case~1.} $\gamma_{n} \subset W_1$ for all $n \in \N$.
\smallskip

In this case, by Lemma~\ref{Lem:Rep_lemma}, there exists $\wt\k>0$ such that each $\wt\gamma_n$, up to a reparametriation, is $(1,\wt\k)$-geodesic for $(U,K_U)$ such that $(\wt{\gamma}_n)_{n\geq 1}$ eventually avoids every compact subset of $\OM$. This gives a contradiction to weak visibility of $A$ and $B$ with respect to $(U,K_U)$.
\smallskip
    
\noindent{\bf Case~2.} There exists a subsequence $(n_{k})_{k\geq 1}$ such that $\gamma_{n_{k}} \not\subset W_1$ for all $k \geq 1$.
Now, set
\[
t_{n_{k}}\defeq \sup\{t\in [0, a_{n_{k}}]: \gamma_{n_{k}}(s) \in W_1 \, \forall s \in[0,t]\}.
 \]
Since $W_1$ is an open, by construction, $\gamma_{n_{k}}(t_{n_{k}}) \in \OM\setminus W_1$ and $\wt{\gamma}_{n_{k}}\defeq \left.\gamma_{n_{k}}\right|_{[0, t_{n_{k}})} \subset W_1$. By Case~1, $\wt{\gamma}_{n_{k}}$ is $(1,\wt\kappa)$-almost-geodesic of $(U, K_U)$. Therefore, $(\wt{\gamma}_{n_{k}})_{k\geq 1}$ eventually avoids every compact subset of $U$. Note that $\wt\gamma_n(0)\in A\subset W_1$ and $\gamma_{n_{k}}(t_{n_{k}}) \in \OM\setminus W_1$ for all $k \geq 1$. This contradicts the visibility of $W$ and $\OM\setminus W_1$ with respect to $(U,K_U)$.
\end{proof}

\begin{theorem}\label{Thm:Global_to_loc_vis}
Suppose $\OM \subset \C^d$ is a Kobayashi hyperbolic domain and $U \subset \OM$ is a sub-domain. Suppose $W, W_1, W_2, W_3$ are sets such that $W_1$ and $W_3$ are an open set and such that
\begin{itemize}
\item[(i)] $W \subset W_1 \subset W_2\subset W_3 \subset  U$ and $\overline{W_3} \cap \OM \subset U$, and
\item[(ii)] $W$ and $\OM\setminus W_1$,  $W_2$ and $\OM\setminus W_3$, and $W_3$ and $\OM\setminus U$ have weak visibility property with respect to $(\OM, K_\OM)$.
\end{itemize}
 Let $A \subset W$ and $B \subset U$. If $A$ and $B$ have weak visibility property with respect to $(\OM, K_\OM)$, then $A$ and $B$ have weak visibility property with respect to $(U, K_{U})$.
\end{theorem}
\begin{proof}
 By Theorem~\ref{Thm:Global_Loc_Res}, there exists $C>0$ such that for all $z, w \in W_2$
\begin{equation}\label{ineq:Localization_Ineq_in_Main_Thm1_00}
        K_{U}(z,w)\leq K_\OM(z,w) + C.
\end{equation}
Assume, to get a contradiction, that $A$ and $B$ does not have visibility property with respect to $(U, K_{U})$. Then there exists a sequence $\gamma_{n}:[0, a_{n}] \to U$ of $(1,\kappa)$-almost-geodesic of $(U, K_{U})$ with $\gamma_{n}(0) \in A$ and $\gamma_{n}(a_{n}) \in B$ such that $(\gamma_{n})_{n \geq 1}$ eventually avoids all compact subset of $U$. We may have following two cases:
\smallskip
    
\noindent{\bf Case~1.} $\gamma_{n} \subset W_1$ for all $n \in \N$.
\smallskip

By~\eqref{ineq:Localization_Ineq_in_Main_Thm1_00}, for all $s, t\in [0,a_{n}]$
\begin{align*}
K_\OM(\gamma_{n}(s), \gamma_{n}(t)) &\geq K_{U}(\gamma_{n}(s), \gamma_{n}(t)) - C\\
                                                        &\geq |s-t| -\kappa -C.
\end{align*}
Moreover, $K_\OM(\gamma_{n}(s), \gamma_{n}(t)) \leq K_{U}(\gamma_{n}(s), \gamma_{n}(t)) \leq |s-t| + \kappa + C$ for all $s, t\in [0, a_{n}]$. That is, $\gamma_{n}$ is $(1,\kappa_{0})$-almost-geodesic of $(\OM, K_\OM)$ where $\kappa_{0}\defeq \kappa + C$. Therefore, $(\gamma_n)_{n\geq 1}$ eventually avoids every compact subset of $\OM$, which contradicts the visibility of $A$ and $B$ with respect to $(\OM, K_\OM)$.

\noindent{\bf Case~2.} There exists a subsequence $(n_{k})_{k\geq 1}$ such that $\gamma_{n_{k}} \not\subset W_1$ for all $k \geq 1$.

Now, set
\[
t_{n_{k}}\defeq \sup\{t\in [0, a_{n_{k}}]: \gamma_{n_{k}}(s) \in W \, \forall s \in[0,t]\}.
 \]
By construction, $\gamma_{n_{k}}(t_{n_{k}}) \in \partial W_1\cap\OM$ for all $k \in \N$ and $\wt{\gamma}_{n_{k}}\defeq \left.\gamma_{n_{k}}\right|_{[0, t_{n_{k}})} \subset W_1$. By Case~1, $\wt{\gamma}_{n_{k}}$ is $(1,\kappa_{0})$-almost-geodesic of $(\OM, K_\OM)$. Note that $\wt\gamma_{n_k}(0)\in W$. Since $W_1$ is an open set, so it follows that $\gamma_{n_{k}}(t_{n_{k}}) \in \OM\setminus W_1$ for all $k \geq 1$. Therefore, $(\wt{\gamma}_{n_{k}})_{k\geq 1}$ eventually avoids every compact subset of $\OM$. This contradicts the weak visibility property of $W$ and $\OM\setminus W_1$ with respect to $(\OM, K_\OM)$. This completes the proof.
\end{proof}

\begin{remark} 
Theorem~\ref{Thm:Global_to_loc_vis} require the global localization result and weak visibility property of $W$ and $\OM\setminus W_1$ with respect to $(\OM,K_\OM)$. On the other hand, Theorem~\ref{Thm:Local_to_global_vis} require the reparametrization of $(1,\k)$-almost-geodesic of $(\OM,K_\OM)$ so that it becomes $(1,\wt\k)$-almost-geodesic for $(U,K_U)$ for some $\wt\k>0$. Moreover, it also required that $K_\OM(W_2,\OM\setminus U)>0$ which is weaker assumption than $\OM$ being hyperbolically embedded (see \cite{Masanta,CGMS2024}).
\end{remark}

Theorem~\ref{Thm:Local_to_global_vis} and Theorem~\ref{Thm:Global_to_loc_vis} establish the equivalence of intrinsic local-global visibility. The equivalence of standard local-global visibility was established in \cite{Masanta, CGMS2024}. In fact, it is easy to see that equivalence of standard local-global visibility can be proved using the equivalence of intrinsic local-global visibility.

We now state the following theorem which shows that additive localization is stronger than multiplicative localization. This theorem is analogous to Theorem~1.8 in \cite{ADS2023} in our setting. We shall use it in the proof of Condition~(M$\star$) implies Condition~($\star$), i.e., Proposition~\ref{Pro:Cond(M*)_implies_Cond(*)}.
\begin{theorem}\label{Thm:Mult_localization}
Suppose $\OM \subset \C^d$ is a Kobayashi hyperbolic domain and $U \subset \OM$ is a sub-domain. Suppose $W,W_1,W_2$ are sets with $W_1$ and $W_2$ are open sets such that $W \subset W_1 \subset W_2 \subset U$, $K_\OM(W_1, \OM \setminus W_2) > 0$, and
there exists a constant $C_0 > 0$ such that for all $z , w \in W$, 
\[
K_U(z,w) \leq K_\OM(z,w ) + C_0.
\]
Then there exists a constant $C\geq 1$ such that for all $z,w\in W$
\[
K_U(z,w) \leq C K_\OM(z,w ).
\]
\end{theorem}
\begin{proof}
The proof is identical to that of Theorem~1.8 in \cite{ADS2023}, and we therefore omit it.
\end{proof}

\section{Busemann horosphere}\label{Sec:Busemann_horo}
We start this section by studying some basic properties of the Busemann horosphere. Given a geodesic ray $\gamma:[0,\infty)\to\OM$ of $(\OM,K_\OM)$, consider the {\em cluster set} of $\gamma$
\[
\mathcal{C}(\gamma):=\{x\in\bdy\OM: \exists (t_n)_{n\geq 1}\subset [0,\infty), t_n\to\infty, \text{ such that }\lim_{n\to\infty}\gamma(t_n)=x\}.
\]
The following lemma says that the cluster set of geodesic ray is subset of closure of the associated Busemann horosphere intersection of the boundary of domain.
\begin{lemma}\label{Lem:Cluster_set_of_geo_sits_in_geo_horosphere}
Suppose $\OM\subset \C^{d}$ is complete hyperbolic domain. If $\gamma:[0,\infty) \to \OM$ is a geodesic ray of $(\OM,K_\OM)$, then $\mathcal{C}(\gamma) \subset \clos{H_{\OM}(\gamma,R)}\cap \partial\OM$ for all $R>0$.
\end{lemma}

\begin{proof}
If $\mathcal{C}(\gamma)=\emptyset$, then we are done. Suppose $\mathcal{C}(\gamma)\not=\emptyset$. In this case, we assume to the contrary that $\mathcal{C}(\gamma)\not\subset\clos{H_{\OM}(\gamma,R)}$ for some $R>0$. Then there exists $x \in \mathcal{C}(\gamma)$ such that $x \not\in \clos{H_{\OM}(\gamma,R)}\cap \partial\OM$. There exists $r > 0$ such that $B(x,r)\cap H_{\OM}(\gamma,R) = \emptyset$. That is, for all $z \in B(x,r)\cap\OM$, we have
\[
\lim_{s\to\infty}[K_\OM(z, \gamma(t)) - K_\OM(\gamma(0), \gamma(t))] \geq \frac{1}{2}\log R.
\]
There exists a sequence $(t_{n})_{n\geq 1} \subset [0,\infty)$ with $t_{n} \to \infty$ as $n \to \infty$ such that $\gamma(t_{n}) \to x$ as $n \to \infty$.Without loss of generality, we may assume that $\gamma(t_{n}) \in B(x,r)\cap\OM$ for all $n \in \N$. Since for all $n \in \N$,
\[
K_\OM(\gamma(t_{n}), \gamma(t)) - K_\OM(\gamma(0), \gamma(t)) = -t_{n}. 
\]
Therefore, for all $n \in \N$, we obtain
\[
\lim_{t\to \infty}[K_\OM(\gamma(t_{n}), \gamma(t)) - K_\OM(\gamma(0), \gamma(t))] = -t_{n} \geq \frac{1}{2}\log R.
\]
Note that $t_{n} \to \infty$ as $n \to \infty$, so we get a contradiction. This completes the proof.
\end{proof}

The next lemma provides the relation between the Busemann horosphere of geodesic ray and the Busemann horosphere of truncated geodesic ray.
\begin{lemma}\label{Lem:Relation_between_origin_geo_horo_and_truncated_geo_horo}
    Suppose $\OM \subset \C^d$ is complete hyperbolic domain and $\gamma:[0,\infty) \to \OM$ is geodesic ray of $(\OM, K_\OM)$. If $\wh{\gamma}:=\left.\gamma\right|_{[T,\infty)}$ for some $T>0$, then for all $R>0$
\begin{equation*}
H_{\OM}(\wh{\gamma}, R)=H_{\OM}(\gamma, e^{-2T}R) \text{ and } H_{\OM}(\wh{\gamma}, e^{2T}R)=H_{\OM}(\gamma, R)
\end{equation*}
\end{lemma}

\begin{proof}
First, we reparametrize $\wh\gamma$ as $\wh{\gamma}(t)= \gamma(t+T)$ for $t \in [0,\infty)$. For all $t \in [0,\infty)$
        \[
            K_\OM(z, \wh{\gamma}(t)) - K_\OM(\wh\gamma(0),\wh{\gamma}(t)) = K_\OM(z, \gamma(t+T)) - K_\OM(\gamma(T),\gamma(t+T)) = K_\OM(z, \gamma(t+T)) - t.
        \]
        Let $s:= t + T$. Therefore
        \[
            K_\OM(z, \gamma(t+T)) - t = K_\OM(z, \gamma(s)) - s + T.
        \]
        By taking limit $t \to \infty$, we get
        \[
             \lim_{t\to\infty}[K_\OM(z, \gamma(t+T)) - t] = \lim_{s\to \infty}[K_\OM(z, \gamma(s)) - s] + T.
        \]
        Hence
\begin{align*}
\lim_{t\to\infty}[K_\OM(z, \gamma(t+T)) - t] < \frac{1}{2}\log R &\iff \lim_{s\to\infty}[K_\OM(z, \gamma(s)) - s] + T< \frac{1}{2}\log R\\
            &\iff \lim_{s\to\infty}[K_\OM(z, \gamma(s)) - s]< \frac{1}{2}\log R -T\\
            &\iff \lim_{s\to\infty}[K_\OM(z, \gamma(s)) - s]< \frac{1}{2}\log e^{-2T}R.
\end{align*}
 Similarly, we get
 \[
 \lim_{t\to\infty}[K_\OM(z, \gamma(t+T)) - t] < \frac{1}{2}\log e^{2T}R\iff \lim_{s\to\infty}[K_\OM(z, \gamma(s)) - s]< \frac{1}{2}\log R.
 \]
Thus, $H_\OM(\wh\gamma,R)=H_\OM(\gamma,e^{-2T}R)$ and $H_{\OM}(\wh{\gamma}, e^{2T}R)=H_{\OM}(\gamma, R)$.  This completes the proof of lemma.
\end{proof}

\begin{lemma}\label{Lem:Asym_compares_two_geo_horo}
Let $\OM \subset\C^{d}$ be a complete hyperbolic domain. Suppose $\gamma_{1}, \gamma_{2}:[0,+\infty) \to \OM$ are two geodesics of $(\OM, K_\OM)$ such that there exists a constant $C>0$ such that for all $t \geq 0$
\begin{equation}\label{Ineq:Asym_of_two_geo_rays}
K_\OM(\gamma_{1}(t), \gamma_{2}(t))\leq C.
\end{equation}
Then there exists a constant $C^{\prime}:= C^{\prime}(C, K_\OM(\gamma_{1}(0), \gamma_{2}(0)))>0$ such that for all $R>0$
\begin{equation*}
H_{\OM}(\gamma_{1}, R) \subseteq H_{\OM}(\gamma_{2}, C^{\prime}R)
\end{equation*}
\end{lemma}
\begin{proof}
For every $t \geq 0$, consider the following inequality which is obtained using the triangle inequality
\begin{equation}\label{Ine:Geo_horo_ine_1}
K_\OM(\gamma_{1}(0), \gamma_{1}(t)) \leq K_\OM(\gamma_{1}(0), \gamma_{2}(0)) + K_\OM(\gamma_{2}(0), \gamma_{2}(t)) + K_\OM(\gamma_{1}(t), \gamma_{2}(t))
\end{equation}
For any $z \in \OM$, using the triangle inequality, \eqref{Ineq:Asym_of_two_geo_rays} and~\eqref{Ine:Geo_horo_ine_1}, we get
\begin{align*}
K_\OM(z, \gamma_{2}(t)) - K_\OM(\gamma_{2}(0), \gamma_{2}(t)) &\leq K_\OM(z, \gamma_{1}(t)) + K_\OM(\gamma_{1}(t),  \gamma_{2}(t)) + K_\OM(\gamma_{1}(0), \gamma_{2}(0))\\
& + K_\OM(\gamma_{1}(t), \gamma_{2}(t)) - K_\OM(\gamma_{1}(0), \gamma_{1}(t))\\
& \leq K_\OM(z, \gamma_{1}(t)) - K_\OM(\gamma_{1}(0), \gamma_{1}(t)) + 2C + K_\OM(\gamma_{1}(0), \gamma_{2}(0)).
\end{align*}
If $z \in H_{\OM}(\gamma_{1}, R)$, then taking $t \to \infty$, we obtain
\[
\lim_{t\to \infty}[K_\OM(z, \gamma_{2}(t)) - K_\OM(\gamma_{2}(0), \gamma_{2}(t))] \leq \frac{1}{2}R + 2C + K_\OM(\gamma_{1}(0), \gamma_{2}(0)).
\]
Therefore, by taking $C^{\prime} := e^{2(2C + K_\OM(\gamma_{1}(0), \gamma_{2}(0)))}$, $H_{\OM}(\gamma_{1}, R) \subseteq H_{\OM}(\gamma_{2}, C^{\prime}R)$.
\end{proof}

In the literature, two geodesic rays $\gamma_1$ and $\gamma_2$ satisfying inequality~\eqref{Ineq:Asym_of_two_geo_rays} are said to be {\em asymptotic}. If $\lim_{t\to\infty}K_\OM(\gamma_1(t),\gamma_2(t))=0$, then $\gamma_1$ and $\gamma_2$ are said to be {\em strongly asymptotic}. The following shows that lemma if two geodesic rays are strongly asymptotic then their corresponding Busemann horospheres are same up to some constant which is a generalization of Lemma~2.3 in \cite{Benini-Bracci}.

\begin{lemma}\label{Lem:Strong_asym_imp_geo_horo_are_same}
Suppose $\OM \subset \C^{d}$ is a complete hyperbolic domain. Suppose that $\gamma_{1}, \gamma_{2} : [0, \infty) \rightarrow \OM$ are geodesic rays of $(\OM, K_\OM)$ such that
\[
\lim_{t \rightarrow \infty}K_\OM(\gamma_{1}(t), \gamma_{2}(t)) = 0.
\]
Then there exists a constant $C:=C(\gamma_{1}, \gamma_{2}) > 0$ such that for all $R > 0$
\[
H_{\OM}(\gamma_{1}, R) = H_{\OM}(\gamma_{2}, CR).
\]
\end{lemma}

\begin{proof}
Fix $p \in \OM$. We make the following claim:
\smallskip

\noindent{\bf Claim.} $\lim_{t \rightarrow \infty}[K_\OM(z, \gamma_{1}(t)) - K_\OM(p, \gamma_{1}(t))] = \lim_{t \rightarrow \infty}[K_\OM(z, \gamma_{2}(t)) - K_\OM(p, \gamma_{2}(t))]$.
\smallskip
        
\noindent {\itshape Proof of Claim.} Fix $t \in [0, \infty)$. By triangle inequality, we obtain
\[
|K_\OM(z, \gamma_{1}(t)) - K_\OM(p, \gamma_{1}(t)) - K_\OM(z, \gamma_{2}(t)) + K_\OM(p, \gamma_{2}(t))| \leq 2 K_\OM(\gamma_{1}(t), \gamma_{2}(t)).
\]
Since $\lim_{t \rightarrow \infty}K_\OM(\gamma_{1}(t), \gamma_{2}(t)) = 0$, by letting $ t \rightarrow \infty$, we obtain the desired equality.
\hfill $\blacktriangleleft$
\smallskip
        
For given geodesic ray $\gamma: [0, \infty) \rightarrow \OM$ of $(\OM, K_\OM)$, define
\[
M(\gamma):= \lim_{t \rightarrow \infty} [K_\OM(\gamma(0), \gamma(t)) - K_\OM(p, \gamma(t))].
\]
Consider
\begin{align*}
\lim_{t \rightarrow \infty}[K_\OM(z, \gamma_{1}(t)) - K_\OM(\gamma_{1}(0), \gamma_{1}(t))] &= \lim_{t \rightarrow \infty}[K_\OM(z, \gamma_{1}(t)) - K_\OM(p, \gamma_{1}(t))] - M(\gamma_{1})\\
&= \lim_{t \rightarrow \infty}[K_\OM(z, \gamma_{2}(t)) - K_\OM(p, \gamma_{2}(t))] - M(\gamma_{1})\\
&= \lim_{t \rightarrow \infty}[K_\OM(z, \gamma_{2}(t)) - K_\OM(\gamma_{2}(0), \gamma_{2}(t))] + M(\gamma_{1}, \gamma_{2})
\end{align*}
where $M(\gamma_{1}, \gamma_{2}):= M(\gamma_{2})- M(\gamma_{1})$ and the second last equality follows from the Claim above. Therefore, for all $R > 0$, we have $H_{\OM}(\gamma_{1}, R) = H_{\OM}(\gamma_{2}, CR)$ where $C:=e^{M(\gamma_{1}, \gamma_{2})}$.
\end{proof}

We now present the proof of Theorem~\ref{Thm:Geo_horosphere_of_small_domain_is_subset_of_geo_horosphere_of_big_domain}.
\begin{proof}[Proof of \Cref{Thm:Geo_horosphere_of_small_domain_is_subset_of_geo_horosphere_of_big_domain}]
Since $\overline{W_1} \cap \OM \subset U, W_{1}$ is an open set, and $W$ and $\OM \setminus W_1$, and $W_1$ and $\OM\setminus U$ have weak visibility property with respect to $(\OM, K_\OM)$, so it follow from Theorem \ref{Thm:Global_Loc_Res} that there exists a constant $C_{0}>0$ such that for all $z,w\in W$,
\begin{equation}\label{ineq:Localization_Ineq_in_Main_Thm1_1}
K_{U}(z,w))\leq K_\OM(z,w) + C_{0}.
\end{equation}
Using inequality \eqref{ineq:Localization_Ineq_in_Main_Thm1_1}, the triangle inequality and inequality \eqref{ineq:Asymp_Condition_for_quasi_geo_1}, we obtain 
    \begin{align*}
         K_\OM(z, \gamma(t)) - K_\OM(\gamma(0),\gamma(t)) &= K_\OM(z, \gamma(t)) - K_{U}(z, \wt{\gamma}(t)) + K_{U}(\wt{\gamma}(0),\wt{\gamma}(t))\\
                                                              &- K_\OM(\gamma(0),\gamma(t)) + K_{U}(z, \wt{\gamma}(t))  - K_{U}(\wt{\gamma}(0),\wt{\gamma}(t))\\
                                                              &\leq K_\OM(z, \gamma(t)) - K_\OM(z, \wt{\gamma}(t)) + K_\OM(\wt{\gamma}(0),\wt{\gamma}(t)) + C_{0}\\
                                                              &- K_\OM(\gamma(0),\gamma(t)) + K_{U}(z, \wt{\gamma}(t))  - K_{U}(\wt{\gamma}(0),\wt{\gamma}(t))\\
                                                              &\leq K_\OM(\gamma(t), \wt{\gamma}(t)) + C_{0} + K_\OM(\wt{\gamma}(0),\wt{\gamma}(t)) - K_\OM(\gamma(0),\gamma(t))\\
                                                              &+ K_{U}(z, \wt{\gamma}(t))  - K_{U}(\wt{\gamma}(0),\wt{\gamma}(t))\\
                                                              &\leq C^{\prime} + C_{0} + K_\OM(\wt{\gamma}(0),\wt{\gamma}(t)) - K_\OM(\gamma(0),\wt{\gamma}(0)) - K_\OM(\gamma(0),\gamma(t))\\
                                                              &+ K_\OM(\gamma(0),\wt{\gamma}(0)) + K_{U}(z, \wt{\gamma}(t))  - K_{U}(\wt{\gamma}(0),\wt{\gamma}(t))\\
                                                              &\leq C^{\prime} + C_{0} + K_\OM(\gamma(0),\wt{\gamma}(t)) - K_\OM(\gamma(0),\gamma(t)) + K_\OM(\gamma(0),\wt{\gamma}(0))\\
                                                              &+ K_{U}(z, \wt{\gamma}(t))  - K_{U}(\wt{\gamma}(0),\wt{\gamma}(t))\\
                                                              &\leq C^{\prime} + C_{0} + K_\OM(\gamma(t), \wt{\gamma}(t)) + K_\OM(\gamma(0),\wt{\gamma}(0)) + K_{U}(z, \wt{\gamma}(t))\\
                                                              &- K_{U}(\wt{\gamma}(0),\wt{\gamma}(t))\\
                                                              &\leq 2C^{\prime} + C_{0} + K_\OM(\gamma(0),\wt{\gamma}(0)) + K_{U}(z, \wt{\gamma}(t)) - K_{U}(\wt{\gamma}(0),\wt{\gamma}(t)).
    \end{align*}
    Let $z \in H_{U}(\wt{\gamma}, R)$. Taking $t \to \infty$ we obtain
    \[
        \lim_{t\to\infty}[K_\OM(z, \gamma(t)) - K_\OM(\gamma(0),\gamma(t))] < 2C^{\prime} + C_{0} +  K_\OM(\gamma(0),\wt{\gamma}(0))+ \frac{1}{2}\log R.
    \]
    That is, $z \in  H_{\OM}(\gamma, CR)$ where $C:= e^{2(C_{0} + 2C^{\prime}) + 2K_\OM(\gamma(0),\wt{\gamma}(0))}$. This completes the proof. 
\end{proof}

The following lemma says that given two sets having intrinsic visibility property and given a geodesic ray and a sequence in the domain such that geodesic is inside one of the set and sequence lies in other set, we have that there is no subsequence of it which can lie inside the corresponding Busemann horosphere.
\begin{lemma}\label{Lem:Tail_of_Seq_cannot_lie_in_geodesic_horosphere}
Suppose $\OM \subset \C^d$ is a complete hyperbolic domain and $A, B \subset \OM$ such that $A$ and $B$ have visibility with respect to $(\OM, K_\OM)$. Suppose that $\gamma : [0, \infty) \rightarrow \OM$ is a geodesic ray of $(\OM, K_\OM)$ such that $\gamma \subset A$. Suppose $(z_{n})_{n \geq 1} \subset B$ is a sequence such that either $\delta_{\OM}(z_{n}) \rightarrow 0$ or $\|z_n\|\to\infty$ as $n \rightarrow \infty$. Then for given $R > 0$, there exists $N \in \N$ such that for all $n \geq N, z_{n} \not \in H_{\OM}(\gamma, R)$. 
\end{lemma}

\begin{proof}
Assume that the conclusion of the lemma is not true. Then there exists $R>0$ and a subsequence $(z_{n_{k}})_{k \geq 1}$ such that $z_{n_{k}} \in H_{\OM}(\gamma, R)$ for all $k \geq 1$. By the abuse of notation, we denote it by $(z_{n})_{n \geq 1}$. Suppose $\gamma_{n} : [0, a_{n}] \rightarrow \OM$ is a geodesic of $(\OM, K_\OM)$ with $\gamma_{n}(0) = z_{n}$ and $\gamma_{n}(a_{n}) = \gamma(t)$. Since $\gamma(t) \in A$ and $A$ and $B$ have visibility property with respect to $(\OM, K_\OM)$, there exists a compact set $K \subset \OM$ such that $\gamma_{n}\cap K \neq \emptyset$. Let $o_{n, t} \in \gamma_{n}\cap K$. Then, up to a subsequence, $(o_{n,k})_{k\in\N}$ converges to $o_n$ such that $o_n\in K$ as $n\in\N$. Let 
\[
S_{n,k}:=K_\OM(o_{n,k},\gamma(t_k))-K_\OM(\gamma(0),\gamma(t_k))
\]
where $t_k\to\infty$ as $k \to\infty$. Since $o_{n,k}\in K$ and $K$ is compact set, it follows that for all $k \in \N, |S_{n,k}|$ is bounded above, say by a constant $C>0$. Therefore, for all $k \geq 0$ and for all $N \in \N$
\begin{align*}
K_\OM(z_{n}, o_{n,k}) &= K_\OM(z_{n}, \gamma(t_k)) - K_\OM(\gamma(0), \gamma(t_k)) - K_\OM(o_{n,t},\gamma(t_k))+K_\OM(\gamma(0),\gamma(t_k))\\
&\leq K_\OM(z_{n}, \gamma(t_k)) - K_\OM(\gamma(0), \gamma(t_k)) + |S_{n,k}|
\end{align*}
Therefore, by taking $k \to \infty$, we obtain for all $n \in \N$
\[
K_\OM(z_{n}, o_{n}) < \frac{1}{2}\log R +C.
\]
Since either $\delta_{\OM}(z_{n}) \rightarrow 0$ or $\|z_n\|\to\infty$ as $n \rightarrow \infty$ and $(\OM,K_\OM)$ is complete, by taking $n \to \infty$, we get a contradiction. This completes the proof.
\end{proof}

We now present the proof of Theorem~\ref{Thm:Geo_horosphere_of_big_domain_is_subset_of_geod_of_small_domain_upto_some_compact_set_and-constant}.
\begin{proof}[Proof of {\Cref{Thm:Geo_horosphere_of_big_domain_is_subset_of_geod_of_small_domain_upto_some_compact_set_and-constant}}] It follows from Theorem~\ref{Thm:Global_to_loc_vis} that there exists a constant $C_0>0$ such that for all $z,w \in W_1$
\begin{equation}\label{ineq:Localization_Ineq_in_Main_Thm1_001}
K_{U}(z,w))\leq K_\OM(z,w) + C_{0}.
\end{equation}

We now make the following claim.
\smallskip
    
\noindent{\bf Claim~1.} There exists a constant $C_{1}>0$ such that $H_{\OM}(\gamma, R)\cap W_1 \subset H_{U}(\wt{\gamma}, C_{1}R)$.
\smallskip
    
\noindent {\itshape Proof of Claim~1.} Let $z \in H_{\OM}(\gamma, R)\cap W_1$. Since $\gamma, \wt{\gamma} \subset W$, using~\eqref{ineq:Localization_Ineq_in_Main_Thm1_001}, and the triangle inequality, we obtain following inequalities 
    \begin{align*}
        K_{U}(z, \wt{\gamma}(t)) - K_{U}(\wt\gamma(0),\wt{\gamma}(t)) &= K_{U}(z, \wt{\gamma}(t)) - K_\OM(z, \gamma(t)) + K_\OM(\gamma(0),\gamma(t))\\
                                                                      &- K_{U}(\wt\gamma(0),\wt{\gamma}(t)) + K_\OM(z, \gamma(t)) - K_\OM(\gamma(0),\gamma(t))\\
                                                                      &\leq  C_{0} + K_\OM(z, \wt{\gamma}(t)) - K_\OM(z, \gamma(t)) + K_\OM(\gamma(0),\gamma(t))\\
                                                                      & - K_\OM(\wt\gamma(0),\wt{\gamma}(t)) + K_\OM(z, \gamma(t)) - K_\OM(\gamma(0),\gamma(t))\\
                                                                      &\leq  C_{0} + K_\OM(\gamma(t), \wt{\gamma}(t)) +  K_\OM(\gamma(0),\gamma(t)) - K_\OM(\gamma(0),\wt{\gamma}(0))\\
                                                                      &- K_\OM(\wt\gamma(0),\wt{\gamma}(t)) + K_\OM(\gamma(0),\wt{\gamma}(0)) + K_\OM(z, \gamma(t))\\
                                                                      &- K_\OM(\gamma(0),\gamma(t))\\
                                                                      &\leq  C_{0} + C^{\prime} + K_\OM(\wt{\gamma}(0),\gamma(t)) - K_\OM(\wt\gamma(0),\wt{\gamma}(t))\\
                                                                      &+ K_\OM(\gamma(0),\wt{\gamma}(0))+ K_\OM(z, \gamma(t)) - K_\OM(\gamma(0),\gamma(t))\\
                                                                      &\leq C_{0} + C^{\prime} + K_\OM(\gamma(t), \wt{\gamma}(t)) + K_\OM(\gamma(0),\wt{\gamma}(0))\\
                                                                      &+ K_\OM(z, \gamma(t)) - K_\OM(\gamma(0),\gamma(t))\\
                                                                      &\leq C_{0} + 2C^{\prime} + K_\OM(\gamma(0),\wt{\gamma}(0)) + K_\OM(z, \gamma(t)) - K_\OM(\gamma(0),\gamma(t)).
    \end{align*}
Taking $t \to \infty$ we obtain
\[
\lim_{t\to\infty}[K_{U}(z, \wt{\gamma}(t)) - K_{U}(\wt\gamma(0),\wt{\gamma}(t))] < C_{0} + 2C^{\prime} +  K_\OM(\gamma(0),\wt{\gamma}(0))+ \frac{1}{2}\log R.
\]
Therefore, $H_{\OM}(\gamma, R)\cap W_1 \subseteq H_{U}(\wt{\gamma}, CR)$ where $C:= e^{2(C_{0} + 2C^{\prime}) + 2K_\OM(\gamma(0),\wt{\gamma}(0))}$.
\hfill $\blacktriangleleft$
\smallskip

\noindent{\bf Claim~2.} $\clos{H_{\OM}(\gamma, R)\setminus W_1}$ is compact subset of $\OM$.
\smallskip

\noindent {\itshape Proof of Claim~2.} Suppose that the claim is not true. There exists a sequence $(z_{n})_{n\geq 1} \subset H_{\OM}(\gamma, R)\setminus W_1$ such that either $\delta_{\OM}(z_{n}) \to 0$ or $\|z_n\|\to\infty$ as $n \to \infty$. Take $A:= W$ and $B:= H_{\OM}(\gamma, R)\setminus W_1$. Since $W$ and $\OM\setminus W_1$ have weak visibility property with respect to ($\OM, K_\OM$), so it follows that $A$ and $B$ have weak visibility property with respect to ($\OM, K_\OM$). Since $\gamma \subset W$, so it follows from Lemma~\ref{Lem:Tail_of_Seq_cannot_lie_in_geodesic_horosphere} that there exists $N \in \N$ such that $z_{n} \not\in H_{\OM}(\gamma, R)\setminus W_1$ for all $n \geq N$. This contradicts to $z_{n} \in H_{\OM}(\gamma, R)\setminus W_1$ for all $n \in \N$.
\hfill $\blacktriangleleft$
\smallskip

Therefore, by Claim~1 and Claim~2, we obtain $H_{\OM}(\gamma, R) \setminus K \subseteq H_{U}(\tilde{\gamma}, CR)$ where $K:=\clos{H_{\OM}(\gamma, R)\setminus W_1}$. This completes the proof.
\end{proof}

We end this section with the following lemma. Fix a holomorphic self map of complete hyperbolic domain. Then the next lemma shows that if two sets have intrinsic visibility property such that subsequence of orbit of map at some point lies in one of the sets, then there is no subsequence of that orbit which lies in the other set. The proof follows the same idea as that of Theorem~8.2 in \cite{BZ2017}, but we present it in our setting.
\begin{lemma}\label{Lem:Whole_Orbit_lies_if_one_of_its_subseq_lies}
Suppose $\OM \subset \C^d$ is a complete hyperbolic domain and $A, B \subset \OM$ such that $A$ and $B$ have visibility with respect to $(\OM, K_\OM)$. Suppose $f \in \hbox{Hol}(\OM, \OM)$ such that $(f^{n})_{n \geq 1}$ is compactly divergent. Suppose $p \in \OM$ and there exists a subsequence $(n_{k})_{k \geq 1}$ with $n_{k} \rightarrow \infty$ as $k \rightarrow \infty$ such that $(f^{n_{k}}(p))_{k \geq 1} \subset A$. Then there exists $N \in \N$ such that $f^{n}(p) \in \OM\setminus B$ for all $n \geq N$.
\end{lemma}
\begin{proof}
Since $(f^{n})_{n \geq 1}$ is compactly divergent and $(\OM, K_\OM)$ is complete,
\[
\lim_{k \rightarrow \infty} K_\OM(p, f^{n_k}(p)) = \infty.
\]
We pick a subsequence $m_{i} \to \infty$ of $(n_k)_{k\geq 1}$ such that
\begin{equation}\label{Eq:Picking_record_seq-01}
K_\OM(f^{n_k}(p), p) \leq K_\OM(f^{m_{i}}(p), p)
\end{equation}
for all $n_k \leq m_{i}$. We now assume that conclusion of lemma is not true. Then there exists a subsequence $(f^{l_j}(p))_{j\geq 1}$ such that $f^{l_j}(p) \in B$ for all $j\geq 1$.
\smallskip

Pick a sequence $i_j \rightarrow \infty$ with $m_{i_j} > l_j$. Let $\sigma_j :[0,a_j] \rightarrow \Omega$ be geodesic of $(\OM, K_\OM)$ with $\sigma_j(0) = f^{m_{i_j}}(p)$ and $\sigma_j(a_j) = f^{l_j}(q)$. There exists a compact set $K\subset \OM$ that $\sigma_{j}\cap K\neq \emptyset$. Let $t_j \in [0,a_j]$ so that $\sigma_{j}(t_{j}) \in \sigma_{j}\cap K$. Hence, for all $j \in \N$, we have
\begin{align*}
K_\OM(p, \sigma_j(t_j)) \leq \diam(K\cup\{p\})=:C<\infty.
\end{align*}
Since $\sigma_{j}$ is a geodesic of $(\OM,K_\OM)$, using triangle inequality, we have
\begin{align*}
K_\OM(f^{m_{i_j}}(p), f^{l_j}(p))& = K_\OM(f^{m_{i_j}}(p),\sigma_j(t_j)) + K_\OM(\sigma_j(t_j), f^{l_j}(p))\\
        & \geq K_\OM(f^{m_{i_j}}(p), p) + K_\OM(p, f^{l_j}(p))-2C
\end{align*}
On the other hand from \eqref{Eq:Picking_record_seq-01}, we get
\begin{align*}
K_\OM(f^{m_{i_j}}(p), f^{l_j}(p)) \leq K_\OM(f^{m_{i_j}-l_j}(p), p) \leq K_\OM(f^{m_{i_j}}(p), p).
\end{align*}
Hence for all $j \in \N$
\begin{align*}
K_\OM(p, f^{l_j}(p)) \leq 2C.
\end{align*}
Therefore, we get a contradiction. This completes the proof of proposition.
\end{proof}

\section{Auxiliary results}\label{Sec:Prelim_and_aux_res}
In this section, we prove auxiliary results required for the proofs of characterization of domain having the Denjoy-Wolff property that satisfies Condition~(M$\star$) for all boundary points, and outside the totally disconnected set.

\subsection{Relation between Condition~(M$\star$), Condition~($\star$) and Condition $(\star\star)$}\label{Sec:Prop_of_Con(*)_and_relation_bw_Cond(*)_and_Cond(**)} We start with exploring the relation between Condition~(M$\star$), Condition~($\star$) and Condition $(\star\star)$. Firstly, we record some results which shall be used in the proof of Condition~(M$\star$) implies Condition~($\star$).

\begin{result}[Fatou's Theorem]\label{Res:Fatou_thm}
Let $f: \D\to\C$ be a bounded holomorphic function. Then for almost all $\theta$
\[
f^{\star}(e^{i \theta})= \lim_{r \to 1}f(re^{i \theta})
\]
exists.
\end{result}

\begin{result}{\cite[Theorem~2.5]{Collingwood-Lowater_Book}}\label{Thm:Lusin-Privalov_thm}
Let $f:\D \to \C$ be a bounded holomorphic function. If the set of points $e^{i \theta}$ for which $\lim_{r \to 1}f(re^{i \theta})=0$ has positive measure, then $f$ is identically zero.
\end{result}

\begin{proposition}\label{Pro:Cond(M*)_implies_Cond(*)}
Let $\OM \subset \C$ be a hyperbolic domain. If $\OM$ satisfies Condition~(M$\star$), then it satisfies Condition~($\star$).
\end{proposition}
\begin{proof}
Let $x \in\partial\OM$ and a sequence $(z_{n})_{n \geq 1} \subset \OM$ that converges to $x$. Then there exists an injective holomorphic map $\varphi : \D \rightarrow \OM$ obtained by $x$ and $(z_{n})_{n \geq 1}$ as in Condition~(M$\star$). Suppose $U$ is a bi-geodesic region and $(z_{n_{k}})_{k \geq 1}$ is a subsequence as in Condition~(M$\star$). Without loss of generality, we may assume that $\varphi^{-1}(z_{n_{k}}) \in U$ for all $k \in \N$. Therefore, $z_{n_{k}} \in \varphi(U)$ for all $k \geq 1$. Let $\Gamma$ be a component of $\partial\OM$ such that $x \in \Gamma$. We may have the following cases:
\smallskip

\noindent{\bf Case~1.} Assume that $\Gamma$ is bounded.
\smallskip

In this case, we first make the following claim.
\smallskip

\noindent{\bf Claim.} There exists a bi-geodesic region $U^{\prime}$ such that $\varphi(U^{\prime})$ is bounded. 
\smallskip

\noindent{\itshape Proof of Claim.} Suppose that the conclusion of the claim is not true. Therefore, there exists a sequence of bi-geodesic regions $(U_{k})_{k \geq 1}$ such that $U_{k+1} \subset U_{k}$ and $\varphi(U_{k})$ is unbounded for all $k \geq 1$. We may find a sequence $(\zeta_{k})_{k \geq 1} \subset U_{k}$ such that $\delta_{\D}(\zeta_{k}) \to 0$ as $k \to \infty$ and $\varphi(\zeta_{k}) \to \infty$ as $k \to \infty$. We may assume, if necessary pass to a subsequence, that $(\zeta_{k})_{k \geq 1}$ converges to some $\zeta \in \partial\D$. Without loss of generality, we may assume that $\zeta_{k} \in U$ for all $k \geq 1$. Note that $\varphi^{-1}(z_{n_{k}}) \in U$ for all $k \geq 1$. Suppose $\sigma_{k}:[0,a_{k}] \to \D$ is a continuous curve that joins $\zeta_{k}$ and $\varphi^{-1}(z_{n_{k}})$ such that $\sigma_{k} \subset U$. Note that, if necessary by passing to a subsequence, we may assume that $(\sigma_{k})_{k \geq 1}$ converges in the Hausdorff distance to an arc joining $1$ and $\zeta$. Moreover, $\varphi(\sigma_{k})$ is compact and connected for all $k\geq 1$. There exists a set $S$ such that it is compact in $\C_\infty$ and such that, if necessary by passing to a subsequence, $(\varphi(\sigma_{k}))_{k\geq 1}$ converges in the Hausdorff distance to $S$. Note that, by construction, $x, \infty \in S$. Note that $S$ is connected, $\Gamma$ is connected and $\Gamma\cap S\not=\emptyset$. Hence $S \subset \Gamma$. This contradicts the fact that $\Gamma$ is bounded because it contains unbounded set $S$.
\hfill $\blacktriangleleft$
\smallskip

We now only need to show that for every two bi-geodesic regions $U_{1}$ and $U_{2}$ with $U_{1} \subsetneq U_{2} \subseteq U^{\prime}$, $K_\OM(\varphi(U_{1}), \varphi({\D \setminus U_{2}})) > 0$. 
Suppose that, to get a contradiction, $K_\OM(\varphi(U_{1}), \varphi(\D \setminus U_{2})) = 0$. Then, there exist sequences $(w_{n})_{n \geq 1} \subset \varphi(U_{1})$ and $(w_{n}^{\prime})_{n \geq 1} \subset \varphi(\D)\setminus\varphi(U_{2})$ such that $K_\OM(w_{n}, w^{\prime}_{n}) \to 0$ as $n \to \infty$. Therefore, we may assume, up to a subsequences, that $(w_{n})_{n \geq 1}$ and $(w_{n}^{\prime})_{n \geq 1}$ converge to $y \in \partial \OM$. Suppose $\gamma_{n}:[0,a_{n}] \to \OM$ is geodesic of $(\OM, K_\OM)$ joining $w_{n}$ and $w_{n}^{\prime}$. Since each $\gamma_n$ is a geodesic of $(\OM,K_\OM)$, so it follows that for any $\wt{w}_{n}, \wt{w}_{n}^{\prime} \in \gamma_{n}, K_\OM(\wt{w}_{n}, \wt{w}_{n}^{\prime}) \to 0$ as $n \to \infty$. Also, for each $n \in \N$, let
\begin{align*}
s_{n}&:= \sup\{t\in[0,a_{n}]:\gamma_{n}(s)\in \varphi(U_{1}) \forall s\in [0,t]\} \,, \text{and}\\
        t_{n}&:=\inf\{t\geq s_{n}:\gamma_{n}(t) \not \in \varphi(U_{2})\}.
\end{align*}
Note that, by construction, $\gamma_{n}(s_{n}) \in \partial\varphi(U_{1})\cap \OM, \gamma_{n}(t_{n}) \in \partial\varphi(U_{2})\cap\OM$ and $\left.\gamma_{n}\right|_{[s_{n},t_{n}]} \subset \clos{\varphi(U_{2})}\setminus \varphi(U_{1})$. Let $\wt{\gamma}_{n}:= \left.\gamma_{n}\right|_{[s_{n},t_{n}]}$. Clearly, $(\wt{\gamma}_{n})_{n \geq 1}$ converges in the Hausdorff distance to a point $y$. Moreover, $\wt{\sigma}_{n}:=\varphi^{-1}(\wt{\gamma}_{n})$ also converges to $\wt{\sigma}$ in the Hausdorff distance where $\wt{\sigma}$ is a component of $\clos{U_{2}}\cap\partial\D \setminus \clos{U_{1}}\cap \partial\D$. Note that $U^{\prime}$ is a Jordan domain, so it follows from the Riemann mapping theorem and the Carath\`eodory extension theorem that there exists a homeomorphism $\wt{f}: \clos{\D} \to \clos{U^{\prime}}$ such that $f:=\left.\wt{f}\right|_\D: \D \to U^{\prime}$ is a biholomorphism. Let $\sigma^{\prime}_{n}:= \wt{f}^{-1}\circ \wt{\sigma}_{n}$. Since $\wt{f}$ is a homeomorphism, therefore, $(\sigma^{\prime}_{n})_{n \geq 1}$ converges in the Hausdorff distance to an arc $\sigma^{\prime}$ of boundary of the unit disc. Note that, by Claim above, $g:=\varphi\circ f$ is bounded holomorphic function on $\D$, so it follows from Result~\ref{Res:Fatou_thm} that radial limit of $g$ exists almost everywhere. In particular, radial limit of $g$ exists for every $\xi\in int(\sigma^{\prime})$ where $int(\sigma^\prime)$ denotes the interior of $\sigma^\prime$. Suppose $(\zeta_{n})_{n \geq 1}\ \subset \sigma^{\prime}_{n}\cap L_{\xi}$ that converges to $\xi\in int(\sigma^{\prime})$ where $L_{\xi}$ is a radial line joining $0$ and $\xi$. Since $(\wt{\gamma}_{n})_{n \geq 1}$ converges in the Hausdorff distance to a point $y$, $(g(\zeta_{n}))_{n \geq 1}$ converges to $y$. Therefore, radial limit of $g$ at almost every $\xi \in int(\sigma^{\prime})$ is $y$. Consider a function $F: \D \to \C$ given by $F(\zeta)= g(\zeta) - y$. By the Claim, $F$ is bounded holomorphic function. Moreover, radial limit of $F$ at almost every $\xi\in int(\sigma^{\prime})$ is $0$. Therefore, by Result~\ref{Thm:Lusin-Privalov_thm}, $F \equiv 0$. That is, $g(\zeta) = y$ for all $\zeta \in \D$. This contradicts injectivity of the map $g$.
\smallskip

\noindent{\bf Case~2.} Assume that $\Gamma$ is unbounded.
\smallskip

Since $U$ is a bounded simply connected domain such that its boundary is the Jordan curve, so it follows from the Riemann mapping theorem and the Carath\`eodory extension theorem that there exists a homeomorphism $\wt{g}: \clos{\D} \to \clos{U}$ such that $g\defeq\left.\wt{g}\right|_{\D}: \D \to D$ is a biholomorphism. Let $p \in \OM\setminus\varphi(U)$ and $\OM_{p}\defeq \OM\setminus\{p\}$. Consider a map $T(z) \defeq 1/(z-p)$. Moreover, it is easy to see that $T(x) \in \partial T(\OM_{p})$ and $(T(z_{n}))_{n \geq 1}\subset T(\OM_{p})$ converges to $T(x)$. Note that $T(\OM_p)$ satisfies Condition~(M$\star$) at $T(x)$. Therefore, we take an injective map $\varphi_{p}:\D \to T(\OM_{p})$ given by $\varphi_{p}:=T\circ\varphi$ obtained by $T(x)$ and $(T(z_{n}))_{n \geq 1}$ as in Condition~(M$\star$). Note that $T(\Gamma)$ is bounded. By Case~1, we get that there exists a bi-geodesic region $U_p$ such that $\varphi_p(U_p)$ is bounded and
\begin{equation}\label{ineq:separation_under_analytic_disc-1}
K_{T(\OM_{p})}(\varphi_{p}(U_{1}), \varphi_{p}(\D)\setminus \varphi_p(U_{2}))>0
\end{equation}
for every two bi-geodesic regions $U_{1}$ and $U_{2}$ with $U_{1} \subsetneq U_{2}\subsetneq U_{p}$. We want to show that $K_\OM(\varphi(U_{1}), \varphi(\D)\setminus \varphi(U_{2}))>0$. Note that 
\begin{equation}\label{ineq:decreasing_prop_of_Kobayashi_dist-1}
K_{T(\OM_{p})}(\varphi_{p}(U_{1}), \varphi_p(\D)\setminus \varphi(U_{2})) \leq K_{\OM_{p}}(\varphi(U_{1}), \varphi(\D)\setminus \varphi(U_{2})).
\end{equation}
We choose $r>0$ such that $\clos{B(p,3r)}\cap \varphi(U_0)=\emptyset$. Let $W\defeq \OM_p\setminus \clos{B(p,3r)}, W_1\defeq \OM_p\setminus \clos{B(p,2r)}$ and $W_2\defeq\OM_p\setminus \clos{B(p,r)}$. It easy to see that $W$ and $\OM\setminus W_1$, and $W_1$ and $\OM\setminus \OM_p$ have weak visibility property with respect to $(\OM,K_\OM)$. Therefore, by Theorem~\ref{Thm:Global_Loc_Res}, there exists a constant $C_0>0$ such that for all $z,w\in W$
\[
K_{\OM_p}(z,w)\leq K_\OM(z,w)+ C_0.
\]
Note that it is clear that $K_\OM(W_1,\OM\setminus W_2)>0$. Thus, by Theorem~\ref{Thm:Mult_localization}, there exists a constant $C\geq 1$ such that 
\begin{equation}\label{Ineq:Mult_loc_ineq}
K_{\OM_{p}}(\varphi(U_{1}), \varphi(\D\setminus U_{2})) \leq CK_\OM(\varphi(U_{1}), \varphi(\D)\setminus \varphi(U_{2})).
\end{equation}
Then, using~\eqref{Ineq:Mult_loc_ineq} and~\eqref{ineq:decreasing_prop_of_Kobayashi_dist-1}, we get
\[
(1/C)K_{T(\OM_{p})}(\varphi_{p}(U_{1}), \varphi_p(\D)\setminus \varphi(U_{2})) \leq K_\OM(\varphi(U_{1}), \varphi(\D)\setminus \varphi(U_{2})).
\]
Therefore, using previous inequality and~\eqref{ineq:separation_under_analytic_disc-1}, we get $K_\OM(\varphi(U_{1}), \varphi(\D)\setminus \varphi(U_{2}))>0$.

Hence, we get a bi-geodesic region $\wt U$ such that $K_\OM(\varphi(U_{1}), \varphi(\D)\setminus \varphi(U_{2}))>0$ for every two bi-geodesic regions $U_{1}$ and $U_{2}$ with $U_{1} \subsetneq U_{2}\subsetneq \wt U$. This completes the proof.
\end{proof}

We now discuss the relation between Condition~(M$\star$) and Condition $(\star\star)$.
\begin{lemma}\label{Lem:Contions_double_star_implies_mild_star}
Let $\OM \subset \C$ be a hyperbolic domain. If $\OM$ satisfies Condition $(\star\star)$ then it satisfies Condition $M(\star)$.    
\end{lemma}
\begin{proof}Suppose $x \in \partial\OM$ and $(z_{n})_{\geq 1}\subset \OM$ is a sequence that converges to $x$. Suppose $\Gamma \subset \partial\OM$ is a component of $\bdy\OM$ such that $x \in \Gamma$. According to Condition $(\star\star)$, there exists an injective holomorphic map $f:A(0,r,1)\to\OM$ such that for every $p\in\Gamma$ and for every sequence $(w_n)_{n\geq 1}$ in $\OM$ that converges to $p$, we have $|f^{-1}(w_n)|\to \bdy\D$ as $n\to\infty$. Thus, without loss of generality, we may assume that there exists a subsequence $(z_{n_k})_{k\geq 1}$ such that $f^{-1}(z_{n_k})\to 1$ as $k\to \infty$. Fix $\epsilon>0$ such that $\clos{B(0,r)} \cap \clos{B(1,\epsilon)} = \emptyset$. Let $D\defeq A(0,r,1) \cap B(1,\epsilon)$. Without loss of generality, we may assume that $f^{-1}(z_{n_k})\in D$ for all $k \geq 1$. Note that $D \subset \C$ is a bounded simply connected domain such that its boundary is the Jordan curve, so it follows from the Riemann mapping theorem and the Carath\`eodory extension theorem that there exists a homeomorphism $\wt{g}: \clos{\D} \to \clos{D}$ such that $g\defeq \left.\wt{g}\right|_{\D}: \D \to D $ is a biholomorphism. Without loss of generality, we may assume that $\wt{g}(1) = 1$. Let $\varphi: \D \to \OM$ defined by $\varphi\defeq f \circ g$. Note that $\varphi$ is an injective holomorphic map. Moreover, $(\varphi^{-1}(z_{n_{k}}))_{k \geq 1}$ converges to $1$. Since $\wt g$ is a continuous map, so we may choose a bi-geodesic region $U$ such that $\wt{g}(\clos{U}\cap \partial \D) \subseteq \partial \D\cap \clos{B(1,\epsilon)}$.  Let $(\zeta_{n})_{n \geq 1} \subset U$ such that $\delta_{\D}(\zeta_{n}) \to 0$ as $n \to \infty$. Note that $\distance(\wt{g}(\zeta_{n}), \partial\D\cap\clos{B(1,\epsilon)}) \to 0$ as $n \to \infty$ where $\distance(\wt{g}(\zeta_{n}), \partial\D\cap\clos{B(1,\epsilon)})$ denotes the Euclidean distance of $\wt{g}(\zeta_{n})$ to the set $\partial\D\cap\clos{B(1,\epsilon)}$. Since either $\delta_\OM(f(\zeta_n))\to 0$ as $n\to\infty$ or $|f(\zeta_n)|\to \infty$ as $n\to\infty$, so it follows that either $\delta_{\OM}(\varphi(\zeta_{n})) \to 0$ as $n \to \infty$ or $|\varphi(\zeta_{n}))| \to \infty$ as $n \to \infty$. This completes the proof of lemma. 
\end{proof}
We recall the following topological result.
\begin{result}[{Zoretti's theorem; \cite[Corollary~VI.3.11]{Whyburn:topology}}]\label{Res:Zoretti_Thm}
Let $E \subset \C$ be a compact set and $S$ be a component of $E$. For each $\epsilon > 0$ there exists a Jordan region $U$ such that $S \subset U$, $\partial U \cap E = \emptyset$, and $\partial U \subset N_{\C}(S,\epsilon)$.
\end{result}

\begin{lemma}\label{Lem:Equivance_of_Con(*)andCon(**)_for_bdd_compo}
Let $\OM \subset \C$ be a hyperbolic domain such that its every boundary component is bounded. If $\OM$ satisfies Condition~(M$\star$) then it satisfies Condition $(\star\star)$.   
\end{lemma}
\begin{proof}
Let $\Gamma$ be a connected component of $\partial\OM$ and let $\OM_{\Gamma}$ be simply connected domain obtained by filling all the connected components of $\C\setminus\OM$ except which contains $\Gamma_1$. Then, by the Riemann mapping theorem, there exists a biholomorphism $f:\D\to \OM_\Gamma$. Note that $\bdy\OM_\Gamma=\Gamma$. Let $x \in \Gamma$ and $(z_{n})_{n \geq 1}\subset \OM$ be a sequence that converges to $x$. Therefore, we get an injective holomorphic map $\varphi_x: \D \to \OM$ obtained by $x$ and $(z_{n})_{n \geq 1}$ as in Condition~(M$\star$). Let $D_x\defeq \varphi_x(\D)$. Then the set $\{\clos{D}_x:x \in \Gamma\}$ is a cover of $\Gamma$. Since $\Gamma$ is bounded, so there exists $x_1,\ldots,x_n\in \Gamma$ such that $\Gamma\subseteq \cup_{i=1}^{n}\clos{D}_{x_i}$. Choose $\delta>0$ such that for all $z \in \OM\setminus\displaystyle\bigcup_{i=1}^{n}\clos{D}_{y_{i}}$, we have $\distance(z, \Gamma)>\delta$. We choose $\epsilon<\delta$ sufficiently small such that $N_\C(\Gamma,\epsilon)\cap\OM_\Gamma\subset\cup_{i=1}^{n}D_i$. By Result~\ref{Res:Zoretti_Thm}, there exists a Jordan region $V$ such that $\Gamma \subset V$, $\partial V\cap \Gamma = \emptyset$ and $\partial V \subset N_{\C}(\Gamma,\epsilon)$. We shall show that $V\cap\OM$ is doubly connected. Note that $\partial V \subset \cup_{i=1}^{n}D_{x_{i}} \subset \OM$. Note that $\partial V \subset \cup_{i=1}^{n}D_{x_{i}} \subset \OM$. Therefore, $\C_{\infty}\setminus V\cap\OM$ has only one component in $V$ which is $\Gamma$. By the Jordan curve theorem, $\partial V$ divides $\C_{\infty}$ in two connected components both having common boundary $\partial V$. Hence, $\C_{\infty}\setminus V\cap\OM$ has only two components. This implies that $V\cap \OM$ is doubly connected. Therefore, $f^{-1}(V\cap\OM)$ is doubly connected. Without loss of generality, we may assume that $f^{-1}(V\cap\OM)=A(0,r,1)$ for some $r>0$. Therefore, we get an injective holomorphic map $f_0\defeq\left.f\right|_{A(0,r,1)}:A(0,r,1)\to\OM$ such that $f_0(A(0,R,1))=V\cap\OM$ and such that for every sequence $(\zeta_n)_{n\geq 1}\in A(0,r,1)$ with $|\zeta_n|\to 1$ as $n\to\infty$, we have either $\delta_\OM(f(\zeta_n))\to 0$ as $n\to\infty$ or  $|f(\zeta_n)|\to \infty$ as $n\to\infty$. This completes the proof.
\end{proof}

\begin{remark}
By establishing a result analogous to Result~\ref{Res:Zoretti_Thm}, we may prove Lemma~\ref{Lem:Equivance_of_Con(*)andCon(**)_for_bdd_compo} without assuming that the boundary components are bounded. However, we omit such proof, as we only need to establish that Condition $(\star\star)$ implies Condition~(M$\star$) for our purposes.
\end{remark}

\subsection{Embedded disc preserves intrinsic visibility property}\label{Sub:Embedded_disc_preserve_vis}
It is easy to see that if $U_1$ and $U_2$ are two bi-geodesic regions with $U_1\subsetneq U_2$, then $U_1$ and $\D\setminus U_2$ have visibility property with respect to $(\D, K_\D)$. In this subsection, we prove that visibility property of the image under the attached analytic disc of these two sets is preserved.
\begin{lemma}\label{Lem:Vis_of_strips_wrt_phi(D)}
Suppose $\OM \subset \C$ is a hyperbolic domain and $\varphi : \D \rightarrow \OM$ is an injective holomorphic map. If $U_{1}$ and $U_{2}$ are two bi-geodesic regions with $U_{1} \subsetneq U_{2}$, then $\varphi(U_{1})$ and $\varphi(\D) \setminus\varphi( U_{2})$ have weak visibility property with respect to $(\varphi(\D), K_{\varphi(\D)})$.
\end{lemma}
\begin{proof}
Assume to the contrary that $\varphi(U_{1})$ and $\varphi(\D) \setminus\varphi( U_{2})$ does not have weak visibility property with respect to $(\varphi(\D), K_{\varphi(\D)})$. Then there exist $\k>0$ and a sequence $\gamma_{n}: [0,a_{n}] \rightarrow \varphi(\D)$ of $(1,\k)$-almost-geodesic of $(\varphi(\D),K_{\varphi(\D)})$ with $\gamma_{n}(0) \in \varphi(U_{1}) $ and $\gamma_{n}(a_{n}) \in \varphi(\D) \setminus\varphi( U_{2})$ such that $(\gamma_{n})_{n \geq 1}$ eventually avoids all compact subset of $\varphi(\D)$. Thus, $\varphi^{-1}\circ \gamma_{n}$ is $(1,\k)$-almost-geodesic of $(\D, K_{\D})$ with $\varphi^{-1}(\gamma_{n}(0)) \in U_{1}$ and $\varphi^{-1}( \gamma_{n}(a_{n})) \in \D \setminus U_{2}$ such that $(\varphi^{-1}\circ \gamma_{n})_{n\geq1}$ evetually avoids every compact subset of $\D$. This contradicts the weak visibility of $U_{1}$ and $\D\setminus U_{2}$ with respect to $(\D, K_{\D})$. This completes the proof of lemma.
\end{proof}

\begin{lemma}\label{Cor:Vis_of_strips_wrt_Omega}
Suppose $\OM \subset \C$ is a hyperbolic domain and $\varphi : \D \rightarrow \OM$ is an injective holomorphic map. Fix a bi-geodesic region $U$ such that for every two bi-geodesic regions $U_1$ and $U_2$ with $U_1\subsetneq U_2\subseteq U$, we have $K_\OM(\varphi(U_1), \varphi(\D) \setminus \varphi(U_2)) > 0$. If $U_0\subsetneq U_0^\prime\subseteq U$ then $\varphi(U_0)$ and $\varphi(\D) \setminus \varphi(U_0^\prime)$ have weak visibility property with respect to $(\OM, K_\OM)$.
\end{lemma}
\begin{proof}
Let $U^\prime, U^\star$ and $U_W$ be bi-geodesic regions such that $U^\prime \subsetneq U^\star\subsetneq U_W\subsetneq U$ and such that $U_0\subsetneq U^\prime$. Therefore, $\varphi(U^\prime)\subsetneq \varphi(U^\star)\subsetneq \varphi(U_W)$ and $\varphi(U_0)\subsetneq \varphi(U^\prime)$. By Lemma~\ref{Lem:Vis_of_strips_wrt_phi(D)}, $\varphi(U^\prime)$ and $\varphi(\D)\setminus\varphi(U^\star)$ and also $\varphi(U^\star)$ and $\varphi(\D)\setminus\varphi(U_W)$ have weak visibility property with respect to $(\varphi(\D),K_{\varphi(\D)})$. Since for every two bi-geodesic regions $U_1$ and $U_2$ with $U_1\subsetneq U_2\subseteq U$ we have $K_\OM(\varphi(U_1), \varphi(\D) \setminus \varphi(U_2)) > 0$, so it follows that $K_\OM(\varphi(U_W),\OM\setminus\varphi(\D))>0$. Moreover, by Lemma~\ref{Lem:Vis_of_strips_wrt_phi(D)}, $\varphi(U_0)$ and $\varphi(\D)\setminus\varphi(U_0^\prime)$ have weak visibility property with respect to $(\varphi(\D),K_{\varphi(\D)})$. Note that all the hypothesis's of Theorem~\ref{Thm:Local_to_global_vis} are satisfied. Therefore, by Theorem~\ref{Thm:Local_to_global_vis}, $\varphi(U_0)$ and $\varphi(\D)\setminus\varphi(U_0^\prime)$ have weak visibility property with respect to $(\OM,K_\OM)$. This completes the proof of the Corollary.
\end{proof}

\begin{lemma}\label{Lem:Point_outside_U}
Suppose $\OM \subset \C$ is a hyperbolic domain and $\varphi : \D \rightarrow \OM$ is an injective holomorphic map. Suppose $U$ is a bi-geodesic region such that for every sequence $(\zeta_{n})_{n \geq 1}$ in $U$ such that $\delta_{\D}(\zeta_{n}) \rightarrow 0$ as $n \rightarrow \infty$, either $\delta_{\OM}(\varphi(\zeta_{n})) \rightarrow 0$ or $|\varphi(\zeta_{n})| \rightarrow \infty$ as $n \rightarrow \infty$. Suppose $U_1$ and $U_2$ are bi-geodesic regions with $U_1\subsetneq U_2 \subseteq U$ and $\gamma : [0, a] \rightarrow \OM$ is geodesic of $(\OM, K_\OM)$ such that $\gamma(0) \in \varphi(U_1)$ and $\gamma(a) \in \OM\setminus\varphi(U_2)$. Then there exists a point $\zeta \in \D\setminus U_2$ such that $\varphi^{-1}\circ \gamma(s) = \zeta$ for some $s \in (0, a^{\prime})$ where 
\[
a^{\prime} := \sup \{s\in [0, a] : \left.\gamma\right|_{[0,s]} \subseteq \varphi(\D)\}.
\]
\end{lemma}
\begin{proof}
Let $\wh\gamma:=\left.\gamma\right|_{[0,a^\prime]}$.
Assume, to get a contradiction, that $\eta := \varphi^{-1}\circ\wh\gamma \subset U_2$. Let $(s_{k})_{k \geq 1} \subset (0, a^{\prime})$ be an increasing sequence such that $s_{k} \rightarrow a^{\prime}$ as $k \rightarrow \infty$. Without loss of generality, if necessary passing to a subsequence, we assume $(\eta(s_{k}))_{k \geq 1}$ converges to some $\xi \in \clos{U_2}$. We now make the following claim.
\smallskip

\noindent {\bf Claim.} $\xi \in \clos{U_2} \setminus \partial{\D}$.
\smallskip

\noindent {\itshape Proof of Claim.} Assume that $\xi \in \clos{U_2} \cap \partial{\D}$. Note that $\varphi\circ\eta=\wh\gamma$. By the continuity of $\wh\gamma, \wh\gamma(s_{k}) \rightarrow \wh\gamma(a^{\prime})$ as $k \rightarrow \infty$. Since $\delta_{\OM}$ is a continuous function, $\delta_{\OM}(\wh\gamma(s_{k})) \rightarrow \delta_{\OM}(\wh\gamma(a^{\prime}))$ as $k \rightarrow \infty$. Moreover, by assumption, $\delta_{\OM}(\wh\gamma(s_{k})) \rightarrow 0$ or $|\wh\gamma(s_{k})| \rightarrow \infty$ as $k \rightarrow \infty$. If $\delta_{\OM}(\wh\gamma(s_{k})) \rightarrow 0$ as $k \rightarrow \infty$, then $\delta_{\OM}(\wh\gamma(a^{\prime})) = 0$ which is not possible as $\OM$ is an open set and $\wh\gamma(a^{\prime}) \in \OM$, so $\delta_{\OM}(\wh\gamma(a^{\prime})) > 0$. Moreover, $|\wh\gamma(s_{k})| \rightarrow \infty$ as $k \rightarrow \infty$ is not possible because $\wh\gamma$ is compact. Therefore, we get contradiction in both the cases. This completes the proof of claim.
\hfill $\blacktriangleleft$
\smallskip

By continuity of $\varphi$, $\wh\gamma(s_k) \to \varphi(\xi)$. Now, by continuity of $\wh\gamma$, $\wh\gamma(a^\prime)=\varphi(\xi)$. Let $r>0$ such that $B(\varphi(\xi),r)\subset\varphi(\D)$. By continuity of $\gamma$, there exists $\delta>0$ such that $|\gamma(s)-\gamma(a^\prime)|<r$ whenever $|s-a^\prime|<\delta$. Choose $s^\prime>a^\prime$ with $|s-a^\prime|<\delta$. Therefore, $|\gamma(s^\prime)-\gamma(a^\prime)|<r$. That is, $\left.\gamma\right|_{[a^\prime,s^\prime]} \subset B(\varphi(\xi),r)$. Therefore, $\left.\gamma\right|_{[a^\prime,s^\prime]} \subset\varphi(\D)$. This gives a contradiction to the definition of $a^\prime$. This completes the proof.
\end{proof}

\begin{lemma}\label{Lem:Vis_of_image_of_bigeodesic_region_and_image_of_unit_disc}
Suppose $\OM \subset \C$ is a hyperbolic domain. Suppose that $\varphi : \D \rightarrow \OM$ is an injective holomorphic map. Let $U$ be a bi-geodesic region such that following holds:
\begin{enumerate}
\item[(i)] if $(\zeta_{n})_{n \geq 1}\subset U$ is a sequence such that $\delta_{\D}(\zeta_{n}) \rightarrow 0$ as $n \rightarrow \infty$, then either $\delta_{\OM}(\varphi(\zeta_{n})) \rightarrow 0$ or $|\varphi(\zeta_{n})| \rightarrow \infty$ as $n \rightarrow \infty$, and
\item[(ii)] if $U_1$ and $U_2$ are two bi-geodesic regions with $U_1\subsetneq U_2\subseteq U$ then $K_\OM(\varphi(U_1), \varphi(\D) \setminus \varphi(U_2)) > 0$.
\end{enumerate}
If $U^\prime$ is bi-geodesic region with $U^\prime\subsetneq U$, then $\varphi(U^\prime)$ and $\OM\setminus \varphi(\D)$ have weak visibility property with respect to $(\OM,K_\OM)$.
\end{lemma}
\begin{proof}
Assume to the contrary that $\varphi(U^\prime)$ and $\OM \setminus \varphi(\D)$ does not have weak visibility property with respect to $(\OM, K_\OM)$. Then there exists a sequence $ \gamma_{n}: [0,a_{n}] \rightarrow \OM$ of geodesic of $(\OM,K_\OM)$ with $\gamma_{n}(0) \in \varphi(U^\prime) $ and $\gamma_{n}(a_{n}) \in \OM\setminus \varphi(\D)$ such that $(\gamma_n)_{n\geq1}$ eventually avoids every compact subset of $\OM$. Fix, for each $n\in \N$,
\[
t_{n}:= \sup \{s\in[0, a_{n}]: \gamma_{n}(t) \in \varphi(\D) \forall t \in [0,s]\}
\]
By construction, for all $n \in \N, \left.\gamma_n\right|_{[0,t_n)}\subset \varphi(\D)$. Thus, for all $n \in \N, \sigma_{n}:=\varphi^{-1}\circ\left.\gamma_n\right|_{[0,t_{n})} \subset \D$. Let $U_0^\prime$ be a bi-geodesic region such that $U^\prime\subsetneq U_0^\prime\subsetneq U$. By Lemma~\ref{Lem:Point_outside_U}, there exists $\zeta_{n} \in \D\setminus U_0^\prime$ such that $\sigma_{n}(s_{n}) = \zeta_{n}$ for some $s_n \in (0, t_{n})$. Note that $\gamma_{n}(0) \in \varphi(U^\prime)$ and $\gamma_{n}(s_{n}) \in \varphi(\D)\setminus\varphi(U_0^\prime)$. Therefore, $(\left.\gamma_n\right|_{[0,s_n]})_{n\geq 1}$ eventually avoids every compact subset of $\varphi(\D)$. By Lemma~\ref{Cor:Vis_of_strips_wrt_Omega}, we get a contradictions to the fact that $\varphi(U^\prime)$ and $\varphi(\D)\setminus\varphi(U_0^\prime)$ weak visibility property of $(\varphi(\D),K_{\varphi(\D)})$. This completes the proof of lemma. 
\end{proof}

The next lemma shows that embedded disc preserves intrinsic visibility property.
\begin{lemma}\label{Lem:Visi_of_images_of_bi-geodesic_reg_inside_U}
Suppose $\OM \subset \C$ is a hyperbolic domain and $\varphi: \D \rightarrow \OM$ is an injective holomorphic map. Let $U$ be a bi-geodesic region such that following holds:
\begin{enumerate}
\item[(i)] if $(\zeta_{n})_{n \geq 1}\subset U$ is a sequence such that $\delta_{\D}(\zeta_{n}) \rightarrow 0$ as $n \rightarrow \infty$, then either $\delta_{\OM}(\varphi(\zeta_{n})) \rightarrow 0$ or $|\varphi(\zeta_{n})| \rightarrow \infty$ as $n \rightarrow \infty$, and
\item[(ii)] if $U_1$ and $U_2$ are two bi-geodesic regions with $U_1\subsetneq U_2\subseteq U$ then $K_\OM(\varphi(U_1), \varphi(\D) \setminus \varphi(U_2)) > 0$.
\end{enumerate}
If $U^\prime$ and $U^{\prime\prime}$ are bi-geodesic regions with $U^\prime \subsetneq U^{\prime\prime} \subseteq U$, then $\varphi(U^\prime)$ and $\OM \setminus \varphi(U^{\prime \prime})$ have visibility property with respect to $(\OM, K_\OM)$.
\end{lemma}
\begin{proof}
Assume to the contrary that $\varphi(U^\prime)$ and $\OM \setminus \varphi(U^{\prime \prime})$ does not have visibility property with respect to $(\OM, K_\OM)$. Then there exists a sequence $ \gamma_{n}: [0,a_n] \to(\OM, K_\OM)$ of geodesic of $(\OM,K_\OM)$ with $\gamma_{n}(0) \in \varphi(U^\prime)$ and $\gamma_{n}(a_n) \in \OM \setminus(U^{\prime\prime})$ such that $(\gamma_n)_{n \geq 1}$ eventually avoids every compact subset of $\OM$. Consider
\[
t_n:= \sup \{s\in[0, a_{n}]: \gamma_{n}(t) \in \varphi(\D) \forall t \in [0,s]\}.
\]
By definition, for all $n \in \N, \left.\gamma_n\right|_{[0,t_{n})} \subset \varphi(\D)$. Thus, for all $n \in \N, \sigma_{n}:=\varphi^{-1}\circ\left.\gamma_n\right|_{[0,t_{n})} \subset \D$. Let $U_0$ be a bi-geodesic region such that $ U^\prime \subsetneq U_0\subsetneq U^{\prime \prime}$. By Lemma~\ref{Lem:Point_outside_U}, there exists $\zeta_{n} \in \D\setminus U_0$ such that $\sigma_{n}(s_{n}) = \zeta_{n}$ for some $s_n \in (0, t_{n})$. Note that $\gamma_{n}(0) \in \varphi(U^\prime)$ and $\gamma_{n}(s_{n}) \in \varphi(\D)\setminus\varphi(U_0)$. Moreover, $(\left.\gamma_n\right|_{[0,s_n]})_{n\geq 1}$ eventually avoids every compact subset of $\OM$. Therefore, $\varphi(U^\prime)$ and $\varphi(\D)\setminus \varphi(U_0)$ does not have visibility property with respect to $(\OM, K_\OM)$ which contradicts Lemma~\ref{Cor:Vis_of_strips_wrt_Omega}. This completes the proof.
\end{proof}

We shall end this subsection by proving one additional lemma which states that intrinsic visibility property allows to truncate a geodesic ray inside some set. 
\begin{lemma}[Truncation lemma]\label{Lem:Part_of_Geo_lies_in}
Suppose $\OM \subset \C^d$ is a complete hyperbolic domain and $A, B \subset \OM$ such that $A$ and $B$ have visibility with respect to $(\OM, K_\OM)$. Suppose that $\gamma : [0, \infty) \rightarrow \OM$ is a geodesic ray of $(\OM, K_\OM)$ such that there exists a sequence $(t_{n})_{n \geq 1}$ with $t_{n} \rightarrow \infty$ as $n \rightarrow \infty$ such that $\gamma(t_{n}) \in A$ for all $n \in \N$. Then there exists $T > 0$ such that $\left.\gamma\right|_{[T, \infty)} \subset \OM \setminus B$. 
\end{lemma}

\begin{proof}
Assume that the conclusion of the lemma is not true. Suppose that there exists a sequence $(s_{n})_{n \geq 1}$ with $s_{n} \rightarrow \infty$ such that $\gamma(s_{n}) \in B$ for all $n \in \N$. Without loss of generality we assume that $s_{n} \leq t_{n}$ for all $n \in \N$. Since $A$ and $B$ have visibility with respect to $(\OM, K_\OM)$, there exists a compact set $K \subset \OM$ such that $\left.\gamma\right|_{[s_{n}, t_{n}]}\cap K \neq \emptyset$. Let $o_{n} \in \left.\gamma\right|_{[s_{n}, t_{n}]}\cap K \neq \emptyset$. Note that $o_{n} =\gamma(a_{n})$ for some $a_{n} \in [s_{n}, t_{n}]$. Hence, for all $n \in \N$
\[
K_\OM(\gamma(0), \gamma(a_{n})) \leq \diam{(K \cup \{\gamma(0)\})}.
\]
Therefore, $\limsup_{n \rightarrow \infty}{K_\OM(\gamma(0), \gamma(a_{n}))} < \infty$. Note that $\gamma$ is a geodesic ray of $(\OM, K_\OM)$, so $K_\OM(\gamma(0), \gamma(a_{n})) = a_{n}$ for all $n \in \N$. Hence, $\limsup_{n \rightarrow \infty}{K_\OM(\gamma(0), \gamma(a_{n}))} = \infty$. This contradiction concludes the proof.
\end{proof}

\subsection{Properties of Condition~(M$\star$)}\label{Sub:Prop_of_M*}
In this subsection, we shall prove few properties of domains satisfying Condition~(M$\star$) at every point and as well as for domains satisfying Condition~(M$\star$) outside totally disconnected set.
\begin{proposition}\label{Prop:Vis_of_all_req_sets_which_works_everywhere}
Let $\OM\subset \C$ be a hyperbolic domain and let $(z_n)_{n \geq 1}$ be a sequence that converges to some $x\in\bdy\OM$. Suppose that $\OM$ satisfies Condition~(M$\star$) at $x$. Let $\varphi:\D \to \OM$ be an injective holomorphic map, $U$ be a bi-geodesic regions and $(z_{n_k})_{k\geq 1}$ be a subsequence as in Condition~(M$\star$) obtained by $(z_n)_{n \geq 1}$ and $x\in\bdy\OM$. Then there exist bi-geodesic regions $U_1\subsetneq U_2\subsetneq U_3 \subsetneq U$ such that $\varphi(U_1)$ and $\OM\setminus\varphi(U_2)$, $\varphi(U_2)$ and $\OM\setminus\varphi(U_3)$, and $\varphi(U_3)$ and $\OM\setminus\varphi(\D)$ have visibility property with respect to $(\OM, K_\OM)$.
\end{proposition}
\begin{proof}
Let $U_1$ and $U_2$ be two bi-geodesic regions such that $U_1\subsetneq U_2\subsetneq U$. Since $\OM$ satisfies Condition~(M$\star$), every hypotheses of Lemma~\ref{Lem:Visi_of_images_of_bi-geodesic_reg_inside_U} are satisfied. Therefore, by Lemma~\ref{Lem:Visi_of_images_of_bi-geodesic_reg_inside_U}, $\varphi(U_{1})$ and $\OM \setminus \varphi(U_{2})$ have visibility property with respect to $(\OM, K_\OM)$. Let $U_3$ be a bi-geodesic region such that $U_2\subsetneq U_3\subsetneq U$. Then $\varphi(U_2)$ and $\OM\setminus\varphi(U_3)$ have visibility property with respect to $(\OM, K_\OM)$. By Lemma~\ref{Lem:Vis_of_image_of_bigeodesic_region_and_image_of_unit_disc}, $\varphi(U_3)$ and $\OM\setminus\varphi(\D)$ have visibility property with respect to $(\OM,K_\OM)$. This completes the proof of proposition.
\end{proof}

\begin{remark}
Although, in previous subsection, we have proved much stronger result, see Lemma~\ref{Lem:Visi_of_images_of_bi-geodesic_reg_inside_U}, than the proposition above, but we want to mention the proposition above in this particular form so that we can directly use it while applying localization results.
\end{remark}
The next lemma allows to construct a geodesic ray such that it is asymptotic to the image of the radial line joining $0$ and $1$ under the attached disc at a boundary point.
\begin{lemma}\label{Lem:Truncation_of_geo_and_its_asymptoticity_with_radial_line}
Let $\OM \subset\C$ be a hyperbolic domain and let $\xi\in\bdy\OM$ such that $\OM$ satisfies Condition~(M$\star$) at $\xi$. Let $(z_{n})_{n \geq 1}\subset \OM$ be a sequence converging to $\xi \in \bdy\OM$. Suppose $\varphi: \D \to \OM$ is an injective holomorphic map, $U$ be a bi-geodesic region and $(z_{n_{k}})_{k \geq 1}$ be a subsequence obtained by $(z_{n})_{n \geq 1}$ and $\xi$ as in Condition~(M$\star$). Let $U_1$ be a bi-geodesic region with $U_1\subsetneq U$ and let $\sigma:[0,\infty) \to \D$ be a radial line that joins $0$ and $1$ such that there is $S>0$ with $\wt{\sigma}:=\left.\sigma\right|_{[S, \infty)}\subset U_1$, by reparameterization we may assume $\wt{\sigma}:[0, \infty) \to \D$. Then there exist a geodesic ray $\eta:[0,\infty)\to\OM$ of $(\OM,K_\OM)$ and a constant $C>0$ such that $\eta\subset \varphi(U_1)$ and such that for all $t\geq 0$
\[
K_\OM(\varphi\circ\wt{\sigma}(t), \eta(t))\leq C.
\]
\end{lemma}
\begin{proof}
Suppose $\eta_{k}:[0,a_{k}] \to \OM$ is a sequence of geodesic of $(\OM,K_\OM)$ with $\eta(0)= \varphi(0)$ and $\eta(a_k)=z_{n_{k}}$. Let $U_0$ be bi-geodesic region with $U_0\subsetneq U_1$. Without loss of generality, we assume $\varphi(z_{n_{k}}) \in \varphi(U_0)$ for all $k \geq 1$. Fix
\[
s_{k}:= \inf \{s\in[0,a_{k}]: \eta_{k}(t) \in \varphi(U_1) \forall t \in [s, a_{k}]\}.
\]
By definition, $\left.\eta_{k}\right|_{(s_{k},a_{k}]} \subset \varphi(U_1)$ for all $k \geq 1$. We now make the following claim.
\smallskip
         
\noindent{\bf Claim~1.} $(s_{k})_{k \geq 1}$ is bounded above.
\smallskip
        
\noindent {\itshape Proof of Claim~1.} Assume to the contrary that there exists a monotonically increasing sequence $(s_{k_{l}})_{l \geq 1}$ such that $s_{k_{l}} \to \infty$ as $l \to \infty$. Note that $\eta_{k_{l}}(s_{k_{l}}) \in \partial\varphi(U_1) \cap\OM$ for all $l \geq 1$ and $\eta_{k_{l}}(a_{k_{l}}) \in \varphi(U_0)$ for all $l \geq 1$. Since $\left.\eta_{k_{l}}\right|_{[s_{k_{l}}, a_{k_{l}}]}$ is geodesic of $(\OM, K_\OM)$ and, by Lemma \ref{Lem:Visi_of_images_of_bi-geodesic_reg_inside_U}, $\varphi(U_0)$ and $\OM \setminus\varphi(U_1)$ have visibility property with respect to $(\OM, K_\OM)$, so there exists a compact set $K\subset \OM$ such that $\left.\eta_{k_{l}}\right|_{[s_{k_{l}}, a_{k_{l}}]} \cap K \neq \emptyset$. Let $o_{k_{l}} \in \left.\eta_{k_{l}}\right|_{[s_{k_{l}}, a_{k_{l}}]} \cap K$. Note that $\eta_{k_{l}}(a^{\prime}_{k_{l}}) = o_{k_{l}}$ where $a^{\prime}_{k_{l}} \in [s_{k_{l}}, a_{k_{l}}]$. Since $K_\OM(\eta_{k_{l}}(0), \eta_{k_{l}}( a_{k_{l}}^{\prime})) \leq \diam(K\cup\{\eta_{k_{l}}\})$ for all $l \geq 1$, $\limsup_{l \to \infty}K_\OM(\eta_{k_{l}}(0), \eta_{k_{l}}( a_{k_{l}}^{\prime})) < \infty$. Note that $K_\OM(\eta_{k_{l}}(0), \eta_{k_{l}}( a_{k_{l}}^{\prime})) = a^{\prime}_{k_{l}}$ for all $l \geq 1$, so $\limsup_{l \to \infty}K_\OM((\eta_{k_{l}}(0), \eta_{k_{l}}(a_{k_{l}}^{\prime})) = \infty$. This gives a contradiction. 
\hfill $\blacktriangleleft$
\smallskip

By Claim~1, there exists $T\in (0, \infty)$ with $T > s_{k}$ for all $k \in \N$ such that $\left.\eta_{k}\right|_{[T, a_{k}]} \subset \varphi(U_1)$ for all $k \in \N$. Let $\wh{\eta}_{k}:=\left.\eta_{k}\right|_{[T, a_{k}]}$. Since $T < \infty$ and $(\OM,K_\OM)$ is complete, $(\wh{\eta}_{k})_{k\geq 1}$ converges locally uniformly on compact sets to a geodesic ray $\eta:[0,\infty)\to\OM$ of $(\OM, K_\OM)$. Note that $\eta \subset \varphi(U_1)$.
\smallskip
    
By Lemma \ref{Lem:Vis_of_image_of_bigeodesic_region_and_image_of_unit_disc}, $\varphi(U_1)$ and $\OM\setminus \varphi(\D)$ have visibility property with respect to $(\OM,K_\OM)$. Then, by Lemma \ref{Lem:Quasi-geo_wrt_Small_Domain}, there exists $\kappa_{0}>0$ such that for all $k \geq 1$, $\wh{\eta}_{k}$ is $(1, \kappa_{0})$-quasi-geodesic of $(\varphi(\D), K_{\varphi(\D)})$ and $\eta$ is a $(1,\kappa_{0})$-quasi-geodesic of $(\varphi(\D), K_{\varphi(\D)})$. Hence, $\varphi^{-1}\circ\wh{\eta}_{k}$ is $(1,\kappa_{0})$-almost-geodesic of $(\D, K_{\D})$ and $\varphi^{-1}\circ\eta$ is a $(1,\kappa_{0})$-almost-geodesic of $(\D, K_{\D})$. Moreover, $(\varphi^{-1}\circ\wh{\eta}_{k})_{k \geq 1}$ converges uniformly on compact subset to a $(1,\kappa_{0})$-almost-geodesic $\wh{\sigma}:[0,\infty)\to\D$ of $(\D, K_{\D})$. Note that $(\varphi^{-1}(z_{n_{k}}))_{k \geq 1}$ converges to $1$. It is well known that $(\D, K_{\D})$ is Gromov hyperbolic. We refer the reader to \cite[Chapter~III.H]{Bridson_Haefliger} for further details on Gromov hyperbolic spaces. Thus, there exists a constant $C^\prime>0$ such that for all $t \geq 0$
\[
K_{\D}(\wt{\sigma}(t), \wh{\sigma}(t))\leq C^{\prime}.
\]
Note that $(\varphi^{-1}\circ\wh{\eta}_{k})_{k \geq 1}$ converges to $\varphi^{-1}\circ\eta$ uniformly on compact subset of $[0,\infty)$. Therefore, by uniqueness of the limit, $\wh{\sigma}(t) = \varphi^{-1}\circ\eta(t)$ for all $t \geq 0$. Therefore, for all $t \geq 0$
\[
K_\OM(\varphi\circ\wt{\sigma}(t),\eta(t)) \leq K_{\varphi(\D)}(\varphi\circ\wt\sigma(t), \eta(t)) = K_{\D}(\wt\sigma(t), \varphi^{-1}\circ\eta(t))\leq C^{\prime}.
\]
This completes the proof.
\end{proof}

\begin{lemma}\label{Lem:Geo_of_Horo_can_be_truncated_and_it_is_asyp_wrt_limit_of_geo_seq}
Let $\OM\subset \C$ be a hyperbolic domain and let $R>0$. Let $\gamma:[0,\infty) \to \OM$ be a geodesic ray of $(\OM,K_\OM)$ such that $\clos{H(\gamma,R)}\cap\bdy\OM\neq\emptyset$. Let $x\in \clos{H(\gamma,R)}\cap\bdy\OM\neq\emptyset$ and $(z_n)_{n\geq 1}\in H(\gamma,R)$ be a sequence such that $z_n\to x$ as $n\to\infty$. Suppose that $\OM$ satisfies Condition~(M$\star$) at $x$. Let $\varphi:\D \to\OM$ be an injective holomorphic map, $U$ be a bi-geodesic region and $(z_{n_k})_{k\geq 1}$ be a subsequence obtained by $x$ and $(z_n)_{n\geq 1}$ as in Condition~(M$\star$). Let $U_1\subsetneq U$ be a bi-geodesic region. Then
\begin{itemize}
\item[(i)] there exists $T>0$ such that $\wh\gamma:=\left.\gamma\right|_{[T,\infty)}\subset \varphi(U_1)$, and
\item[(ii)] there exists a constant $C_1>0$ such that for all $t\in [0,\infty)$.
\[
K_\OM(\varphi\circ\wt\sigma(t),\wh\gamma(t))\leq C_1.
\]
\end{itemize}
where up to reparameterization we assume $\wh{\gamma}:[0,\infty) \to \OM$. In particular, there exists a constant $C_2>0$ such that
\[
K_\OM(\eta(t),\wh\gamma(t))\leq C_2.
\]
where $\eta:[0,\infty)\to\OM$ is a geodesic ray of $(\OM,K_\OM)$ as in Lemma~\ref{Lem:Truncation_of_geo_and_its_asymptoticity_with_radial_line}.
\end{lemma}
\begin{proof}
Let $\sigma:[0,\infty) \to \D$ be a radial line joining $0$ and $1$. Let $S>0$ such that $\wt{\sigma}:=\left.\sigma\right|_{[S, \infty)}\subset U_1$. Up to reparameterization we may assume $\wt{\sigma}:[0, \infty) \to \D$. By Lemma~\ref{Lem:Truncation_of_geo_and_its_asymptoticity_with_radial_line}, there exists a constant $C^{\prime}>0$ such that for all $t\geq 0$
\begin{equation}\label{Ineq:Asym_of_two_rays-02}
K_\OM(\varphi\circ\wt{\sigma}(t), \eta(t))\leq C^\prime.
\end{equation}

\noindent{\bf Claim.} For each $m\in \N$ and bi-geodesic region $U^m\subsetneq U_1$, there exists $t_m\to\infty$ such that $\gamma(t_m)\in\varphi(U^m)$.
\smallskip

\noindent{\itshape Proof of Claim.} Assume that the claim is not true. Then there exists a subsequence $(t_{m_n})_{n\geq 1}\subset[0,\infty)$ such that $\gamma(t_{m_n})\not\in\varphi(U^{m_n})$. Therefore, for each $n\in\N$, $\gamma(t_{m_n})\in\OM\setminus\varphi(U^{m_n})$. For each $n\in\N$, consider a bi-geodesic region $U^{m_n}_0\subsetneq U^{m_n}$. Then, by Lemma~\ref{Lem:Visi_of_images_of_bi-geodesic_reg_inside_U}, for each $n\in \N$, $\varphi(U^{m_n}_0)$ and $\OM\setminus\varphi(U^{m_n})$ have visibility property with respect to $(\OM,K_\OM)$. Let $A:=\cup_{n\geq 1}\{\gamma(t_{m_n})\}$ and $B:=\cup_{n\geq 1}\{z_{m_n}\}$. Therefore, $A$ and $B$ have visibility property with respect to $(\OM,K_\OM)$. Note that $\gamma$ need not be subset of $A$. By following exact the same proof as the proof of Lemma~\ref{Lem:Tail_of_Seq_cannot_lie_in_geodesic_horosphere}, we get that there exists $N\in\N$ such that $z_{m_n}\not\in H_\OM(\gamma,R)$ for all $n\geq N$. This contradicts to the fact that $z_n\in H_\OM(\gamma,R)$ for all $n\geq 1$. This completes the proof of the claim.
\smallskip
\hfill $\blacktriangleleft$

Since each $U^m$ contains $1$ as an interior point in its boundary, it follows that $\varphi^{-1}(\gamma(t_{m}))\to1$ as $m \to \infty$. Let $U_0$ be a bi-geodesic region such that $U_0\subsetneq U_1$. Then, without loss of generality, we may assume that $\varphi^{-1}(\gamma(t_{m}))\in U_0$ for all $m\geq 1$. Therefore, by Lemma~\ref{Lem:Visi_of_images_of_bi-geodesic_reg_inside_U}, $\varphi(U_0)$ and $\OM \setminus \varphi(U_{1})$ have visibility property with respect to $(\OM, K_\OM)$. Then, by Truncation Lemma~\ref{Lem:Part_of_Geo_lies_in}, there exists $T>0$ such that $\left.\gamma\right|_{[T^,\infty)}\subset\OM\setminus(\OM\setminus\varphi(U_1)) = \varphi(U_1)$. By Lemma~\ref{Lem:Vis_of_image_of_bigeodesic_region_and_image_of_unit_disc} and Lemma~\ref{Lem:Quasi-geo_wrt_Small_Domain}, there exists $\k>0$ such that $\wh\gamma$ is a $(1,\k)$-quasi-geodesic of $(\varphi(\D),K_{\varphi(\D)})$. Thus, $\varphi^{-1}\circ\wh\gamma$ is a $(1,\kappa)$-quasi-geodesic ray of $(\D,K_\D)$. Since $(\D, K_\D)$ is Gromov hyperbolic space, so it follows that there exists a constant $C_1>0$ such that 
\begin{equation}\label{Ineq:Ineq:Asym_of_two_rays-03}
K_\OM(\varphi\circ\wt{\sigma}(t), \wh{\gamma}(t))\leq C_1. 
\end{equation}
Therefore, by triangle inequality and inequalities~\eqref{Ineq:Ineq:Asym_of_two_rays-03} and~\eqref{Ineq:Asym_of_two_rays-02}, we get for all $t\geq 0$
\[
K_\OM(\eta(t), \wh{\gamma}(t))\leq C^\prime +C_1.
\]
This completes the proof if we take $C_2\defeq C^\prime + C_2$.
\end{proof}

We now record the following remark; details can be found in \cite{CGMS2024}.
\begin{remark}\label{Rem:Conv_res_of_geo}
Let $\OM\subset\C$ be a hyperbolic domain and let $S\subset\bdy\OM$ be a totally disconnected set. Suppose that $\gamma_n : [a_n,b_n] \to \OM$ is a sequence of geodesic of $(\OM, K_\OM)$ such that $(\gamma_n(a_n))_{n\geq 1}$ and $(\gamma_n(b_n))_{n\geq 1}$ converge to distinct points $\xi_1\in\bdy\OM$ and $\xi_2\in\bdy_{\infty}\OM$ respectively, and such that $(\gamma_{n})_{n\geq1}$ eventually avoids every compact subset of $\OM$. Then there exist $x\in\bdy\OM\setminus S$ with $x \neq \xi_1\neq\xi_2$, $R>0$ and a sequence $(t_n)_{n\geq 1}\subset[a_n,b_n]$ such that $(\gamma_n(t_n))_{n\geq1}$ converges to $x$ and such that each $\gamma_n$ intersects $\OM\setminus B(x,R)$. 
\end{remark}

\begin{lemma}\label{Lem:Seq_of_geo_lies_in_image_of_bi-geo_region}
Let $\OM\subset \C$ be a hyperbolic domain and let $S\subset\bdy\OM$ be a totally disconnected set such that $\OM$ satisfies Condition~(M$\star$) outside the set $S$. Suppose that $\gamma_n : [a_n,b_n] \to \OM$ is a sequence of geodesic of $(\OM, K_\OM)$ such that $(\gamma_n(a_n))_{n\geq 1}$ and $(\gamma_n(b_n))_{n\geq 1}$ converge to distinct points $\xi_1\in\bdy\OM$ and $\xi_2\in\bdy\OM$ respectively, and such that $(\gamma_{n})_{n\geq1}$ eventually avoids every compact subset of $\OM$. Let $x\in\bdy\OM\setminus S$ as in Remark~\ref{Rem:Conv_res_of_geo}. Let $\varphi:\D \to \OM$ be an injective holomorphic map and $U$ be a bi-geodesic region obtained by the boundary point $x$ and $(\gamma_n(t_{n}))_{n\geq 1}$ as in Condition~(M$\star$). If $U_{1}$ is a bi-geodesic region such that $U_1\subsetneq U$, then there exists $N\in \N$ such that $\gamma_{n}\subset \varphi(U_1)$ for all $n\geq N$. Furthermore, both sequences $(\varphi^{-1}(\gamma_{n}(a_{n})))_{n \geq N}$ and $(\varphi^{-1}(\gamma_{n}(b_{n})))_{n \geq N}$ converge to $1$.
\end{lemma}

\begin{proof}
Assume to the contrary. Then there exist $s_i\in [a_i,b_i]$ and $m_i\in\N$ such that $m_i\to\infty$ as $i\to\infty$ and $\gamma_{m_{i}}(s_{i}) \not\in \varphi(U_{1})$ for all $i \geq 1$. Hence $\gamma_{m_{i}}(s_{i}) \in \OM\setminus\varphi(U_{1})$. Let $(\gamma_{n_{k}}(t_{n_{k}}))_{k \geq 1} \subset U$ be a subsequence such that $\varphi^{-1}(\gamma_{n_{k}}(t_{n_{k}})) \to 1$ as $k \to \infty$ as in Condition~(M$\star$). Pick a sequence $i_k \rightarrow \infty$ with $m_{i_k} > n_k$. Without loss of generality, we may assume that $t_{n_{k}}\leq s_{i_{k}}$ for all $k \geq 1$. Let $U_0$ be a bi-geodesic region such that $U_0 \subsetneq U_1$. We may assume that, if necessary passing to a subsequence, $\gamma_{m_{i_k}}(t_{n_{k}}) \in \varphi(U_0)$ for all $k \geq 1$. Therefore, $(\gamma_{m_{i_k}}|_{[t_{n_{k}}, s_{i_{k}}]})_{k \geq 1}$ eventually  avoids all compact subsets of $\OM$. This implies $\varphi(U_0)$ and $\OM \setminus\varphi(U_1)$ does not have visibility property with respect to $(\OM,K_\OM)$ which contradicts Lemma~\ref{Lem:Visi_of_images_of_bi-geodesic_reg_inside_U}. Therefore, there exists $N\in \N$ such that $\gamma_{n}\subset \varphi(U_1)$ for all $n\geq N$.
\smallskip

Assume that $(\varphi^{-1}(\gamma_{n}(a_{n})))_{n \geq N}$ does not converge to $1$. Then there exist a sequence $(n_l)_{l\geq 1}$ with $n_l\to\infty$ as $l\to\infty$ and $\mu\in\bdy\D$ such that $(\varphi^{-1}(\gamma_{n_l}(a_{n_l})))_{l \geq 1}$ converges to $\mu$. Since $\gamma_n\subset\varphi(U_1)$ for all $n\geq N$, so it follows that $\mu\in\bdy\D\cap\bdy U_1$. Without loss of generality, we assume that $\mu \in(\bdy U_1\cap\bdy\D)\setminus \clos{\bdy U_1\cap \D}$. Let $U_0^\prime\subsetneq U_1$ be a bi-geodesic region such that $\mu \not\in\clos{U_0^\prime}\cap\bdy\D$. Let $\wt U_0$ be a bi-geodesic region with $\wt U_0\subsetneq U_0^\prime$. Without loss of generality, we may assume $\gamma_{n_{k}}(t_{n_{k}}) \in \varphi(\wt U_0)$ for all $k \geq 1$. We may choose a sequence $l_k \to \infty$ such that $n_k\leq n_{l_k}$ for all $k\in \N$. Thus, $(\gamma_{n_{k}}|_{[t_{n_{k}}, n_{l_{k}}]})_{k \geq 1}$ eventually avoids all compact subset of $\OM$. Therefore, $\varphi(\wt U_0)$ and $\OM\setminus\varphi(U_0)$ does not have visibility property with respect to $(\OM,K_\OM)$ which contradicts Lemma~\ref{Lem:Visi_of_images_of_bi-geodesic_reg_inside_U}. Hence $(\varphi^{-1}(\gamma_{n}(a_{n})))_{n \geq N}$ converges to $1$.
\smallskip

Similarly, we may prove that $(\varphi^{-1}(\gamma_{n}(b_{n})))_{n \geq N}$ converges to $1$. This completes the proof of lemma.
\end{proof}

\section{Denjoy-Wolff theorem for planar domains}\label{Sec:Proofs_of_main_theorems}
In this section, we present the proof of Theorem~\ref{Thm:Characterization_of_Condition(*)}, Corollary~\ref{Cor:Finitely_conn_dom_charac_for_DWP}, Corollary~\ref{Cor:Well_seperated_domains_charac} and Theorem~\ref{Thm:Main_thm_outside_totally_disc_set}. We shall first collect some results which will be used immediately in the proof of Theorem~\ref{Thm:Characterization_of_Condition(*)} and in the upcoming proofs. We recall the following definition in \cite{Benini-Bracci}.
\begin{definition}
Let $f:\D\to\C$ be a holomorphic function and $\xi\in\bdy\D$. Then the {\em H-cluster set} of radius $R>0$ of $f$ at $\xi$ is defined as
\[
\Gamma_H(f,\xi,R):\{p\in\C_\infty: \exists (\zeta_n)_{n\geq 1}\subset H(\xi,R), \lim_{n\to\infty}\zeta_n = \xi, \lim_{n\to\infty}f(\zeta_n) = p\}.
\]
\end{definition}
Benini and Bracci proved that H-cluster set of a holomorphic function at any boundary point is independent of radius. We also recall it.
\begin{result}[{\cite[Lemma~2.5]{Benini-Bracci}}]\label{Res:H-cluster_set_ind_of_R}
Let $\varphi:\D\to \C$ be holomorphic. Then, for every $R, R'>0$
\[
\Gamma_H(\varphi,1, R)=\Gamma_H(\varphi,1, R').
\]
In particular,  $\varphi$ has $H$-limit at $1$ if and only if $\Gamma_H(\varphi,1, R)$ is a point for some--and hence any--$R>0$.
\end{result}
Therefore, using Result~\ref{Res:H-cluster_set_ind_of_R}, consider the H-cluster set of holomorphic map $\varphi:\D\to\OM$ at $\xi\in\bdy\D$ as
\[
\Gamma_H(\varphi, \xi)\defeq \Gamma_H(\varphi, \xi, R)
\]
where $R>0$ is any positive real number. 

The following proposition is a generalization of Proposition~2.6 in \cite{Benini-Bracci} for the domain satisfying the Condition~(M$\star$).
\begin{proposition}\label{Prop:Relation_bewtween_H-cluster_and_geo_horosphere}
Suppose $\OM \subset \C$ is a domain and $x\in\bdy\OM$. Suppose that $\OM$ satisfies Condition~(M$\star$) at $x$. Suppose $\gamma :[0, \infty) \rightarrow \OM$ is a geodesic ray of $(\OM, K_\OM)$ such that there exists a sequence $(t_{n})_{n\geq1}\subset [0,\infty)$ with $t_{n} \to \infty$ as $n \to \infty$ and $\gamma(t_{n}) \rightarrow x \in \partial{\OM}$ as $n \rightarrow \infty$. Suppose $\varphi: \D \rightarrow \OM$ is an injective holomorphic map obtained by $x$ and $(\gamma(t_{n}))_{n \geq 1}$ as in Condition~(M$\star$). Then for all $R > 0$,
\[
\clos{H_{\OM}(\gamma, R)} \cap \partial{\OM} = \Gamma_{H}(\varphi, 1).
\]
\end{proposition}

\begin{proof}
Let $U$ be a bi-geodesic region and $(\gamma(t_{n_{k}}))_{k\geq 1}$ be a subsequence as in Condition~(M$\star$). By Proposition~\ref{Prop:Vis_of_all_req_sets_which_works_everywhere}, there exist bi-geodesic regions $U_1\subsetneq U_2\subsetneq U_3 \subsetneq U$ such that $\varphi(U_1)$ and $\OM\setminus\varphi(U_2)$, $\varphi(U_2)$ and $\OM\setminus\varphi(U_3)$, and $\varphi(U_3)$ and $\OM\setminus\varphi(\D)$ have visibility property with respect to $(\OM, K_\OM)$. Let $U_0\subsetneq U_1$ be a bi-geodesic region. Without loss of generality, we may assume $(\gamma(t_{n_{k}}))_{k\geq1} \in \varphi(U_0)$ for all $k \in \N$. Then, by Lemma~\ref{Lem:Part_of_Geo_lies_in}, there exists $T>0$ such that $\wh\gamma:= \left.\gamma\right|_{[T,\infty)}\subset\varphi(U_1)$. By reparameterization, we may assume that $\wh\gamma:[0,\infty) \to \OM$. By Lemma~\ref{Lem:Quasi-geo_wrt_Small_Domain}, there exists $\kappa_{0}>0$ such that $\wh\gamma$ is a ($1,\kappa_{0}$)-quasi-geodesic of $(\varphi(\D), K_{\varphi(\D)})$. Then $\wh\sigma:= \varphi^{-1}\circ\wh\gamma$ is $(1,\kappa_{0})$-quasi-geodesic of $(\D, K_{\D})$.
\smallskip

Since $(\varphi^{-1}\circ\wh\gamma(t_{n_{k}}))_{k \geq 1}$ converges to $1$, so $(\wh\sigma(t_{n_{k}}))$ converges to $1$. Suppose $\sigma:[0,\infty) \to \D$ is a radial line that joins $0$ and $1$. Let $S>0$ such that $\wt{\sigma}:=\left.\sigma\right|_{[S, \infty)} \subset U_1$. Up to a reparameterization, we may assume $\wt{\sigma}:[0,\infty) \to \D$.  Since $(\D, K_{\D})$ is Gromov hyperbolic space, there exists a constant $C_{0}>0$ such that for all $t \geq 0$
\begin{equation}\label{Ineq:Asymp_Condition_for_quasi_geo_3}
K_{\D}(\wh\sigma(t),\wt{\sigma}(t)) \leq C_{0}.
\end{equation}
Note that $\wt{\gamma} := \varphi\circ\wt{\sigma}$ is a geodesic ray of $(\varphi(\D), K_{\varphi(\D)})$. Using~\eqref{Ineq:Asymp_Condition_for_quasi_geo_3}, we get for all $t \in [0,\infty)$
\begin{equation}
K_\OM(\wh\gamma(t), \wt{\gamma}(t)) \leq K_{\varphi(\D)}(\wh\gamma(t), \wt{\gamma}(t))= K_\D(\wh\sigma(t),\wt{\sigma}(t)) \leq C_{0}.
\end{equation}
Then, by Theorem~\ref{Thm:Geo_horosphere_of_big_domain_is_subset_of_geod_of_small_domain_upto_some_compact_set_and-constant}, there exist a constant $C_{1}$ and compact set $K\subset \OM$ such that
\begin{equation*}
H_{\OM}(\wh\gamma, e^{2T}R) \setminus K \subseteq H_{\varphi(\D)}(\wt{\gamma}, C_{1}e^{2T}R).
\end{equation*}
Therefore, by Lemma~\ref{Lem:Relation_between_origin_geo_horo_and_truncated_geo_horo}, we get
\begin{equation}\label{ineq:Geo_horosphere_of_big_domain_is_subset_of_geod_of_small_domain_upto_some_compact_set_and_constant_1}
H_{\OM}(\gamma, R) \setminus K \subseteq H_{\varphi(\D)}(\wt{\gamma}, C_{1}e^{2T}R).
\end{equation}
Moreover, by Theorem~\ref{Thm:Geo_horosphere_of_small_domain_is_subset_of_geo_horosphere_of_big_domain}, there exists a constant $C_{2}>0$ such that for all $R>0$
\begin{equation*}
H_{\varphi(\D)}(\wt{\gamma}, e^{2T}R) \subseteq H_{\OM}(\wh\gamma, e^{2T}C_{2}R).
\end{equation*}
Therefore, by Lemma~\ref{Lem:Relation_between_origin_geo_horo_and_truncated_geo_horo}, we get that for all $R>0$
\begin{equation}\label{ineq:Geo_horosphere_of_small_domain_is_subset_of_geo_horosphere_of_big_domain_1}
H_{\varphi(\D)}(\wt{\gamma}, e^{2T}R) \subseteq H_{\OM}(\gamma, C_{2}R).
\end{equation}
\noindent{\bf Claim~1.} $\clos{H_{\OM}(\gamma, R)} \cap \partial\OM \subseteq \Gamma_{H}(\varphi,1)$.
\smallskip
        
\noindent {\itshape Proof of Claim~1.} Let $y \in \clos{H_{\OM}(\gamma, R)} \cap \partial\OM$. There exists a sequence $(z_{n})_{n \geq 1} \subset H_{\OM}(\gamma, R)$ such that $z_{n} \to y $ as $n \to \infty$. Without loss of generality, we may assume that $(z_{n})_{n \geq 1} \in H_{\OM}(\gamma, R) \setminus K$. Therefore, by~\eqref{ineq:Geo_horosphere_of_big_domain_is_subset_of_geod_of_small_domain_upto_some_compact_set_and_constant_1}, we have $z_{n} \in H_{\varphi(\D)}(\wt{\gamma}, C_{1}e^{2T}R)$ for all $n \geq 1$. There exists $\zeta_{n} \in H_{\D}(\wt{\sigma}, C_{1}R)$ with $\zeta_{n} \to 1$ as $n \to \infty$ such that $\varphi(\zeta_{n}) = z_{n}$. By Lemma~\ref{Lem:Relation_between_origin_geo_horo_and_truncated_geo_horo}, $H_{\D}(\sigma, e^{2(T-S)}C_1R) = H_{\D}(\wt{\sigma}, e^{2T}C_{1}R)$. Therefore, $y \in \Gamma_{H}(\varphi,1, e^{2(T-S)}C_{1}R)$. This completes the proof of the claim by Result \ref{Res:H-cluster_set_ind_of_R}.
\hfill $\blacktriangleleft$
\smallskip

\noindent{\bf Claim~2.} $\Gamma_{H}(\varphi,1) \subseteq \clos{H_{\OM}(\gamma, R)} \cap \partial\OM$.
\smallskip
        
\noindent {\itshape Proof of Claim~2.} Let $y \in \Gamma_{H}(\varphi,1)$. There exists a sequence $(\zeta_{n})_{n\geq1}\subset H_{\D}(\sigma, e^{2(T-S)}R)$ with $\zeta_{n} \to 1$ as $n \to \infty$ such that $\varphi(\zeta_{n}) \to y$ as $n \to \infty$. Therefore, by Lemma \ref{Lem:Relation_between_origin_geo_horo_and_truncated_geo_horo}, $\varphi(\zeta_{n}) \in H_{\varphi(\D)}(\wt\gamma, e^{2T}R)$ for all $n \in \N$. By~\eqref{ineq:Geo_horosphere_of_small_domain_is_subset_of_geo_horosphere_of_big_domain_1}, $\varphi(\zeta_{n}) \in H_{\OM}(\gamma, C_{2}R)$ for all $n \in \N$. Note that $(\varphi(\zeta_n))_{n\geq 1}\subset H(\gamma,C_2R)$ for all $R>0$.  This completes the proof of the claim.
\hfill $\blacktriangleleft$
\smallskip

This completes the proof of proposition.
\end{proof}

The next lemma shows that if closed Busemann horosphere intersects the boundary of domain at single point for some radius $R_0>0$, then it intersects the boundary at the same single point for every radius $R>0$.
\begin{lemma}\label{Lem:Horo_single_for_some_radius_implies_single_for_all_radius}
Suppose $\OM\subset\C$ satisfies Condition~(M$\star$) and $\gamma:[0,\infty)\to\OM$ is a geodesic ray of $(\OM,K_\OM)$. If there exist $x\in\bdy\OM$ and $R_0>0$ such that $\clos{H(\gamma,R_0)}\cap\bdy\OM=\{x\}$, then $\clos{H(\gamma,R)}\cap\bdy\OM=\{x\}$ for all $R>0$.
\end{lemma}
\begin{proof}
Since $\clos{H(\gamma,R_0)}\cap\bdy\OM=\{x\}$, so it follows by Lemma~\ref{Lem:Cluster_set_of_geo_sits_in_geo_horosphere} that $\mathcal{C}(\gamma)=\{x\}$. Let $(t_n)_{n\geq 1}\subset[0,\infty)$ such that $t_n\to\infty$ as $n\to\infty$ and such that $\gamma(t_n)\to x$ as $n\to\infty$. Let $\varphi:\D\to\OM$ be an injective holomorphic map obtained by $(\gamma(t_n))_{n\geq 1}$ and $x$. Therefore, by Proposition~\ref{Prop:Relation_bewtween_H-cluster_and_geo_horosphere}, we get for all $R>0$
\begin{equation}\label{Ineq:Geo_horo_and_H-cluster_set}
\clos{H_{\OM}(\gamma, R)} \cap \partial\OM = \Gamma_{H}(\varphi,1).
\end{equation}
Since $\clos{H_{\OM}(\gamma, R_0)} \cap \partial\OM = \{x\}$, so it follows from~\eqref{Ineq:Geo_horo_and_H-cluster_set} that $\Gamma_{H}(\varphi,1)=\{x\}$. Using~\eqref{Ineq:Geo_horo_and_H-cluster_set} again, we get that $\clos{H_{\OM}(\gamma, R)} \cap \partial\OM=\{x\}$ for all $R>0$. This completes the proof.
\end{proof}

\begin{remark}\label{Rem:Closed_geo_horo_empty_for_some_radius_implies_it_is_empty_for_all}
Using an argument similar to that in the proof of Lemma~\ref{Lem:Horo_single_for_some_radius_implies_single_for_all_radius}, we obtain the following: Suppose $\OM\subset\C$ satisfies Condition~(M$\star$) and $\gamma:[0,\infty)\to\OM$ is a geodesic ray of $(\OM,K_\OM)$. If there exists $R_0>0$ such that $\clos{H(\gamma,R_0)}\cap\bdy\OM=\emptyset$, then $\clos{H(\gamma,R)}\cap\bdy\OM=\emptyset$ for all $R>0$.
\end{remark}

We now record few more definitions and standard results which are required in the proof of Theorem~\ref{Thm:Characterization_of_Condition(*)}.
\begin{definition}
Let $\OM \subset \C^{d}$ be a complete hyperbolic domain and $(z_{n})_{n\geq 1} \subset \OM$. The {\em sequential horosphere of $(z_{n})_{n\geq 1}$ of radius $R>0$ having pole at $p \in \OM$} is
\[
H_{p}((z_{n})_{n\geq 1}, R):= \{z\in\OM: \limsup_{n\to \infty}[K_\OM(z, z_{n})-K_\OM(p, z_{n})] < \frac{1}{2}\log R\}
\]
\end{definition}

\begin{result}[{\cite[Lemma~3.1]{Bracci-Gaussier}}]\label{Res:Relation_between_seq_horo_with_diff_base_point} Suppose $\OM \subset \C^{d}$ is complete Kobayashi hyperbolic domain, $p, q\in \OM$ and $(z_{n})_{n \geq 1}$. Then there exist $\alpha, \beta>0$ such that for all $R>0$
\[
H_{q}((z_{n})_{n \geq 1}, \alpha R) \subseteq H_{p}((z_{n})_{n \geq 1}, R) \subseteq H_{q}((z_{n})_{n \geq 1}, \beta R).
\]
\end{result}

\begin{lemma}\label{Lem:Sequential_horo_and_geo_horo_are_comparable_in_unit_disc}
Let $(z_{n})_{n\geq 1} \subset \D$ be a sequence that converges to $1$ and $\sigma$ be a radial line that joins $0$ and $1$. Then there exists a constant $\beta>0$ such that for all $R>0$
\[
H_{p}((z_{n})_{n \geq 1}, R) \subseteq  H_{\D}(\sigma, \beta R).
\]
\end{lemma}
\begin{proof}
The proof follows from the fact that $H_{0}((z_{n})_{n \geq 1}, R) = H_{\D}(\sigma, R)$ for all $R>0$  and Result~\ref{Res:Relation_between_seq_horo_with_diff_base_point}. Therefore, we omit the proof.
\end{proof}

\begin{result}[Lehto-Virtanen, {\cite[Theorem~3.3.1]{BCM}}]\label{Res:Lehto-Virtanen_Thm} Let $f:\D\to \C$ be an injective holomorphic map and $-\infty <a<b\leq +\infty$. Let $\gamma:[a,b) \to \D$ be a continuous curve such that $\lim_{t \to b}\gamma(t) = \xi\in\bdy\D$. Suppose that $\lim_{t\to b}f(\gamma(t)) = L$ for some $L \in \C_\infty$. Then non-tangential of $f$ at $\xi$ exists and it is equal to $L$.
\end{result}

\begin{proof}[Proof of {\Cref{Thm:Characterization_of_Condition(*)}}]
$(2) \Rightarrow (3)$ This is trivial.
\smallskip

$(3)\implies (2)$. Let $\gamma:[0,\infty)\to\OM$ be a geodesic of $(\OM,K_\OM)$. If $\clos{H_{\OM}(\gamma,R)}\cap \partial\OM = \emptyset$ for some $R>0$, then, by Remark~\ref{Rem:Closed_geo_horo_empty_for_some_radius_implies_it_is_empty_for_all}, we are done. Therefore, we assume that $\clos{H_{\OM}(\gamma,R)}\cap \partial\OM\neq\emptyset$ for some $R>0$. Then we shall prove that there exists $x\in\bdy\OM$ such that $\clos{H_{\OM}(\gamma,R)}\cap \partial\OM=\{x\}$. Thus, it follows from Lemma~\ref{Lem:Horo_single_for_some_radius_implies_single_for_all_radius} that $\clos{H_{\OM}(\gamma,R)}\cap \partial\OM = \{x\}$ for all $R>0$. Since $\clos{H_{\OM}(\gamma,R)}\cap \partial\OM\neq\emptyset$, so let $x \in \clos{H_{\OM}(\gamma,R)}\cap \partial\OM$. There exists a sequence $(z_n)_{n\geq }\subset \OM$ such that it converges to $x$. Let $\varphi:\D\to\OM$ be an injective holomorphic map, $U$ be a bi-geodesic region and $(z_{n_k})_{k\geq 1}$ be a subsequence as in Condition~(M$\star$). By Proposition~\ref{Prop:Vis_of_all_req_sets_which_works_everywhere}, there exist bi-geodesic regions $U_1\subsetneq U_2\subsetneq U_3 \subsetneq U$ such that $\varphi(U_1)$ and $\OM\setminus\varphi(U_2)$, $\varphi(U_2)$ and $\OM\setminus\varphi(U_3)$, and $\varphi(U_3)$ and $\OM\setminus\varphi(\D)$ have visibility property with respect to $(\OM, K_\OM)$.

Let $z_0\in\OM$ as in the hypothesis of $(3)$. Then, by Lemma~\ref{Lem:Truncation_of_geo_and_its_asymptoticity_with_radial_line}, there exists a geodesic ray $\eta:[0,\infty)\to\OM$ of $(\OM,K_\OM)$ such that $\eta(0)=z_0$ and such that $\eta\subset \varphi(U_1)$. By Lemma~\ref{Lem:Geo_of_Horo_can_be_truncated_and_it_is_asyp_wrt_limit_of_geo_seq}, there exist $T>0$ and a constant $C^\prime>0$ such that $\wh\gamma:=\left.\gamma\right|_{[T,\infty)}\subset \varphi(U_1)$ and such that for all $t\in [0,\infty)$, $K_\OM(\eta(t),\wh\gamma(t))\leq C^\prime$.
Therefore, by Lemma~\ref{Lem:Asym_compares_two_geo_horo}, there exists a constant $C>0$ such that
\[
H_{\OM}(\wh\gamma,e^{2T}R)\subseteq H_{\OM}(\eta,e^{2T}CR).
\]
Therefore, by Lemma~\ref{Lem:Relation_between_origin_geo_horo_and_truncated_geo_horo}, we get
\begin{equation}\label{Ineq:Asym_implies_horo_are_same}
H_{\OM}(\gamma,R)\subseteq H_{\OM}(\eta,e^{2T}CR).
\end{equation}
Since $z_{n_k}\in H_{\OM}(\gamma,R)$ for all $k\geq 1$, so it follows from \eqref{Ineq:Asym_implies_horo_are_same} that $\clos{H_{\OM}(\eta,e^{2T}CR)}\cap\bdy\OM \neq \emptyset$. By hypothesis, if there exists $R>0$ such that $\clos{H_\OM(\eta,R)}\cap\bdy\OM=\emptyset$, then using Remark~\ref{Rem:Closed_geo_horo_empty_for_some_radius_implies_it_is_empty_for_all} we get that $\clos{H_\OM(\eta,R)}\cap\bdy\OM=\emptyset$ for all $R>0$. This is not possible due to the fact that $x \in \clos{H_\OM(\eta,e^{2T}CR)}\cap\bdy\OM$. Thus, by hypothesis there exist $y\in\bdy\OM$ and $R_0>$ such that
\begin{equation*}\label{ineq:Horo_int_at_single}
\clos{H_{\OM}(\eta,R_0)}\cap \partial\OM = \{y\}.
\end{equation*}
By Lemma~\ref{Lem:Horo_single_for_some_radius_implies_single_for_all_radius}, we have $\clos{H_{\OM}(\eta,R)}\cap \partial\OM = \{y\}$ for all $R>0$. Therefore, $\clos{H_\OM(\eta, e^{2T}CR)}\cap \partial\OM = \{y\}$. Hence, using~\eqref{Ineq:Asym_implies_horo_are_same}, we get
\[
\clos{H_\OM(\gamma, R)}\cap\bdy\OM\subseteq \clos{H_\OM(\eta, e^{2T}CR)}\cap\bdy\OM  = \{y\}.
\]
Since $z_n\in H_\OM(\gamma, R)$ for all $n\geq 1$, so it follows that $x=y$. Therefore, $\clos{H_\OM(\gamma, R)}\cap\bdy\OM = \{x\}$. This completes the proof of $(3) \implies (2)$.
\smallskip

$(2) \Rightarrow (4)$ Let $x\in\bdy\OM$ and $(z_{n})_{n \geq 1}$ be a sequence that converges to $x$. Let $\varphi :\D\to \OM$ is injective holomorphic map obtained by $x$ and $(z_{n})_{n \geq 1}$ as in Condition~($\star$). Let $U$ be a bi-geodesic region and $(z_{n_{k}})_{k\geq 1}$ be a subsequence with $(z_{n_{k}})_{k\geq 1} \subset U$ as in Condition~(M$\star$). By Proposition~\ref{Prop:Vis_of_all_req_sets_which_works_everywhere}, there exist bi-geodesic regions $U_1\subsetneq U_2\subsetneq U_3 \subsetneq U$ such that $\varphi(U_1)$ and $\OM\setminus\varphi(U_2)$, $\varphi(U_2)$ and $\OM\setminus\varphi(U_3)$, and $\varphi(U_3)$ and $\OM\setminus\varphi(\D)$ have visibility property with respect to $(\OM, K_\OM)$. Without loss of generality, we may assume $\varphi^{-1}(z_{n_{k}}) \in U_{1}$ for all $k \geq 1$.
\smallskip
Let $\sigma:[0,\infty) \to \D$ be a radial line joining $0$ and $1$. Let $S>0$ such that $\wt{\sigma}:=\left.\sigma\right|_{[S, \infty)}\subset U_1$. Up to reparameterization we may assume $\wt{\sigma}:[0, \infty) \to \D$. By Lemma~\ref{Lem:Truncation_of_geo_and_its_asymptoticity_with_radial_line}, there exist a geodesic ray $\eta:[0,\infty)\to\OM$ of $(\OM,K_\OM)$ and a constant $C^\prime>0$ such that $\eta\subset \varphi(U_1)$ and such that for all $t\geq 0$
\[
K_\OM(\varphi\circ\wt{\sigma}(t), \eta(t))\leq C^\prime.
\]
By Theorem~\ref{Thm:Geo_horosphere_of_small_domain_is_subset_of_geo_horosphere_of_big_domain}, there exists a constant $C>0$ such that for all $R>0$
\[
H_{\varphi(\D)}(\varphi\circ\wt{\sigma}, R)\subseteq H_{\OM}(\eta, CR)
\]
Therefore, Lemma~\ref{Lem:Relation_between_origin_geo_horo_and_truncated_geo_horo}, we get that for all $R>0$
\begin{equation}\label{ineq:Geo_horosphere_of_small_domain_is_subset_of_geo_horosphere_of_big_domain_0}
H_{\varphi(\D)}(\varphi\circ\sigma, e^{-2S}R) \subseteq H_{\OM}(\eta, CR).
\end{equation}
Hence, by \eqref{ineq:Geo_horosphere_of_small_domain_is_subset_of_geo_horosphere_of_big_domain_0}, we get that for all $R>0$
\begin{equation}\label{Eq:Comp-of_geo_horo}
\clos{H_{\varphi(\D)}(\varphi\circ\sigma, e^{-2S}R)}\cap\partial\OM \subseteq \clos{H_{\OM}(\eta, CR)}\cap\partial\OM.
\end{equation}
By the hypothesis, either $\clos{H_{\OM}(\eta, R)}\cap\partial\OM = \emptyset$ for all $R>0$ or there exists $y \in \partial \OM$ such that $\clos{H_{\OM}(\eta, R)}\cap\partial\OM = \{y\}$ for all $R>0$. If $\clos{H_{\OM}(\eta, R)}\cap\partial\OM = \emptyset$ for some $R>0$, so it follows from Remark~\ref{Rem:Closed_geo_horo_empty_for_some_radius_implies_it_is_empty_for_all} that $\clos{H_{\OM}(\eta, R)}\cap\partial\OM = \emptyset$ for all $R>0$. Therefore, $\clos{H_{\OM}(\eta, CR)}\cap\partial\OM = \emptyset$. So, by~\eqref{Eq:Comp-of_geo_horo}, we get that $\clos{H_{\varphi(\D)}(\varphi\circ\sigma, e^{-2S)}R)}\cap\partial\OM=\emptyset$. That is, H-limit of $\varphi$ at 1 is infinity. If $\clos{H_{\OM}(\eta, CR)}\cap\partial\OM=\{y\}$, then, by~\eqref{Eq:Comp-of_geo_horo} and Result~\ref{Res:H-cluster_set_ind_of_R}, the H-limit of $\varphi :\D\to \OM$ at $1$ exists and it is finite. This completes the proof of $(2) \Rightarrow (4)$.
\smallskip
        
$(4) \Rightarrow (2)$ We continue the same proof as $(3)\implies (2)$ till we get a geodesic ray $\eta:[0,\infty)\to\OM$ of $(\OM,K_\OM)$. After this we now again continue the proof as $(3)\implies (2)$, till we get that there exists a constant $C^\prime>0$ and a constant $C>0$ such that for all $t\geq0$
\begin{equation}\label{Ine:both_are_asym}
K_\OM(\varphi\circ\wt\sigma(t),\eta(t))\leq C^\prime
\end{equation}
and such that
\begin{equation}\label{Ineq:Asym_implies_horo_are_same_01}
H_{\OM}(\gamma,R)\subseteq H_{\OM}(\eta,e^{2T}CR).
\end{equation}
where $\sigma:[0,\infty) \to \D$ be a radial line joining $0$ and $1$, and $S>0$ such that $\wt{\sigma}:=\left.\sigma\right|_{[S, \infty)}\subset U_1$. Up to reparameterization we may assume $\wt{\sigma}:[0, \infty) \to \D$. By Theorem~\ref{Thm:Geo_horosphere_of_big_domain_is_subset_of_geod_of_small_domain_upto_some_compact_set_and-constant}, there exist a constant $C_1>0$ and a compact set $K\subset \OM$ such that
\begin{equation}\label{Ineq:Geo_horo_comp_upto_compact_set-01}
H_\OM(\eta, e^{2T}CR)\setminus K\subset H_{\varphi(\D)}(\varphi\circ\wt\sigma, C_1e^{2T}CR).
\end{equation}
Note that, by Lemma~\ref{Lem:Relation_between_origin_geo_horo_and_truncated_geo_horo}, we have $H_{\varphi(\D)}(\varphi\circ\wt\sigma, C_1e^{2T}CR)=H_{\varphi(\D)}(\varphi\circ\sigma, e^{2(T-S)}CC_1R)$. Since $K$ is a compact set, $(z_{n_k})_{k\geq1}$ converges to the boundary point $x$ and $z_{n_k}\in H_\OM(\eta,e^{2T}CR)$ for all $k\geq1$, so without loss of generality we may assume that $z_{n_k}\in H_\OM(\eta,e^{2T}CR)\setminus K$ for all $k\geq1$. It follows from~\eqref{Ineq:Geo_horo_comp_upto_compact_set-01} that $z_{n_k}\in H_{\varphi(\D)}(\varphi\circ\sigma, e^{2(T-S)}CC_1R)$ for all $k\geq 1$. Note that $\varphi^{-1}(z_{n_k})_{k\geq 1}$ converges to $1$. This shows that H-limit of $\varphi$ at $1$ cannot be infinity. Therefore, by hypothesis, H-limit of $\varphi$ at $1$ exists and it is finite. In fact, H-limit of $\varphi$ at $1$ is equal to $x$. Thus, $\clos{H_{\varphi(\D)}(\varphi\circ\sigma, e^{2(T-S)}CC_1R)}\cap\bdy\OM=\{x\}$. By~\eqref{Ineq:Asym_implies_horo_are_same_01} and~\eqref{Ineq:Geo_horo_comp_upto_compact_set-01}, we get $\clos{H_{\OM}(\gamma,R)}\cap\bdy\OM=\{x\}$. This completes the proof of $(4) \implies (2)$.
\smallskip

$(2) \Rightarrow (1)$ Let $f \in \hbox{Hol}(\Omega,\Omega)$. If for all $z\in\OM$ the $(f^{n}(z))_{n\geq 1}$ is relatively compact in $\OM$, then we are done. Assume that there exists $z_{0} \in \OM$ such that the orbit $(f^n(z_0))$ is not relatively compact in $\OM$. Then there exists a subsequence $(f^{n_l}(z_0))_{l\geq1}$ eventually avoids every compact subset of $\OM$. Then, up to a subsequence, either $(f^{n_l}(z_0))_{l\geq1}$ converges to some boundary point $x \in \partial\OM$ or $|f^{n_l}(z_0)|\to\infty$ as $l\to\infty$. In the first case, we shall prove that $(f^{n})_{n\geq 1}$ converges uniformly on compact subset of $\OM$ to $x$. Therefore, it will force that if $|f^{n_l}(z_0)|\to\infty$ as $l\to\infty$, then $(f^n)_{n\geq1}$ converges uniformly on compact set to $\infty$. 

Let $\varphi :\D\to \OM$ be an injective holomorphic map, $U$ be a bi-geodesic region and $(f^{n_{l_k}}(z_0))_{k\geq 1}$ be a subsequence obtain by $x$ and $(f^{n_l}(z_{0}))_{l \geq 1}$ as in Condition~(M$\star$). By Proposition~\ref{Prop:Vis_of_all_req_sets_which_works_everywhere}, there exist bi-geodesic regions $U_1\subsetneq U_2\subsetneq U_3 \subsetneq U$ such that $\varphi(U_1)$ and $\OM\setminus\varphi(U_2)$, $\varphi(U_2)$ and $\OM\setminus\varphi(U_3)$, and $\varphi(U_3)$ and $\OM\setminus\varphi(\D)$ have visibility property with respect to $(\OM, K_\OM)$. Let $U_0\subsetneq U_1$ be a bi-geodesic region. By Lemma~\ref{Lem:Visi_of_images_of_bi-geodesic_reg_inside_U}, $\varphi(U_0)$ and $\OM\setminus\varphi(U_1)$ have visibility property with respect to $(\OM,K_\OM)$. Since $(f^{n_{l_k}}(z_0))_{k\geq1}$ converges to $1$, so without loss of generality, we may assume $\varphi^{-1}(f^{n_{l_k}}(z_0)) \in U_0$ for all $k \geq 1$. Then, by Lemma \ref{Lem:Whole_Orbit_lies_if_one_of_its_subseq_lies}, there exists $N\in \N$ such that $f^{n}(z_{0}) \in \OM\setminus(\OM\setminus\varphi(U_1)) = \varphi(U_1)$ for all $n \geq N$. For simplicity, we assume $f^{n}(z_{0}) \in \varphi(U_1)$ for all $n\in \N$. Now, by Theorem~\ref{Thm:Global_Loc_Res}, there exists a constant $C_0>0$ such that for all $k\in\N$ and $n\in\N$
\begin{equation}\label{ineq:Localization_Ineq_in_Main_Result-1}
K_{\varphi(\D)}(f^n(z_{0}), f^{n_{l_k}}(z_{0})) \leq K_\OM(f^n(z_{0}), f^{n_{l_k}}(z_{0})) + C_{0}.
\end{equation}
Let $f(z_{0}) = \varphi(\zeta_{0})$ for some $\zeta_{0} \in \D$. We now make the following claim.
\smallskip
        
\noindent{\bf Claim.} Let $\sigma$ be a radial line joining $0$ and $1$. Then, there exists a constant $C_{1} > 0$ such that for all $n\in\N$,
\[
\varphi^{-1}\circ f^{n}(z_{0}) \in H_{\D}(\sigma, C_{1}).
\]
\noindent {\itshape Proof of Claim.} First, we shall show that there exists $C_{2} >0$ such that for all $n \geq 1$ 
\[
\varphi^{-1}\circ f^{n}(z_{0}) \in H_{f(z_{0})}((\varphi^{-1}\circ f^{n_{l_k}}(z_{0}))_{k \geq 1}, C_{2}).
\]
For all $n\geq 1,$ if $n_{l_k} > n$ (that is, for $k$ sufficiently large), by \eqref{ineq:Localization_Ineq_in_Main_Result-1} and triangle inequality,  
\begin{align*}
K_{\D}(\varphi^{-1}(f^n(z_{0})), \varphi^{-1}(f^{n_{l_k}}(z_{0}))) &- K_{\D}(\zeta_{0},\varphi^{-1}(f^{n_{l_k}}(z_{0})))\\ 
&= K_{\varphi(\D)}(f^n(z_{0}),f^{n_{l_k}}(z_{0}))) - K_{\varphi(\D)}(f(z_0), f^{n_{l_k}}(z_{0}))\\
&\leq K_\OM(f^n(z_{0}), f^{n_{l_k}}(z_{0})) + C_{0} - K_\OM(f(z_0), f^{n_{l_k}}(z_{0}))\\
&\leq K_\OM(f^{n_{l_k}-n}(z_{0}),z_{0}) - K_\OM(f(z_0), f^{n_{l_k}}(z_{0})) + C_{0}\\
&\leq K_\OM(f^{n_{l_k}}(z_{0}),z_{0}) - K_\OM(f(z_0), f^{n_{l_k}}(z_{0})) + C_{0}\\
& \leq K_\OM(f(z_{0}),z_{0}) + C_{0}.
\end{align*}
Therefore, by taking $C_2\defeq e^{2K_\OM(f(z_0),z_{0}) + 2C_{0}}$, we obtain for all $n \geq \N$
\[
\varphi^{-1}\circ f^{n}(z_{0}) \in H_{f(z_{0})}((\varphi^{-1}\circ f^{n_{l_k}}(z_{0}))_{k \geq 1}, C_{2}).
\]
Now using the fact that the sequential horosphere of sequence converging to $1$ at pole $0$ and Busemann horosphere of radial line $\sigma$ are same for the unit disc and by Lemma~\ref{Lem:Sequential_horo_and_geo_horo_are_comparable_in_unit_disc}, there exists a constant $\beta >0$ such that for all $m\geq  1$
\[
H_{f(z_{0})}((\varphi^{-1}\circ f^{n_{l_k}}(z_{0}))_{k \geq 1}, C_{2}) \subseteq H_{\D}(\sigma, \beta C_{2}).
\]
Therefore, by taking $C_{1}:= \beta C_{2}$, the proof of the claim is complete.
\hfill $\blacktriangleleft$
\smallskip
        
Let $S>0$ such that $\wt{\sigma}\defeq\left.\sigma\right|_{[S, \infty)}\subset U_1$. By the reparameterization, we may assume that $\wt{\sigma}:[0, \infty) \to \D$. By Lemma~\ref{Lem:Truncation_of_geo_and_its_asymptoticity_with_radial_line}, there exists a geodesic ray $\eta:[0,\infty)\to\OM$ of $(\OM,K_\OM)$ and a constant $C^\prime>0$ such that $\eta\subset \varphi(U_1)$ and such that for all $t\geq0$
\[
K_\OM(\varphi\circ\wt{\sigma}(t), \eta(t))\leq C^\prime.
\]
By Theorem~\ref{Thm:Geo_horosphere_of_small_domain_is_subset_of_geo_horosphere_of_big_domain}, there exists a constant $C>0$ such that
\begin{equation}\label{Eq:comp_of_two_geo_horo}
H_{\varphi(\D)}(\varphi\circ\wt{\sigma}, e^{2S}C_1) \subseteq H_{\OM}(\eta, Ce^{2S}C_1).
\end{equation}
By Lemma~\ref{Lem:Relation_between_origin_geo_horo_and_truncated_geo_horo}, we get that
\begin{equation}\label{ineq:Geo_horosphere_of_small_domain_is_subset_of_geo_horosphere_of_big_domain_2}
\clos{H_{\varphi(\D)}(\varphi\circ\sigma,C_1)}\cap\bdy\OM \subseteq \clos{H_{\OM}(\eta, Ce^{2S}C_1)}\cap\bdy\OM.
\end{equation}
By the Claim, we have that for all $n \geq 1$
\begin{equation}\label{Orbit_lies_in_geo_horo_wrt_unit_disc}
f^n(z_{0}) \in H_{\varphi(\D)}(\varphi\circ\sigma, C_{1}).
\end{equation}
Using~\eqref{Eq:comp_of_two_geo_horo} and \eqref{Orbit_lies_in_geo_horo_wrt_unit_disc}, we get that $f^n(z_{0}) \in H_{\OM}(\eta,e^{2S}CC_{1})$ for all $n\geq1$. Thus, by~\eqref{ineq:Geo_horosphere_of_small_domain_is_subset_of_geo_horosphere_of_big_domain_2}, we get $x\in\clos{H_{\OM}(\eta, Ce^{2S}C_1)}\cap\bdy\OM$. Therefore, by the hypothesis, there exists $y \in \partial\OM$ such that
\begin{equation}\label{ineq:Closed_geo_horo_intersect_bdry_single_point_1}
\clos{H_{\OM}(\eta, Ce^{2S}C_1)}\cap\bdy\OM = \{y\}.
\end{equation}
Therefore using~\eqref{ineq:Closed_geo_horo_intersect_bdry_single_point_1}, the limit set of $\{f^n(z_{0}):n\in\N\}$ is $\{y\}$. Since $(f^{n_k}(z_{0}))_{k \geq 1}$ converges to $x$, so $x=y$. Therefore, $(f^n(z_{0}))_{n\geq 1}$ converges to $x$. By using the fact that the Kobayashi distance is non-expanding under holomorphic map, we get that for every $w_0$ the orbit $(f^{n}(w_{0}))_{n \geq 1}$ also converges to $x$. This completes the proof of $(2) \implies (1)$.
\smallskip

$(1) \Rightarrow (2)$ We continue the same proof as $(4) \Rightarrow (2)$ till we get that there exists a geodesic ray $\eta:[0,\infty)\to\OM$ of $(\OM,K_\OM)$ such that for all $t\geq 0$
\begin{equation}\label{Ineq:Asymp-01}
K_\OM(\varphi\circ\wt\sigma(t),\eta(t))\leq C_0
\end{equation}
where $\sigma:[0,\infty) \to \D$ is a radial line joining $0$ and $1$, and $S>0$ such that $\wt{\sigma}:=\left.\sigma\right|_{[S, \infty)}\subset U_1$. Up to reparameterization, we may assume that $\wt{\sigma}:[0, \infty) \to \D$.
Note that we have also obtained an injective holomorphic map $\varphi: \D \to \OM$ obtained by $x$ and $(z_n)_{n \geq 1}$ as in Condition~(M$\star$). By Lemma~\ref{Lem:Geo_of_Horo_can_be_truncated_and_it_is_asyp_wrt_limit_of_geo_seq}, there exist $T>0$, a constant $C^{\prime}>0$ and $C_1>0$ such that $\wh\gamma\defeq\left.\gamma\right|_{[T,\infty)}\subset \varphi(U_1)$ and such that for all $t\geq 0$
\begin{equation}\label{Ineq:Main_asymp_ineq-01}
K_\OM(\varphi\circ\wt{\sigma}(t), \wh\gamma(t))\leq C^\prime.
\end{equation}
Let $\Gamma$ be a component of $\bdy\OM$ such that $x \in \Gamma$. Let $\OM_{\Gamma}$ be simply connected domain obtained by filling all components of $\C\setminus\OM$ except which contains $\Gamma$. Therefore,  By Riemann mapping theorem, there exists a biholomorphism $h:\D \to \OM_{\Gamma}$. Let $h_{1} : \D \to \K$ be a biholomorphism given by $h_{1}(z) = \frac{1+z}{1-z}$. Let $\wt{h}:= h\circ h^{-1}_{1}: \K \to \OM_{\Gamma}$. We now make the following claim.
\smallskip

\noindent{\bf Claim~1.} $\wh{\gamma}$ is a $(1,\kappa)$-quasi-geodesic of $(\OM_{\Gamma}, K_{\OM_{\Gamma}})$. 
\smallskip

\noindent{\itshape Proof of Claim~1.} Since $\OM_\Gamma$ satisfies Condition~(M$\star$), so it follows from Lemma~\ref{Lem:Vis_of_image_of_bigeodesic_region_and_image_of_unit_disc} and Lemma~\ref{Lem:Visi_of_images_of_bi-geodesic_reg_inside_U} that $\varphi(U_1)$ and $\OM_\Gamma\setminus\varphi(U_2)$, $\varphi(U_2)$ and $\OM_\Gamma\setminus\varphi(U_3)$, and $\varphi(U_3)$ and $\OM_\Gamma\setminus\varphi(\D)$ have visibility property with respect to $(\OM_\Gamma, K_{\OM_\Gamma})$. Note that $\wh\gamma$ is a $(1,\kappa_0)$-quasi-geodesic of $(\varphi(\D),K_{\varphi(\D)})$. Thus, by Theorem~\ref{Thm:Global_Loc_Res}, there exists a constant $C_2>0$ such that for all $s, t\in [0,\infty)$
\begin{equation}\label{ineq:Localization_Ineq_in_Main_Result-03}
K_{\varphi(\D)}(\wh{\gamma}(s), \wh{\gamma}(t)) \leq K_{\OM_{\Gamma}}(\wh{\gamma}(s), \wh{\gamma}(t)) + C_2.
\end{equation}
Using~\eqref{ineq:Localization_Ineq_in_Main_Result-03}, we get that for all $s, t \in [0,\infty)$
\[
    |s-t|-\kappa_0 - C_2 \leq K_{\OM_{\Gamma}}(\wh{\gamma}(s), \wh{\gamma}(t)) \leq |s-t| + \k_0+ C_2.
\]
Therefore, $\wh{\gamma}$ is a $(1,\kappa)$-quasi-geodesic of $(\OM_{\Gamma}, K_{\OM_{\Gamma}})$ for $\k:=\k_0+C_2$.
\hfill $\blacktriangleleft$
\smallskip

By Claim~1, $\wt{h}^{-1}\circ \wh{\gamma}$ is a $(1,\kappa)$-quasi-geodesic of $(\K, K_{\K})$. Let $(t_k)_{k\geq 1}$ be a sequence with $t_k\to\infty$. Without loss of generality, we may assume, up to subsequence, $(h^{-1}\circ \wh{\gamma}(t_k))_{k\geq 1}$ converges to $1$. Then $(h_{1}\circ h^{-1} \circ \wh{\gamma}(t_k))_{k\geq 1}$ converges to $\infty$. That is, $(\wt{h}^{-1}\circ \wh{\gamma}(t_k))_{k\geq 1}$ converges to $\infty$. Let $L_{\R}: [0,\infty) \to \K$ be a geodesic of $(\K, K_{\K})$ given by $L_{\R} = \{z\in\R:\rprt z\geq a\}$ where $a>0$. Since $(\K, K_{\K})$ is Gromov hyperbolic space, there exists a constant $C_3>0$ such that for all $t \geq 0$
\[
K_{\K}(L_{\R}(t), \wt{h}^{-1}\circ \wh{\gamma}(t))\leq C_3.
\]
Therefore, for all $t \in [0,\infty)$
\begin{equation}\label{ineq:Asymp_Condition_for_quasi_geo_in_filled_domain_5}
    K_{\OM_{\Gamma}}(\wt{h}\circ L_{\R}(t), \wh{\gamma}(t)) \leq C_3.
\end{equation}

Let $U_0\subsetneq U_1$ be a bi-geodesic region. Since $\wh\gamma\subset \varphi(U_1)$, we may assume that $\wh\gamma(t_k)\in \varphi(U_0)$ for all $k\geq 1$. We make the following claim
\smallskip
        
\noindent{\bf Claim~2.} There exists $M_0>a$ such that $\wt{h}\circ L_{\R}([M_0, \infty))\subset \varphi(U_1)$. 
\smallskip
        
\noindent{\itshape Proof of Claim~2.} Assume that the claim is not true. There exists a $(s_k)_{k \geq 1}$ with $s_k \to \infty$ as $k \to \infty$ such that $\wt{h}(L_{\R}(s_k)) \in \OM_{\Gamma}\setminus \varphi(U_1)$ for all $k \geq 1$. Suppose that $\alpha_k:[0,b_{k}] \to \OM_{\Gamma}$ is a geodesic of $(\OM_{\Gamma}, K_{\OM_{\Gamma}})$ with $\alpha_{k}(0) = \wh\gamma_{k}(t_{k})$ and $\alpha_{k}(b_k) = \wt{h}(L_{\R}(s_k))$. Since $\varphi(U_0)$ and $\OM_{\Gamma}\setminus\varphi(U_1)$ have visibility property with respect to $(\OM_{\Gamma}, K_{\OM_{\Gamma}})$, there exists a compact set $K\subset \OM_{\Gamma}$ such that $\alpha_{k}\cap K\neq \emptyset$. Let $o_k \in \alpha_{k}\cap K$. Hence, for all $k \geq 1$
\[
 K_{\OM_\Gamma}(\wh\gamma(t_k), \wt{h}(L_{\R}(s_k))) = K_{\OM_{\Gamma}}(\wh\gamma(t_k), o_k) + K_{\OM_{\Gamma}}(o_k, \wt{h}(L_{\R}(s_k))).
\]
Then, $K_{\OM_{\Gamma}}(\wh\gamma(t_k), o_k) < K_{\OM_\Gamma}(\wh\gamma(t_k), \wt{h}(L_{\R}(s_k)))$ for all $k \geq 1$. Note that $\delta_{\OM_\Gamma}(\wh\gamma(t_k)) \to 0$ or $\|\wh\gamma(t_k)\|\to\infty$ as $k \to \infty$ and $K$ is compact set, so it follows that $K_{\OM_{\Gamma}}(\wh\gamma(t_k), o_k) \to \infty$ as $k \to \infty$. Therefore, $K_{\OM_\Gamma}(\wh\gamma(t_k), \wt{h}(L_{\R}(s_k))) \to \infty$ as $k \to \infty$ which contradicts \eqref{ineq:Asymp_Condition_for_quasi_geo_in_filled_domain_5}. This completes the proof of claim.
\hfill $\blacktriangleleft$
\smallskip

Let $\wt{L}_{\R}:= \left.L_{\R}\right|_{[M_0,\infty)}$. We may assume, by reparameterization, $\wt{L}_{\R}:[0,\infty) \to \K$. Since $\left.L_{\R}\right|_{[a,M_0]}$ is compact, so  by~\eqref{ineq:Asymp_Condition_for_quasi_geo_in_filled_domain_5}, there exists $C^{\prime\prime}>0$ such that for all $t \geq 0$
\begin{equation}\label{ineq:Asymp_Condition_for_quasi_geo_06}
    K_{\OM_{\Gamma}}(\wt{h}\circ \wt{L}_{\R}(t), \wh{\gamma}(t)) \leq C^{\prime\prime}.
\end{equation}
Therefore, using~\eqref{Ineq:Main_asymp_ineq-01}, \eqref{ineq:Asymp_Condition_for_quasi_geo_06} and triangle inequality, for all $t \in [0,\infty)$
\[
K_{\OM_{\Gamma}}(\wt{h}\circ\wt{L}_{\R}(t), \varphi\circ \wt{\sigma}(t)) \leq C_1 + C^{\prime\prime} =: C_4.
\]
By Claim~1, we have $\wt h\circ \wt L_\R\subset \varphi(U_1)$. Note that $\varphi(U_1)$ and $\OM_\Gamma\setminus\varphi(U_2)$, $\varphi(U_2)$ and $\OM_\Gamma\setminus\varphi(U_3)$, and $\varphi(U_3)$ and $\OM_\Gamma\setminus\varphi(\D)$ have visibility property with respect to $(\OM_\Gamma, K_{\OM_\Gamma})$. Since $\varphi\circ\wt\sigma\subset \varphi(U_1)$, by Theorem~\ref{Thm:Geo_horosphere_of_small_domain_is_subset_of_geo_horosphere_of_big_domain}, there exists a constant $C_5>0$ such that for all $R>0$
\begin{equation}\label{ineq:Geo_horosphere_of_small_domain_is_subset_of_geo_horosphere_of_big_domain-1}
    H_{\varphi(\D)}( \varphi\circ \wt{\sigma}, R)\subseteq H_{\OM_{\Gamma}}(\wt{h}\circ\wt{L}_{\R}, C_5R).
\end{equation}
Fix $R>0$. Then, by Theorem~\ref{Thm:Geo_horosphere_of_big_domain_is_subset_of_geod_of_small_domain_upto_some_compact_set_and-constant}, there exist a compact set $K\subset \OM_{\Gamma}$ and a constant $C_6>0$ such that
\begin{equation}\label{ineq:Geo_horosphere_of_big_domain_is_subset_of_geod_of_small_domain_upto_some_compact_set_and_constant-0}
    \wt{h}(H_{\K}(\wt{L}_{\R}, R))\setminus K \subseteq H_{\varphi(\D)}( \varphi\circ \wt{\sigma}, C_6R).
\end{equation}
Therefore,
\[
H_{\K}(\wt{L}_{\R}, R) \setminus \wt{h}^{-1}(K) \subseteq \wt{h}^{-1}(H_{\varphi(\D)}( \varphi\circ \wt{\sigma}, C_6R)).
\]
By Lemma~\ref{Lem:Relation_between_origin_geo_horo_and_truncated_geo_horo}, we get
\begin{equation}\label{ineq:Geo_horosphere_of_big_domain_is_subset_of_geod_of_small_domain_upto_some_compact_set_and_constant-1}
H_{\K}(L_{\R}, e^{-2M_0}R) \setminus \wt{h}^{-1}(K) \subseteq \wt{h}^{-1}(H_{\varphi(\D)}( \varphi\circ \wt{\sigma}, C_6R)).
\end{equation}
Fix $M>0$. Let $T_{M}: \K \to \K$ given by $T_{M}(z)= z+ i +M$. Since $\wt{h}^{-1}(K)$ is a compact set, for sufficiently large $M>0$
\[
T_{M}(\K) \subseteq H_{\K}(L_{\R}, e^{-2S}R)\setminus \wt{h}^{-1}(K).
\]
Therefore, by~\eqref{ineq:Geo_horosphere_of_big_domain_is_subset_of_geod_of_small_domain_upto_some_compact_set_and_constant-1}, for sufficiently large $M>0$
\begin{equation}\label{ineq:TM(K)_is_subset_of_preimage_of_horo_wrt_analy_disc-1}
T_{M}(\K) \subseteq \wt{h}^{-1}(H_{\varphi(\D)}( \varphi\circ \wt{\sigma}, C_6R)).
\end{equation}
Note that $H_{\varphi(\D)}( \varphi\circ \wt{\sigma}, C_6R) \subseteq \varphi(\D)$, so by~\eqref{ineq:TM(K)_is_subset_of_preimage_of_horo_wrt_analy_disc-1}, $T_{M}(\K)\subseteq \wt{h}^{-1}(\varphi(\D))$ for sufficiently large $M>0$.
\smallskip

Fix $M>0$ such that $T_{M}(\K)\subseteq \wt{h}^{-1}(\varphi(\D))$. Let $D:= \wt{h}^{-1}(\OM)$ and $\wt{h}_D:= \left.\wt{h}\right|_D$. Define $\psi: \OM \to \OM$ by $\psi:= \wt{h}_D\circ T_{M}\circ \wt{h}_D^{-1}$. Therefore, $\psi$ is a well defined injective holomorphic self map on $\OM$ such that $\psi(\OM) \subseteq \varphi(\D)$.
\smallskip

Since $\OM = \wt{h}_D(D)$, $\varphi(\D) \subseteq \wt{h}_D(D)$. For sufficiently small $R^{\prime}>0$ with $R^{\prime}< R$,
\[
    H_{\K}(L_{\R}, R^{\prime}) \subseteq H_{\K}(L_{\R}, e^{-2M_0}R) \setminus \wt{h}(K).
\]
Note that by~\eqref{ineq:Geo_horosphere_of_big_domain_is_subset_of_geod_of_small_domain_upto_some_compact_set_and_constant-1}, $H_{\K}(L_{\R}, e^{-2S}R) \setminus \wt{h}(K) \subseteq \wt{h}^{-1}_D(\varphi(\D))$, and so $H_{\K}(L_{\R}, R^{\prime}) \subseteq \wt{h}^{-1}_D(\varphi(\D))$. Hence, $H_{\K}(L_{\R}, R^{\prime}) \subseteq D$. Therefore, both $\wt{h}$ and $\wt{h}_D$ make sense in $H_{\K}(L_{\R}, R^{\prime})$.
\smallskip
        
Let $T:\K \to \K$ given by $T(z) = z+ i$. Let $\phi: \OM_{\Gamma} \to \OM_{\Gamma}$ defined by $\phi:= \wt{h}\circ T\circ \wt{h}^{-1}$. Note that $\left.\phi\right|_\OM$ need not be self map of $\OM$. Let $f: \OM \to \OM$ given by $f:=\psi^{-1}\circ \phi \circ \psi$. We need to show that $f$ is a self map of $\OM$. It suffices to show that $\phi(\psi(\OM)) \subset \psi(\OM)$ for the same. By~\eqref{ineq:TM(K)_is_subset_of_preimage_of_horo_wrt_analy_disc-1}, $T_{M}(D) \subseteq \wt{h}^{-1}_D(H_{\varphi(\D)}( \varphi\circ \wt{\sigma}, C_6R))$. Since $\wt{h}_D(T_{M}(D)) \subseteq H_{\varphi(\D)}( \varphi\circ \wt{\sigma}, C_6R)$, so $\psi(\OM)\subseteq H_{\varphi(\D)}( \varphi\circ \wt{\sigma}, C_6R)$. Now, using~\eqref{ineq:Geo_horosphere_of_small_domain_is_subset_of_geo_horosphere_of_big_domain-1} and Lemma~\ref{Lem:Relation_between_origin_geo_horo_and_truncated_geo_horo}, we get  
\begin{align*}
    \phi(\psi(\OM)) &\subseteq \phi(H_{\varphi(\D)}( \varphi\circ \wt{\sigma}, C_5R))\\
                    &\subseteq \phi(\wt{h}(H_{\K}(\wt{L}_{\R}, C_5C_6R)))\\
                    &\subseteq \phi(\wt{h}(H_{\K}(L_{\R}, e^{-2S}C_5C_6R)))\\
                    &\subseteq \wt{h}(T(H_{\K}(L_{\R}, e^{-2S}C_5C_6R)))\\
                    &\subseteq \wt{h}(H_{\K}(L_{\R}, e^{-2S}C_5C_6R))\\
                    &\subseteq \varphi(\D).
\end{align*}
where last inclusion follows from~\eqref{ineq:Geo_horosphere_of_big_domain_is_subset_of_geod_of_small_domain_upto_some_compact_set_and_constant-0} and second last inclusion follows from the fact that $T(H_\K(L_{\R}, e^{-2S}C_5C_6R))=H_\K(L_{\R}, e^{-2S}C_5C_6R)$. Hence, $f$ is a well defined holomorphic self map of $\OM$. Let $p \in \OM$. Note that $f^{n}(p) = \wt{h}_D(\wt{h}^{-1}_D(p)+in)$. Since $\wt{h}^{-1}_D(p)+in \to \infty$ as $n \to \infty$ and $\wt{h}_D$ is a proper map, $\delta_{\OM}(\wt{h}(\wt{h}^{-1}_D(p)+in)) \to 0$ or $\|\wt{h}(\wt{h}^{-1}_D(p)+in)\|\to\infty$ as $n \to \infty$. Therefore, $K_\OM(f^{n}(p), p) \to \infty$ as $n\to\infty$. Thus, for each $p\in \OM$ the orbit $(f^n(p))_{n\geq1}$ is not relatively compact on $\OM$.

Therefore, by hypothesis, there exists $\xi_0 \in \bdy_{\infty} \OM$ such that $(f^{n})_{n \geq 1}$ converges uniformly on compact subset of $\OM$ to $\xi_0$. Note that $\xi_0$ can be equal to $\infty$. Let $z_{0} \in H_{\K}(L_{\R}, R^{\prime}) \cap \R$. Using the fact that $\wt{h}_D = \left.\wt{h}\right|_D$ and $T^{-1}_{M}(w) = w - i -M$, for all $n \in \N$, we obtain
\begin{align*}
    f^{n}(h_D(z_{0})) &= \wt{h}_D(T_{M}^{-1}(T^{n}(T_{M}(z_{0}))))\\
                        &= \wt{h}_D(T_{M}^{-1}(T^{n}(z_{0}+i + M)))\\
                        &= \wt{h}_D(T_{M}^{-1}(z_{0} + i(n + 1) + M)\\
                        &= \wt{h}_D(z_{0} + in).
\end{align*}
Hence, $\wt{h}_D( z_{0} + in) \to \xi_0$ as $n \to \infty$. Let $R_{0}>0$ such that $R_{0} < R^{\prime}$. Therefore, denoting by $[t]$ and $\{t\}$ the integer part and fractional part of a real number $t$, respectively, we get
\begin{align*}
    K_{\OM_{\Gamma}}(\wt{h}(R_{0}+it), \wt{h}(R_{0}+i[t])) &= K_\K(R_{0}+it, R_{0}+i[t])\\
                                                            &= K_\K(T^{-[t]}(R_{0}+it), T^{-[t]}(R_{0}+i[t]))\\
                                                            &= K_\K(R_{0}+it - i[t], R_{0}+i[t] - i[t])\\
                                                            &= K_\K(R_{0}+i\{t\} , R_{0})
\end{align*}
where second equality follows from the fact that $T^{-1}(z)= z-i$ is an automorphism of $\K$. Therefore, for all $t \geq 0, K_{\OM_{\Gamma}}(\wt{h}(R_{0}+it), \wt{h}(R_{0}+i[t])) < \infty$. Suppose $(t_{n})_{n \geq 1} \subset [0,\infty)$ is a sequence that converges to $\infty$. Hence, for all $n \in \N$
\begin{equation}\label{ineq:K-distance_images_of_sequences_in_vertical_line_is_fin-1}
K_{\OM_{\Gamma}}(\wt{h}(R_{0}+it_{n}), \wt{h}(R_{0}+i[t_{n}])) = K_\K(R_{0}+i\{t_n\} , R_{0}) < \infty.
\end{equation}
Since $\wt{h}_{0}(R_{0} + in) \to \xi_0$ as $n \to \infty$, so $\wt{h}(R_{0}+i[t_{n}])\to\infty$ as $n\to\infty$. Since $R_{0}+it_{n} \to \infty$ as $n \to \infty$ and $\wt h$ is a proper map, by~\eqref{ineq:K-distance_images_of_sequences_in_vertical_line_is_fin-1}, we get that $\wt{h}(R_{0}+it_{n}) \to \xi_0$ as $n \to \infty$. Therefore, 
\begin{equation}\label{eq:Limit_of_h_along_vertical_line_exists-1}
\lim_{t\to\infty}\wt{h}(R_{0}+it) = \xi_0.
\end{equation}
By Result~\ref{Res:Lehto-Virtanen_Thm} for $\K$, the non-tangential limit of $\wt h$ at $\infty$ exists and it is equal to $\xi_0$. In particular
\begin{equation}\label{eq:Limit_of_h_along_horizontal_line_exists-1}
\lim_{r \to \infty}\wt{h}(R_{0}+r) = \xi_0.
\end{equation}   
\noindent{\bf Claim~4.} For every sequence $(z_{n}^{\prime})_{n \geq 1}\subset D$ converging to $\infty$ with $\rprt(z_{n}^{\prime}) \geq R_{0}$ and $\iprt z_{n}^{\prime}\geq 0$,
\[
\lim_{n \to \infty}\wt{h}(z_{n}^{\prime}) = \xi_0.
\]
\noindent{\itshape Proof of Claim~4.} Let $\Sigma$ be the Jordan curve given by $\clos{\{R_{0}+it:t\geq0\}}^{\infty}\cup\clos{\{R_{0}+r:r\geq 0\}}^{\infty}$. Let $D^\prime$ be a domain in $\C_{\infty}$ bounded by $\Sigma$ which contains $R_{0}+1+i$. Since $H_{\K}(L_{\R}, R^{\prime}) \subseteq D$, $\clos{D^\prime}^{\infty}\setminus \{\infty\} \subseteq D$. By~\eqref{eq:Limit_of_h_along_vertical_line_exists-1} and~\eqref{eq:Limit_of_h_along_horizontal_line_exists-1}, $\wt{h}(\Sigma)$ is a Jordan curve in $\C_{\infty}$. Note that one of the connected component of $\C_{\infty}\setminus\wt{h}(\Sigma)$ is $\wt{h}(D^\prime)$. Since $D^\prime$ and $\wt{h}(D^\prime)$ are Jordan domains in $\C_{\infty}$ and $\left.\wt{h}\right|_{D^\prime}: D^\prime\to \wt{h}(D^\prime)$, by the Carath{\'e}odory Extension Theorem, by the abuse of notation, $\wt{h}: \clos{D^\prime}^\infty \to \clos{\wt{h}(D^\prime)}^{\infty}$ is a homeomorphism. Note that $\wt{h}^{-1}(\xi_0) = \infty$. Therefore, for every $(z_{n}^{\prime})_{n \geq 1}$ converging to $\infty$ with $\rprt(z_{n}^{\prime}) \geq R_{0}$ and $\iprt z_{n}^{\prime} \geq 0, \lim_{n \to \infty}\wt{h}(z^{\prime}_{n})=\xi_0$.
\hfill $\blacktriangleleft$
\smallskip

Let $g: \OM \to \OM$ given by $g:= \psi^{-1}\circ\phi^{-1}\circ\psi$. Note that, by arguing same as in the case of $f$, $g$ is well defined holomorphic self map of $\OM$. Moreover, we repeat the same construction for $g$ and we show that there exists $\wt{\xi}_0 \in \partial \OM$ such that $\lim_{n \to \infty}\wt{h}(z_{n}^{\prime}) = \wt{\xi}_0$ for all $(z_{n}^{\prime})_{n \geq 1}$ converging to $\infty$ with $\rprt z_{n}^{\prime} \geq R_{0}$ and $\iprt z_{n}^{\prime} \leq 0$. Note that $H_\K(L_{\R}, R^{\prime}) \subset D$, so there exists $N_{0} \in \N$ such that $n \in D$ for all $n \geq N_{0}$. Hence, by taking $z_{n}^{\prime}=n$, we get that $\xi_0 = \wt{\xi}_0$. Therefore, 
\begin{equation}\label{eq:H-limit_of_h_at_infinity_exists-1}
            \lim_{n \to \infty}\wt{h}(z_{n}^{\prime}) = \xi_0
\end{equation}
for all $(z_{n}^{\prime})_{n \geq 1}$ converging to $\infty$ with $\rprt z_{n}^{\prime} \geq R_{0}$. 
\smallskip

By~\eqref{eq:H-limit_of_h_at_infinity_exists-1}, for all $R>0$, we get
\begin{equation}\label{ineq:Closed_geo_horo_intersect_bdry_single_point_2}
\clos{H_{\OM_{\Gamma}}(\wt{h}\circ L_{\R},R)}^{\infty}\cap \bdy_{\infty} \OM_{\Gamma} = \{\xi_0\}.
\end{equation}

By Theorem~\ref{Thm:Geo_horosphere_of_big_domain_is_subset_of_geod_of_small_domain_upto_some_compact_set_and-constant}, there exist a constant $C>0$ and a compact set $K_0\subset\OM$ such that 
\[
H_\OM(\wh\gamma,R)\setminus K_0\subseteq H_{\varphi(\D)}(\varphi\circ\wt\sigma, CR).
\]
Using~\eqref{ineq:Geo_horosphere_of_small_domain_is_subset_of_geo_horosphere_of_big_domain-1}, we get
\begin{equation}\label{ineq:Geo_horosphere_of_small_domain_is_subset_of_geo_horosphere_of_big_domain-4}
    H_{\OM}(\wh\gamma, R)\setminus K_0\subseteq H_{\OM_{\Gamma}}(\wt{h}\circ\wt{L}_{\R}, CC_5R).
\end{equation}
If $\xi_0=\infty$, then using~\eqref{ineq:Closed_geo_horo_intersect_bdry_single_point_2} and~\eqref{ineq:Geo_horosphere_of_small_domain_is_subset_of_geo_horosphere_of_big_domain-4}, we get that $\clos{H_{\OM}(\gamma, e^{-2T}R)}\cap \partial \OM = \emptyset$ which cannot be possible because $x\in \clos{H_{\OM}(\gamma, e^{-2T}R)}\cap \partial \OM$. If $\xi_0\in\bdy\OM$, then using the fact that $\partial\OM_{\Gamma} \subset \partial\OM$, \eqref{ineq:Closed_geo_horo_intersect_bdry_single_point_2} and \eqref{ineq:Geo_horosphere_of_small_domain_is_subset_of_geo_horosphere_of_big_domain-4}, we get
\[
\clos{H_{\OM}(\gamma, e^{-2T}R)}\cap \partial \OM\subset \clos{H_{\OM_{\Gamma}}(\wt{h}\circ L_{\R}, e^{-2S}CR)}\cap \partial \OM_{\Gamma} = \{\xi_0\}.
\]
Since $x\in\clos{H_{\OM}(\gamma, e^{-2T}R)}\cap \partial \OM$, so $x=\xi_0$. Thus, $\clos{H_{\OM}(\gamma, e^{-2T}R)}\cap \partial \OM=\{x\}$. This completes the proof of $(1) \Rightarrow (2)$. This completes the proof of Theorem.
\end{proof}

\begin{proof}[Proof of {\Cref{Cor:Finitely_conn_dom_charac_for_DWP}}]
We shall prove that $\OM$ satisfies Condition $(\star\star)$. By Lemma~\ref{Lem:Contions_double_star_implies_mild_star}, it follows that $\OM$ satisfies Condition~(M$\star$). The proof is then complete by Theorem~\ref{Thm:Characterization_of_Condition(*)}. Let $\Gamma_1,\ldots,\Gamma_n$ be all connected components of the boundary of domain $\OM$. Let $\OM_{\Gamma_1}$ be simply connected domain obtained by filling all the components of $\C\setminus\OM$ except which contains $\Gamma_1$. Then, by the Riemann mapping theorem, there exists a biholomorphism $f:\D\to \OM_{\Gamma_1}$. There exists a compact set $K\subset\D$ such that $f^{-1}(\Gamma_2), \ldots, f^{-1}(\Gamma_n) \subset K$. Therefore, $\Gamma_2,\ldots,\Gamma_n\subset f(K)$. Therefore, there exists a neighbourhood $V$ of $\Gamma_1$ such that $V\cap \OM$ is doubly connected. Hence, $f^{-1}(V\cap\OM)$ is doubly connected. Without loss of generality, we may assume that $f^{-1}(V\cap\OM)=A(0,r,1)$ for some $r>0$. Let $f_0\defeq\left.f\right|_{A(0,r,1)}:A(0,r,1)\to\OM$. Hence, we get an injective holomorphic map $f_0:A(0,r,1)\to\OM$ such that $f_0(A(0,r,1))=V\cap\OM$ and such that for every sequence $(\zeta_n)_{n\geq 1}\in A(0,r,1)$ with $|\zeta_n|\to 1$ as $n\to\infty$, we have either $\delta_\OM(f_0(\zeta_n))\to 0$ as $n\to\infty$ or  $|f_0(\zeta_n)|\to \infty$ as $n\to\infty$. This completes the proof.  
\end{proof}

\begin{proof}[Proof of {\Cref{Cor:Well_seperated_domains_charac}}]
We shall prove that $\OM$ satisfies Condition $(\star\star)$. By Lemma~\ref{Lem:Contions_double_star_implies_mild_star}, it follows that $\OM$ satisfies Condition~(M$\star$). The proof is then complete by Theorem~\ref{Thm:Characterization_of_Condition(*)}. Let $S\defeq\{\Gamma_i:i\in I\}$ be the set of all connected components of $\bdy\OM$. Fix a connected component $\Gamma_1$. Let $\OM_{\Gamma_1}$ be simply connected domain obtained by filling all the components of $\C\setminus\OM$ except which contains $\Gamma_1$. It is clear that $\bdy\OM_{\Gamma_1}=\Gamma_1$. By the Riemann mapping theorem, there exists a biholomorphism $f:\D\to \OM_{\Gamma_1}$. Since there exists $\delta>0$ such that $\distance(\Gamma_i,\Gamma_j)>\delta$ for all $i\neq j$, so it follows that $\clos{N_{\C}(\Gamma_1,\delta/2)}\cap \Gamma_j=\emptyset$ for all $j\in I$ where $\distance(\Gamma_i,\Gamma_j)$ denotes the Euclidean distance of the set $\Gamma_i$ to the set $\Gamma_j$. Since boundary components of domain $\OM$ are bounded, so $\Gamma_1$ is compact set. Hence, by Result~\ref{Res:Zoretti_Thm}, there exists a Jordan region $V$ in $\C$ such that $\Gamma_1 \subset V$, $\bdy V\cap \Gamma_1=\emptyset$ and $\bdy V\subset N_{\C}(\Gamma_1, \delta/2)$. It is clear that, $V\cap\OM$ is doubly connected. Therefore, $f^{-1}(V\cap\OM)$ is doubly connected subset of the unit disc. Without loss of generality, we may assume that $f^{-1}(V\cap\OM)=A(0,r,1)$ for some $r>0$. Let $f_0\defeq\left.f\right|_{A(0,r,1)}:A(0,r,1)\to\OM$. Hence, we get an injective holomorphic map $f_0:A(0,r,1)\to\OM$ such that $f_0(A(0,r,1))=V\cap\OM$ and such that for every sequence $(\zeta_n)_{n\geq 1}\in A(0,r,1)$ with $|\zeta_n|\to 1$ as $n\to\infty$, we have either $\delta_\OM(f_0(\zeta_n))\to 0$ as $n\to\infty$ or  $|f_0(\zeta_n)|\to \infty$ as $n\to\infty$. Since $\Gamma_1$ is arbitrary, so the proof is complete.
\end{proof}

We now prove the following theorem which is required to proof of Theorem~\ref{Thm:Main_thm_outside_totally_disc_set}.
\begin{theorem}\label{Thm:Strengthening_closed_geo_horo_inter_at_single}
Let $\OM \subset\C$ be a hyperbolic domain. Suppose that there exists a totally disconnected set $S \subset \bdy\OM$ such that $\OM$ satisfies Condition $(M \star)$ for every $ x \in \partial\OM\setminus S$. The following are equivalent:
\begin{enumerate}
\item for every sequence $(z_{n})_{n \geq 1}$ in $\OM$ that converges to some boundary point $y \in \bdy\OM\setminus S$ such that if $(\gamma_{k})_{k\geq 1}$ is a sequence of geodesic $(\OM,K_\OM)$ that joins some $p \in \OM$ and $z_{n_{k}}$ such that $(\gamma_{k})_{k\geq 1}$, up to a subsequence, converges uniformly on compact set to a geodesic ray $\gamma$ of $(\OM, K_\OM)$; where $\varphi:\D \to \OM$ is an injective holomorphic map and $(z_{n_k})_{k \geq 1}$ is a subsequence of $(z_{n})_{n \geq 1}$ as in Condition $(M \star)$ obtained by $(z_{n})_{n \geq 1}$ and $y$, then either $\clos{H_{\OM}(\gamma, R)}\cap \bdy\OM =\emptyset$ for all $R>0$ or there exists $x \in \bdy\OM$ such that for all $R>0$
\[
\clos{H_{\OM}(\gamma, R)}\cap \bdy\OM = \{x\};
\]
\item for every geodesic ray $\gamma:[0, \infty) \to \OM$ of $(\OM, K_\OM)$, either $\clos{H_{\OM}(\gamma, R)}\cap \bdy\OM =\emptyset$ for all $R>0$ or there exists $x \in \bdy \OM$ such that for all $R>0$
\[
\clos{H_{\OM}(\gamma, R)}\cap \bdy\OM = \{x\}.
\]
\end{enumerate} 
\end{theorem}
\begin{proof}
It is trivial that $(2) \implies (1)$. We want to prove $(1) \implies (2)$ by the method of contradiction. Let $R>0$ and let $\gamma:[0,\infty)\to\OM$ be a geodesic of $(\OM,K_\OM)$. If $\clos{H_{\OM}(\gamma, R)}\cap \bdy\OM=\emptyset$, then we are done. Assume that $\clos{H_{\OM}(\gamma, R)}\cap \bdy\OM\neq\emptyset$. If it is singleton, then we are done. Assume that $\xi_1,\xi_2\in\clos{H_{\OM}(\gamma, R)}\cap \bdy\OM$ such that $\xi_1\neq\xi_2$. Let $(s_n)_{n\geq 1}\subset [0, \infty)$ be a sequence such that $s_{n} \to \infty$. Now, up to a subsequence, $\gamma(s_{n}) \to \xi_0\in\bdy_{\infty}\OM$. If $\xi_0\in\bdy\OM$, then, without loss of generality we may assume that $\xi_0\neq\xi_1$. There exists a sequence $(z_n)_{n\geq 1}\subset H_{\OM}(\gamma, R)$ such that $z_n\to \xi_1$. Let $\gamma_n:[a_n,b_n]\to\OM$ be a geodesic of $(\OM,K_\OM)$ such that $\gamma_n(a_n)=z_n$ and $\gamma_n(b_n)=\gamma(s_n)$. We now make the following claim.
\smallskip

\noindent{\bf Claim.} $(\gamma_{n})_{n \geq 1}$ eventually avoids every compact subset of $\OM$.
\smallskip

\noindent{\itshape Proof of Claim.} The proof is similar to the proof of Lemma~\ref{Lem:Tail_of_Seq_cannot_lie_in_geodesic_horosphere}. Assume that there exist a compact set $K\subset\OM$ and $N\in\N$ such that $\gamma_n\cap K\neq\emptyset$ for all $n\geq N$. Let $o_{n, j} \in \gamma_{n}\cap K$. Then, up to a subsequence, $(o_{n,j})_{j\in\N}$ converges to $o_n$ such that $o_n\in K$ as $n\in\N$. Let 
\[
S_{n,j}:=K_\OM(o_{n,j},\gamma(s_j))-K_\OM(\gamma(0),\gamma(s_j))
\]
Since $K$ is compact set and $o_{n,j}\in K$, it follows that for all $j \in \N, |S_{n,j}|$ is bounded above, say by a constant $C>0$. Therefore, for all $j \geq 1$ and for all $N \in \N$
\begin{align*}
K_\OM(z_{n}, o_{n,j}) &= K_\OM(z_{n}, \gamma(s_j)) - K_\OM(\gamma(0), \gamma(s_j)) - K_\OM(o_{n,j},\gamma(s_j))+K_\OM(\gamma(0),\gamma(s_j))\\
&\leq K_\OM(z_{n}, \gamma(s_j)) - K_\OM(\gamma(0), \gamma(s_j)) + |S_{n,j}|
\end{align*}
Therefore, by taking $j \to \infty$, we get for all $n \in \N$
\[
K_\OM(z_{n}, o_{n}) < \frac{1}{2}\log R +C.
\]
Since $o_n\in K$ for all $n\in N$ and $z_n\to\xi_1$, so we get a contradiction.
\hfill $\blacktriangleleft$
\smallskip

Note that if $\xi_0=\infty$, then for all $n\in\N$ we may obtain $\wt{\gamma}_n$ which is a truncation of $\gamma_n$ such that $(\wt\gamma_n)_{n\geq 1}$ eventually avoids every compact subset of $\OM$, which can be proved exactly same as the Claim above, and such that the sequence of the end point of $\wt{\gamma}_n$ other than $\gamma_n(a_n)$ converge to $\xi_3\in\bdy\OM$ where $\xi_1\neq\xi_3$. Therefore, without loss of generality, we work with $(\gamma_n)_{n\geq 1}$. Using Claim above and Remark~\ref{Rem:Conv_res_of_geo}, there exist a sequence $(t_n)_{n\geq 1}$ with $t_n\in [a_n,b_n]$ for all $n\geq 1$ and $\xi\in\bdy\OM\setminus S$ such that $\gamma_n(t_n)\to \xi$ as $n\to \infty$. Since $\xi \in \bdy\OM\setminus S$, we get an injective holomorphic map $\varphi:\D \to \OM$, a bi-geodesic region $U$ and a subsequence $(\gamma_{n_{k}}(t_{n_{k}}))_{k \geq 1} \subset U$ such that $\varphi^{-1}(\gamma_{n_{k}}(t_{n_{k}})) \to 1$ as $k \to \infty$. By Proposition~\ref{Prop:Vis_of_all_req_sets_which_works_everywhere}, there exist bi-geodesic regions $U_1\subsetneq U_2\subsetneq U_3 \subsetneq U$ such that $\varphi(U_1)$ and $\OM\setminus\varphi(U_2)$, $\varphi(U_2)$ and $\OM\setminus\varphi(U_3)$, and $\varphi(U_3)$ and $\OM\setminus\varphi(\D)$ have visibility property with respect to $(\OM, K_\OM)$. By Lemma~\ref{Lem:Seq_of_geo_lies_in_image_of_bi-geo_region}, there exists $N \in \N$ such that $\gamma_n \subset \varphi(U_{1})$ for all $n \geq N$ and both sequences $(\varphi^{-1}(\gamma(s_n)))_{n \geq N}$ and $(\varphi^{-1}(z_n)_{n \geq N}$ converge to $1$. Therefore, we may assume $\gamma_{n} \subset \varphi(U_{1})$ for all $n \geq 1$ and both sequences $(\varphi^{-1}(\gamma(s_n)))_{n \geq 1}$ and $(\varphi^{-1}(z_n))_{n \geq 1}$ converge to $1$.
\smallskip 

Let $U_0$ be a bi-geodesic region such that $U_0\subsetneq U_1$. Without loss of generality, we may assume that $\gamma(s_{n}) \in \varphi(U_0)$ for all $n \geq 1$. Note that, using Lemma \ref{Lem:Visi_of_images_of_bi-geodesic_reg_inside_U}, $\varphi(U_0)$ and $\OM\setminus \varphi(U_1)$ have visibility property with respect to $(\OM, K_\OM)$. Therefore, by Lemma \ref{Lem:Part_of_Geo_lies_in}, there exists $T>0$ such that $\wh\gamma:= \left.\gamma\right|_{[T, \infty)} \subset \varphi(U_1)$. By reparameterization we may assume that $\wh\gamma:[0,\infty) \to \OM$. Now, using Lemma \ref{Lem:Quasi-geo_wrt_Small_Domain}, we get that $\wh\gamma$ is quasi-geodesic of $(\varphi(\D),K_{\varphi(\D)})$. Therefore, $\varphi^{-1}\circ \wh\gamma$ is quasi-geodesic of $(\D, K_{\D})$. Therefore, $\varphi^{-1}\circ \wh\gamma$ lands at $1$, i.e. $\lim_{t\to\infty}\varphi^{-1}\circ\wh\gamma(t)=1$. Let $\sigma:[0,\infty)\to\D$ be the radial line joining $0$ and $1$. Let $S>0$ such that $\wt\sigma:= \left.\sigma\right|_{[S, \infty)}\subset U_{1}$. By reparameterization, we may assume that $\wt\sigma:[0,\infty) \to \D$. Since $(\D, K_{\D})$ is Gromov hyperbolic, so it follows that there exists a constant $C_{1}>0$ such that for all $t\geq 0$
\[
K_\OM(\wh\gamma(t), \varphi\circ\wt{\sigma}(t))\leq K_{\varphi(\D)}(\wh\gamma(t), \varphi\circ\wt{\sigma}(t))= K_{\D}(\varphi^{-1}\circ \wh\gamma(t), \wt{\sigma}(t))\leq C_{1}.
\]
Then, by Lemma~\ref{Lem:Truncation_of_geo_and_its_asymptoticity_with_radial_line}, there exist a geodesic ray $\eta:[0, \infty) \to \OM$ of $(\OM, K_\OM)$ such that $\eta\subset U_1$ and there exists a constant $C_2>$ such that all $t\geq 0$
\[
K_\OM(\eta(t), \varphi\circ\wt{\sigma}(t))\leq C_{2}.
\]
Therefore, using triangle inequality, we get that for all $t \geq 0$
\[
K_\OM(\wh\gamma(t), \eta(t)) \leq C_{1} + C_{2}=:C_3.
\]
By Lemma~\ref{Lem:Relation_between_origin_geo_horo_and_truncated_geo_horo} and Lemma~\ref{Lem:Asym_compares_two_geo_horo}, there exists a constant $C_4>0$ such that
\[
\clos{H_{\OM}(\gamma, R)}\cap\bdy\OM \subseteq \clos{H_{\OM}(\eta, e^{2T}C_4R)}\cap\bdy\OM.
\]
Note that $\clos{H_{\OM}(\eta, e^{2T}C_4R)}\cap\bdy\OM\neq \emptyset$ as $\xi_1,\xi_2\in \clos{H_{\OM}(\gamma, R)}\cap\bdy\OM$. Thus, by the hypothesis, there exists $x \in \bdy\OM$ such that
\[
\clos{H_{\OM}(\eta, e^{2T}C_4R)}\cap \bdy\OM = \{x\}.
\]
Hence,
\[
\clos{H_{\OM}(\gamma, R)}\cap \bdy\OM \subseteq \clos{H_{\OM}(\eta, e^{2T}C_4R)}\cap \bdy\OM = \{x\}.
\]
Therefore $\xi_{1}=x = \xi_{2}$ which gives a contradiction to $\xi_1\neq \xi_2$. This completes the proof. 
\end{proof}

We now present the proof of Theorem~\ref{Thm:Main_thm_outside_totally_disc_set}.
\begin{proof}[Proof of {\Cref{Thm:Main_thm_outside_totally_disc_set}}]
We have that $(2) \iff (3)$ by Theorem~\ref{Thm:Strengthening_closed_geo_horo_inter_at_single}.\\
$(3) \implies (1)$ Let $f \in \hbox{Hol}(\OM, \OM)$. Suppose that there exists $z\in\OM$ such that $\{f^n(z):n\geq 1\}$ is not relatively compact in $\OM$. If $f^n(z)\to x\in\bdy_{\infty}\OM$ as $n\to\infty$ then, by the fact that the Kobayashi distance is non-expanding property under holomorphic map, we get that, for all $w\in\OM$, $f^n(w)\to x$ as $n\to\infty$. Therefore, in this case, we are done. If there exists two subsequence $(n_k)_{k\geq 1}$ and $(n_l)_{l\geq 1}$ such that $f^{n_k}(z)\to x_1\in\bdy\OM$ and $f^{n_l}(z)\to x_2\in\bdy_{\infty}\OM$ where $x_1\neq x_2$. Pick a sequence $(l_k)_{k\geq1}$ such that $n_{l_k}>n_k$. Let $\gamma_k:[a_k,b_k]\to\OM$ be a geodesic of $(\OM,K_\OM)$ such that $\gamma_k(a_k)=f^{n_{l_k}}(z)$ and $\gamma_k(b_k)=f^{n_k}(z)$. We now make the following claim.
\smallskip

\noindent{\bf Claim~1.} $(\gamma_k)_{k\geq 1}$ eventually avoids every compact subset of $\OM$.
\smallskip

\noindent{\itshape Proof of Claim~1.} The proof is essentially same as the proof of Lemma~\ref{Lem:Whole_Orbit_lies_if_one_of_its_subseq_lies}. Therefore, we omit the proof.
\smallskip
\hfill $\blacktriangleleft$

Using Claim~1 and by Remark~\ref{Rem:Conv_res_of_geo}, there exist $\xi\in\bdy\OM\setminus S$ and a sequence $t_k \to \infty$ such that $(\gamma_k(t_k))_{k \geq 1}$ converges to $\xi$. Let $\varphi :\D\to \OM$ be an injective holomorphic map, $U$ be a bi-geodesic region and $(\gamma_{k_j}(t_{k_j}))_{j\geq 1}$ be a subsequence obtain by $\xi$ and $(\gamma_k(t_k))_{k \geq 1}$ as in Condition~(M$\star$). By Proposition~\ref{Prop:Vis_of_all_req_sets_which_works_everywhere}, there exist bi-geodesic regions $U_1\subsetneq U_2\subsetneq U_3 \subsetneq U$ such that $\varphi(U_1)$ and $\OM\setminus\varphi(U_2)$, $\varphi(U_2)$ and $\OM\setminus\varphi(U_3)$, and $\varphi(U_3)$ and $\OM\setminus\varphi(\D)$ have visibility property with respect to $(\OM, K_\OM)$. By Lemma~\ref{Lem:Seq_of_geo_lies_in_image_of_bi-geo_region}, there exists $K\in \N$ such that $\gamma_{k}\subset \varphi(U_1)$ for all $k\geq K$ and such that both sequences $(\varphi^{-1}(f^{n_{l_k}}(z)))_{k \geq K}$ and $(\varphi^{-1}(f^{n_k}(z)))_{k \geq K}$ converge to $1$. Without loss of generality, we may assume that $\gamma_{k}\subset \varphi(U_1)$ for all $k\geq 1$ and $(\varphi^{-1}(f^{n_{l_k}}(z)))_{k \geq 1}$ and $(\varphi^{-1}(f^{n_k}(z)))_{k \geq 1}$ converge to $1$.
\smallskip

Let $U_0$ be a bi-geodesic region such that $U_0 \subsetneq U_1$. By Lemma~\ref{Lem:Visi_of_images_of_bi-geodesic_reg_inside_U}, $\varphi(U_0)$ and $\OM\setminus \varphi(U_1)$ have visibility property with respect to $(\OM, K_\OM)$. We may assume that $f^{n_k}(z) \in \varphi(U_0)$ for all $k \geq 1$. By Lemma~\ref{Lem:Whole_Orbit_lies_if_one_of_its_subseq_lies}, there exists $N\in \N$ such that $f^{n}(z) \in \varphi(U_1)$ for all $n \geq N$. Without loss of generality, we may assume that $f^{n}(z) \in \varphi(U_1)$ for all $n \geq 1$. Now, by Theorem~\ref{Thm:Global_Loc_Res}, there exists a constant $C_{1}>0$ such that for all $N \geq 1$ and $k \geq 1$
\begin{equation}\label{ineq:Localization_Ineq_in_Main_Result-01}
K_{\varphi(\D)}(f^{n}(z), f^{n_k}(z)) \leq K_\OM(f^{n}(z), f^{n_k}(z)) + C_{1}.
\end{equation}
Let $f(z) = \varphi(\zeta_0)$ for some $\zeta_0 \in \D$. Let $\sigma:[0,\infty) \to \D$ be radial line that joins $0$ and $1$. We now make the following claim.
\smallskip
        
\noindent{\bf Claim~2.} There exists a constant $C_2 > 0$ such that for all $n \geq 1$,
\[
\varphi^{-1}\circ f^{n}(z) \in H_{\D}(\sigma, C_2).
\]
\noindent{\itshape Proof of Claim~2.} The proof is identical to that of the claim in the proof of $(2)\implies (1)$ of Theorem~\ref{Thm:Characterization_of_Condition(*)}, and we therefore omit it.
\smallskip
\hfill $\blacktriangleleft$

Let $S>0$ such that $\wt\sigma:= \left.\sigma\right|_{[S, \infty)}\subset U_{1}$. By reparameterization, we may assume that $\wt\sigma:[0,\infty) \to \D$. By Lemma~\ref{Lem:Whole_Orbit_lies_if_one_of_its_subseq_lies}, there exist a geodesic ray $\eta:[0,\infty)\to\OM$ of $(\OM,K_\OM)$ and a constant $C_3>0$ such that $\eta\subset \varphi(U_1)$ and such that
\[
K_\OM(\varphi\circ\wt\sigma(t),\eta(t))\leq C_3.
\]
By Lemma~\ref{Lem:Relation_between_origin_geo_horo_and_truncated_geo_horo} and  Theorem~\ref{Thm:Geo_horosphere_of_small_domain_is_subset_of_geo_horosphere_of_big_domain}, there exists a constant $C>0$ such that
\begin{equation}\label{ineq:Geo_horosphere_of_small_domain_is_subset_of_geo_horosphere_of_big_domain_02}
H_{\varphi(\D)}(\varphi\circ\sigma, C_2)) \subseteq H_{\OM}(\eta, e^{2S}CC_2).
\end{equation}
It follows from Claim~2 that $ f^{n}(z) \in H_{\varphi(\D)}(\varphi\circ\sigma, C_{2})$ for all $n \geq 1$. Therefore, by~\eqref{ineq:Geo_horosphere_of_small_domain_is_subset_of_geo_horosphere_of_big_domain_02}, we get that for all $n\geq 1$
\begin{equation}\label{Eq:Orbit_lies_in_geo_horo_wrt_unit_disc_01}
f^{n}(z) \in H_{\OM}(\eta, e^{2S}CC_2).
\end{equation}
By hypothesis, either $\clos{H_\OM(\eta,R)}\cap\bdy\OM=\emptyset$ for all $R>0$ or there exists $x\in\bdy\OM$ such that $\clos{H_\OM(\eta,R)}\cap\bdy\OM=\{x\}$ for all $R>0$. Note that $\clos{H_\OM(\eta,R)}\cap\bdy\OM$ cannot be empty because $x_1\in \clos{H_\OM(\eta,R)}\cap\bdy\OM$. 
Therefore, we get
\[
\clos{H_{\varphi(\D)}(\varphi\circ\sigma, C_2))}\cap\bdy\OM\subseteq \clos{H_\OM(\eta,e^{2S}C_2C)}\cap\bdy\OM= \{x\}.
\]
Thus, we have $f^n(z)\to x$ as $n\to\infty$. Since $f^{n_k}(z)\to x_1$ and $f^{n_l}(z)\to x_2$, so it follows that $x_1=x=x_2$ which gives contradictions to the fact that $x_1\neq x_2$. This completes the proof of $(2)\implies (1)$.
\smallskip

$(1) \implies (2)$. Let $R>0$ and let $\gamma:[0,\infty)\to\OM$ be a geodesic of $(\OM,K_\OM)$. If $\clos{H_{\OM}(\gamma, R)}\cap \bdy\OM=\emptyset$, then we are done. Assume that $\clos{H_{\OM}(\gamma, R)}\cap \bdy\OM\neq\emptyset$. If $\clos{H_{\OM}(\gamma, R)}\cap \bdy\OM$ is singleton, then we are done. Assume that $\xi_1,\xi_2\in\clos{H_{\OM}(\gamma, R)}\cap \bdy\OM$ such that $\xi_1\neq\xi_2$. By following the same argument as in $(1)\implies (2)$ in the proof of Theorem~\ref{Thm:Strengthening_closed_geo_horo_inter_at_single}, we obtain the following:
\begin{itemize}
\item[(i)] There exists $\xi\in\bdy\OM\setminus S$ and a sequence $(w_n)_{n\geq 1}$ in $\OM$ such that $w_n\to \xi$. Therefore, we get an injective holomorphic map $\varphi:\D\to\OM$ and a bi-geodesic region $U$ obtained by $\xi$ and $(w_n)_{n\geq 1}$ is in the Condition~(M$\star$). Note that we also obtain bi-geodesic regions $U_1\subsetneq U_2\subsetneq U_3 \subsetneq U$ such that $\varphi(U_1)$ and $\OM\setminus\varphi(U_2)$, $\varphi(U_2)$ and $\OM\setminus\varphi(U_3)$, and $\varphi(U_3)$ and $\OM\setminus\varphi(\D)$ have visibility property with respect to $(\OM, K_\OM)$.
\item[(ii)] There exist $T>0$ and a constant $C_1>0$ such that $\wh\gamma:= \left.\gamma\right|_{[T, \infty)} \subset \varphi(U_1)$ and such that for all $t\geq 0$,
\[
K_\OM(\wh\gamma(t), \varphi\circ\wt{\sigma}(t))\leq C_1.
\]
Where $\sigma$ is radial line joining $0$ and $1$ and $S>0$ such that $\wt\sigma:=\left.\sigma\right|_{[S,\infty)}\subset U_1$. Note that, by reparameterization, we assume that $\wt\sigma:[0,\infty)\to\D$ and $\wh\gamma:[0,\infty) \to \OM$.
\end{itemize}

Let $\Gamma$ be a component of $\bdy\OM$ such that $\xi \in \Gamma$. Let $\OM_{\Gamma}$ be simply connected domain obtained by filling all the components of $\C\setminus\OM$ except which contains $\Gamma$.

The remainder of the proof follows the same argument as in the proof of $(1)\implies (2)$ of Theorem~\ref{Thm:Characterization_of_Condition(*)}. At the end, we get that there exists $\xi_0\bdy\OM$ such that 
\[
\clos{H(\gamma,R)}\cap\bdy\OM=\{\xi_0\}.
\]
Thus, $\xi_1=\xi_0=\xi_2$. This gives a contradiction to $\xi_1\neq\xi_2$. This completes the proof of $(1)\implies (3)$. This completes the proof of theorem.
\end{proof}

\section{Examples}\label{Sec:Examples}
The first example is an example of planar simply connected domain where Result~\ref{Res:Global_localization_result} is not applicable, however, we can localize the Kobayashi distance using Theorem~\ref{Thm:Global_to_loc_vis}. That is, intrinsic (weak) visibility is more general than the standard (weak) visibility.
\begin{example}\label{Ex:Our_loc_holds_while_previous_one_doesnot}
Let 
\[
\OM:= \{z\in \C: 0<\rprt(z) < 1, 0< \iprt(z)<1\}\setminus \bigcup_{j=2}^{\infty}\{z\in\C:\rprt(z)= 1/j, 0\leq \iprt(z)\leq 1/2\}.
\]
Clearly, $\OM$ is a bounded simply connected domain and $\partial\OM$ is not locally connected. Hence, it follows from \cite[Corollary~3.4]{BNT2022} that $\OM$ is not a visibility domain. Let $a \in \R$ with $a>1/2$. Let 
\[
U:=\{z\in \OM: 0<\rprt(z)<a, 0<\iprt(z)<a\}.
\]
Let $L:=\{\zeta\in \C: \rprt\zeta = 0, 0 \leq \iprt \zeta \leq 1/2\}$. We now make the following claim.
\smallskip
    
\noindent{\bf Claim~1.} Any pair of two distinct points $x, y \in L$ does not satisfy the visibility property with respect to $(\OM,K_\OM)$.
\smallskip
    
\noindent{\itshape Proof of Claim~1.} Without loss of generality, we may assume that $\iprt x< \iprt y$. Consider a sequence $(z_{n})_{n\geq 1}\in \OM$ such that $\iprt z_{n} = x$ and $\rprt z_{n} \in (1/(n+1), 1/n)$ for all $n \geq 1$. Thus, $(z_{n})_{n\geq 1}$ converges to $x$. Let $L_{y}:=\{z\in\C: 0\leq \rprt z\leq 1, \iprt z = \iprt y\}$. Fix a point $p \in \OM$ with $\rprt p =1/2$ and $\iprt p \in (1/2,1)$. Let $\gamma_{n}:[a_{n},b_{n}]\to \OM$ be a geodesic of $(\OM, K_\OM)$ such that $\gamma(a_{n}) = z_{n}$ and $\gamma_{n}(b_{n})=p$. Fix
\[
c_{n}:= \inf\{t\in[a_n, b_n]: \gamma_{n}(t) \in L_y\}.
\]
By construction, we have that, for all $n \in \N$, $\wt\gamma_n:=\left.\gamma_{n}\right|_{[a_{n}, c_{n}]}\subseteq (1/(n+1), 1/n)\times [0,\iprt y]$. Let $w_n:=\gamma_{n}(c_n)$. It is clear that $\rprt w_n \to 0$ as $n \to \infty$. That is, $w_n \to y$ as $n \to \infty$. Let $K\subset \OM$ be an arbitrary compact set. Then there exists $N \in \N$ such that $\wt{\gamma}_{n}\cap K = \emptyset$ for all $n \geq N$. Therefore, $x$ and $y$ does not have visibility property with respect to $(\OM,K_\OM)$. 
\hfill $\blacktriangleleft$
\smallskip
\begin{figure}[htbp]
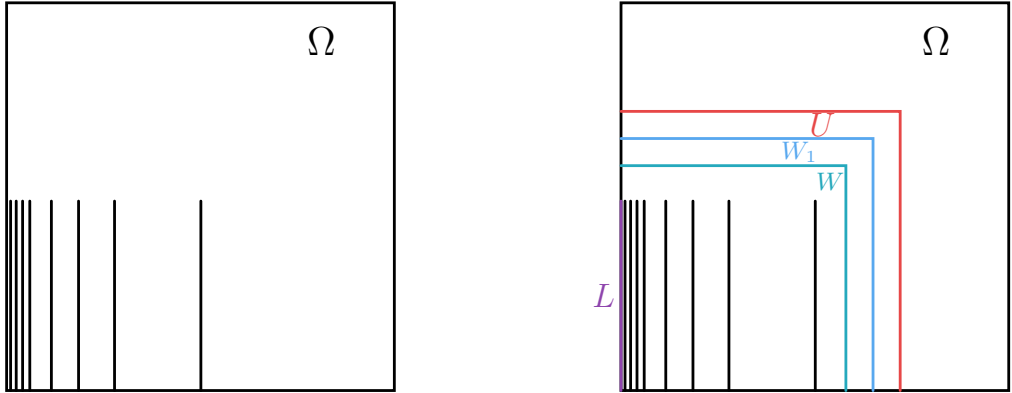

  \centering
  \begin{minipage}[t]{0.50\textwidth}
    \centering
    \resizebox{0.70\linewidth}{!}{\CombDomainPanel{1}{0}}
  \end{minipage}\hfill
  \begin{minipage}[t]{0.50\textwidth}
    \centering
    \resizebox{0.70\linewidth}{!}{\CombDomainPanel{4}{1}}
  \end{minipage}
  \caption{An illustration of domain $\OM$}
  \label{fig:comb-domain-pair}
\end{figure}

It follows from Claim~1 that there exists at least one pair of distinct points in $\bdy U\cap\bdy\OM$ which does not have visibility property with respect to $(\OM,K_\OM)$. Therefore, we cannot localize the Kobayashi distance using Result~\ref{Res:Global_localization_result}.
\smallskip

Let $b,c \in \R$ such that $1/2<c<b<a$. Consider the following sets
\begin{align*}
W_1\defeq&\{z\in \OM: 0<\rprt(z)< b, 0<\iprt(z)<b\} \text{ and }\\
W\defeq&\{z\in \OM: 0<\rprt(z)< c, 0<\iprt(z)<c\}.
\end{align*}
We now make the following claim:
\smallskip

\noindent{\bf Claim~2.} $W$ and $\OM\setminus W_1$ have visibility property with respect to $(\OM,K_\OM)$.
\smallskip

\noindent{\itshape Proof of Claim~2.} Suppose that $W$ and $\OM\setminus W_1$ does not have visibility property with respect to $(\OM,K_\OM)$. There exists a sequence $\gamma_n:[0,a_n]\to\OM$ of geodesic of $(\OM,K_\OM)$ with $\gamma_n(0)\in W$ and $\gamma_n(a_n)\in \OM\setminus W_1$ such that $(\gamma_n)_{n\geq 1}$ eventually avoids compact subset of $\OM$. Up to a subsequence, $(\gamma_n(a_n))$ converges to some $x\in \clos{W_1}\cap\bdy\OM$. Let $r>0$ such that $\clos{B(x,r)}\cap \clos{W}=\emptyset$. Consider
\[
t_n:=\sup\{t\in[0,a_n]:\gamma_n(s)\in B(x,r)\,\forall s\in [t,a_n]\}.
\]
Let $\wh\gamma_n:=\left.\gamma_n\right|_{[t_n,a_n]}$. We may assume $\gamma_n(t_n) \to y\in\bdy B(x,r)\bdy\OM$ as $n \to \infty$. Therefore, $(\wh\gamma_n)_{n\geq 1}$ eventually avoids compact subset of $\OM$ which contradicts visibility property of the pair of points $x$ and $y$ with respect to $(\OM,K_\OM)$. Hence $W$ and $\OM\setminus W_1$ have visibility property with respect to $(\OM, K_\OM)$. This completes the proof of the Claim.
\hfill $\blacktriangleleft$ 
\smallskip

Similarly, we may prove that $W_1$ and $\OM\setminus U$ have visibility property with respect to $(\OM, K_\OM)$. Therefore, we may localize the Kobayashi distance using Theorem~\ref{Thm:Global_to_loc_vis}. 
\end{example}

Let $(X,d)$ be a metric space, and let $A,B\subset X$. We denote by $d^H(A,B)$ the Hausdorff distance between $A$ and $B$ with respect to the metric $d$. For further details on the Hausdorff distance, we refer the reader to \cite{PapaBook}. The proof of the following lemma is directly adapted from the Proposition~10.1.4 in \cite{PapaBook}.
\begin{lemma}\label{Lem:Aysm_&_Haus_dis_finite_r_equi}
Let $(X, d)$ be a metric space and let $\kappa\geq 0$. Let $\gamma_1:[0,\infty) \to X$ be a $(1,\kappa)$-quasi-geodesic ray of $(X,d)$ and $\gamma_2:[0,\infty)\to X$ be a geodesic of $(X,d)$. Then $d_{d}^{H}(\gamma_1,\gamma_2)<\infty$ if and only if there exists a constant $C>0$ such that for all $t \geq 0$
\[
d(\gamma_1(t),\gamma_2(t)) \leq C.
\]
\end{lemma}
\begin{proof}
Suppose that there exists a constant $C>0$ such that $d(\gamma_1(t),\gamma_2(t)) \leq C$ for all $t\geq 0$. Fix $t \geq 0$. Then by definition, $d(\gamma_1(t),\gamma_2)\leq d(\gamma_1(t), \gamma_{2}(t))\leq C$. Similarly $d(\gamma_1,\gamma_2(t))\leq C$. Hence $\sup_{t\geq 0}d(\gamma_1(t),\gamma_2)\leq C$ and $\sup_{t\geq 0}d(\gamma_1,\gamma_2(t))\leq C$. Therefore, $d^{H}(\gamma_1,\gamma_2)\leq C$.
\smallskip

Conversely, assume that $C_{0}:= d^{H}(\gamma_1, \gamma_2)< \infty$. Fix $t\in[0,\infty)$. Let $t^\prime\geq 0$ such that $d(\gamma_1(t),\gamma_2(t^\prime))=d(\gamma_1(t),\gamma_2)\defeq\inf\{d(\gamma_1(t),\gamma_2(s)):s\in[0,\infty)\}$. Therefore, $d(\gamma_1(t),\gamma_2(t^\prime))\leq C_0$. By using $\gamma_2$ is a geodesic and triangle inequality,
\begin{align*}
t^\prime &=d(\gamma_2(t^\prime), \gamma_2(0))\\
&\leq d(\gamma_2(t^\prime),\gamma_1(t))+d(\gamma_1(t),\gamma_1(0))+d(\gamma_1(0),\gamma_2(0))\\
& \leq C_0 + t + d(\gamma_1(0),\gamma_2(0))\\
& \leq C_0 + t +\k + d(\gamma_1(0),\gamma_2(0)).
\end{align*}
Therefore, $t^\prime-t\leq C_0+ \k+d(\gamma_1(0),\gamma_2(0))$. On the other hand, using fact that $\gamma_1$ is $(1,\k)$-quasi-geodesic of $(X,d)$ and triangle inequality, we get
\begin{align*}
t&\leq d(\gamma_1(t),\gamma_1(0))+\k\\
&\leq  d(\gamma_1(t),\gamma_2(t^\prime))+d(\gamma_2(t^\prime),\gamma_2(0))+d(\gamma_2(0),\gamma_1(0))+\k\\
&\leq C_0 + t^\prime + d(\gamma_1(0),\gamma_2(0))+\k.
\end{align*}
Therefore, $-C_0 - t^\prime - d(\gamma_1(0),\gamma_2(0))-\k \leq t^\prime-t$. Hence, we get
\begin{equation}\label{Ineq:Bound_for_modulus}
|t^\prime-t|\leq C_0+ \k+d(\gamma_1(0),\gamma_2(0)).    
\end{equation}
Therefore, using triangle inequality, the fact that $\gamma_2$ is a geodesic of $(X,d)$ and \ref{Ineq:Bound_for_modulus}, for all $t\geq 0$
\begin{align*}
d(\gamma_1(t),\gamma_2(t))&\leq d(\gamma_1(t),\gamma_2(t^\prime))+d(\gamma_2(t^\prime),\gamma_2(t))\\
&\leq C_0 + |t-t^\prime|\\
&\leq 2C_0 + \k+d(\gamma_1(0),\gamma_2(0)).
\end{align*}
Therefore, by taking $C:=2C_0+ \k+d(\gamma_1(0),\gamma_2(0))$, we get for all $t\geq 0$
\[
d(\gamma_1(t),\gamma_2(t))\leq C.
\]
This completes the proof.
\end{proof}

We shall now see an example of infinitely connected domain satisfying the Condition~(M$\star$) such that it does not have the Denjoy–Wolff property. The main idea behind the construction is to consider a domain that has a bounded boundary component such that, after a suitable modification, it becomes a bounded subdomain of the ambient domain and so that subdomain fails to satisfy the Denjoy–Wolff property. We construct a geodesic ray of the ambient domain such that the closure of the associated Busemann horosphere meets the boundary of the domain more than one points. Therefore, using Theorem~\ref{Thm:Characterization_of_Condition(*)}, we get that the ambient domain does not have the Denjoy-Wolff property. The following proposition may be proved similar to Proposition~4.4 in \cite{Benini-Bracci} which essentially provides such subdomain.
\begin{proposition}\label{Prop:Geo_horo_is_nsingleton_for_expo_comb_dom}
Let $D\subset\C$ be domain which is an interior of $U:=\Gamma\cup L$ where $L:=\{z\in \C:|\rprt z|\leq 1, \iprt z= 1\}$; and
\[
\Gamma:= \bigcup_{k =1 }^{\infty}\Big\{z\in \C: x_k\leq |\rprt z| \leq 1, \iprt z= y_k\Big\}\cup\Big(\bigcup_{i=1}^3L_i\Big)
\]
where $x_k = \frac{1}{2^ke^{3k}}$ , $y_k = \frac{1}{2^k}$, and
\begin{align*}
L_{1} &= \{z\in\C:\rprt z= -1, 0\leq \iprt z\leq 1\},\\
L_2 &= \{z\in\C: |\rprt z|\leq 1, \iprt z = 0\} \text{ and }\\
L_{3} &= \{z\in\C:\rprt z= 1, 0\leq \iprt z\leq 1\}.
\end{align*}
Let $0<p<1$ and let $\wt\gamma:[0,\infty)\to D$ be a geodesic of $(D, K_D)$ given by $\wt\gamma(t) := \frac{ip}{1+t}$. Then for all $R>0$
\[
\clos{H_{\OM}(\wt\gamma,R)}\cap\bdy\OM= [-1,1].
\]
\end{proposition}

\begin{example}\label{Ex:Dom_doesnot_have_the_DWP}
Consider a domain
\[
\OM:= \C\setminus \bigcup_{(m,n)\in \Z^2}\Gamma_{m,n}
\]
where $\Gamma_{m,n}:= (m+in)+r\Gamma$ and $r>0$ is chosen such that there exists a $\delta>0$ such that $\distance(\Gamma_{m_1,n_1}, \Gamma_{m_2,n_2})>\delta$ for any $m_1,n_1,m_2,n_2\in \Z$; and
\[
\Gamma:= \bigcup_{k =1 }^{\infty}\Big\{z\in \C: x_k\leq |\rprt z| \leq 1, \iprt z= y_k\Big\}\bigcup_{i=1}^{3}L_{i}
\]
where $x_k := \frac{1}{2^ke^{3k}}$ , $y_k := \frac{1}{2^k}$ and
\begin{align*}
L_{1} &:= \{z\in\C:\rprt z= -1, 0\leq \iprt z\leq 1\},\\
L_2 &:= \{z\in\C: |\rprt z|\leq 1, \iprt z = 0\} \text{ and }\\
L_{3} &:= \{z\in\C:\rprt z= 1, 0\leq \iprt z\leq 1\}.
\end{align*}
We shall prove that $\OM$ does not have Denjoy-Wolff property. It is clear that $\OM$ satisfies Condition~(M$\star$) as it satisfies Condition $(\star\star)$. For simplicity, we may assume $r=1$ otherwise we may work with the domain $D^\prime:=rD$. By Theorem~\ref{Thm:Characterization_of_Condition(*)}, it suffices to prove the following.
\begin{proposition}\label{Prop:Lattice_of_expo_comb_domain_is_not_DW_dom}
Let $\OM$ be a domain as in Example~\ref{Ex:Dom_doesnot_have_the_DWP}. Then there exists a geodesic ray $\gamma:[0,\infty) \to \OM$ of $(\OM,K_\OM)$ such that for all $R>0$
\[
[-1,1]\subseteq \clos{H_{\OM}(\gamma,R)}\cap\bdy\OM.
\]
In particular, $\OM$ does not have Denjoy-Wolff property.
\end{proposition}
\begin{proof}
Let $D$ as in Proposition~\ref{Prop:Geo_horo_is_nsingleton_for_expo_comb_dom}. Let $c, d\in\R$ such that $0<c<d<1$. Consider 
\[
W:=\{z\in D:|\rprt z|<1, 0<\iprt z<c\} \text{ and }W_1:=\{z\in D:|\rprt z|<1, 0<\iprt z<d\}
\]
such that $p \in W$ where $p$ is same as in Proposition~\ref{Prop:Geo_horo_is_nsingleton_for_expo_comb_dom}. Let $\wt\gamma:[0,\infty)\to D$ be a geodesic of $(D,K_D)$ as in Proposition~\ref{Prop:Geo_horo_is_nsingleton_for_expo_comb_dom}. Without loss of generality, we may assume that $\wt\gamma\subset W$. We now make the following claim:
\smallskip

\noindent{\bf Claim~1.} $W$ and $\OM\setminus W_1$, and $W_1$ and $\OM\setminus D$ have visibility property with respect to $(\OM, K_\OM)$.
\smallskip

\noindent{\itshape Proof of Claim~1.} Assume that $W$ and $\OM\setminus W_1$ does not have visibility property with respect to $(\OM, K_\OM)$. There exists a geodesic $\eta_n:[0,T_n]\to\OM$ with $\eta_n(0)\in W$ and $\eta_n(T_n)\in \OM\setminus W_1$ such that $(\eta_n)_{n \geq 1}$ eventually avoids every compact subset of $\OM$. Define
\begin{align*}
a_n^\prime:&=\sup\{t \in [0,T_n]: \eta_n(s)\in W \, \forall s\in[0,t]\} \text{ and }\\
b_n^\prime:&=\inf\{t\geq s_n: \eta_n(t)\not\in W_1\}
\end{align*}
By construction, $\wh\eta_n:=\left.\eta_n\right|_{[a^\prime_n,b^\prime_n]} \subset \clos{W_1\setminus W}$. Therefore, $(\wh\eta_n)_{n \geq 1}$ also eventually avoids every compact subset of $\OM$. We may assume, up to a subsequence, $(\wh\eta_n(a^\prime_n))_{n \geq 1}$ converge to $x^\prime\in \bdy D$ and $(\wh\eta_n(b^\prime_n))_{n \geq 1}$ converge to $y^\prime \in \bdy D$. It is easy to see that $x^\prime \neq y^\prime$. We may have following two cases:
\smallskip

\noindent{\bf Case~1.} Either $x^\prime, y^\prime \in \{z\in C: \rprt z= -1, \iprt z=c \text{ or }d \}$ or $x^\prime, y^\prime \in \{z\in C: \rprt z= 1, \iprt z=c \text{ or }d\}$.

In this case, it is clear that $x^\prime$ and $y^\prime$ have visibility property with respect to $(\OM,K_\OM)$. Therefore, there exist a compact set $K\subset \OM$ and $N\in \N$ such that $\wh\eta_n\cap K\neq \emptyset$ for all $n\geq N$. This gives a contradiction to our starting assumption. 
\smallskip

\noindent{\bf Case~2.} Either $\rprt x^\prime =-1$ and $\rprt y^\prime = 1$ or $\rprt x^\prime =1$ and $\rprt y^\prime = -1$.

In this case, $(\wh\eta_n)_{n \geq 1}$ must intersect $\{z\in \C:\rprt z=0, \iprt z \in [c,d]\}$ which is a contradicts. 

Therefore, $W$ and $\OM\setminus W_1$ have visibility property with respect to $(\OM, K_\OM)$. Similarly, we may prove that $W_1$ and $\OM\setminus D$ have visibility property with respect to $(\OM, K_\OM)$. This completes the proof of claim.
\hfill $\blacktriangleleft$
\smallskip

Let $(z_n)_{n\geq 1}\subset \wt\gamma$ be a sequence that converges to $0$ and let $q\in W$ such that $p\neq q$. Let $\gamma_n:[0,c_n]\to \OM$ be a geodesic of $(\OM,K_\OM)$ such that $\gamma_n(0)=q$ and $\gamma_n(c_n)=z_n$. Since $(\OM, K_\OM)$ is complete, $(\gamma_n)_{n \geq 1}$, upto a subsequence, converges uniformly on compact subsets of $[0,\infty)$ to a geodesic $\gamma:[0,\infty) \to \OM$ of $(\OM,K_\OM)$. Define
\[
T:=\inf\{t\in[0,\infty):\gamma(s)\in W \,\forall s\in [t,\infty)\}.
\]
Without loss of generality, we may assume that $\wh\gamma:=\left.\gamma\right|_{[T,\infty)}\subset W$. We now make the following claim:
\smallskip

\noindent{\bf Claim~2.} $K^{H}_D(\wt\gamma,\wh\gamma ) <+\infty$.
\smallskip

\noindent{\itshape Proof of Claim~2.} Let $\wt\gamma_n:[0,d_n]\to D$ be a geodesic of $(D,K_D)$ such that $\wt\gamma_n(0)=p$ and $\wt\gamma_n(d_n)=z_n$. Therefore, $(\wt\gamma_n)_{n\geq 1}$ converges, up to a subsequence, to a geodesic ray $\wt\gamma$. Let $\eta:[a,b]\to D$ be a geodesic of $(D,K_D)$ such that $\eta(a)=p$ and $\eta(b)=q$. Since $(D,K_D)$ is Gromov hyperbolic, there exists a constant $C_0>0$ such that for all $ n \in \N$
\[
K^H_{D}(\wt\gamma_n\cup\eta,\wh\gamma_n)\leq C_0.
\]
Therefore, by taking $n \to \infty$, we get
\[
K^H_{D}(\wt\gamma\cup\eta,\wh\gamma)\leq C_0.
\]
Since $\eta$ is compact set, there exists a constant $C_1>0$ such that
\[
K^H_{D}(\wt\gamma,\wh\gamma)\leq C_1.
\]
Hence, $K^{H}_D(\wt\gamma,\wh\gamma ) <+\infty$ which completes the proof of claim.
\hfill $\blacktriangleleft$
\smallskip

Since $\wt\gamma \subset W$, by Theorem~\ref{Thm:Global_Loc_Res}, there exists a constant $C_2>0$ such that for all $s, t\in[0,\infty)$,
\[
K_D(\wt\gamma(s),\wt\gamma(t))\leq K_\OM(\gamma_1(s),\gamma_2(t))+ C_2.
\]
This gives that $\wt\gamma$ is a $(1,\k)$-quasi-geodesic of $(\OM,K_\OM)$ where $\k:=C_2$. Therefore, using Claim~2 in Lemma~\ref{Lem:Aysm_&_Haus_dis_finite_r_equi}, there exists a constant $C_3>0$ such that for all $t\geq 0$
\begin{equation}\label{Ineq:Asym_of_two_geo-0001}
K_\OM(\wt\gamma(t), \wh\gamma(t)) \leq C_3.
\end{equation}
Therefore, using~\ref{Ineq:Asym_of_two_geo-0001} in Theorem~\ref{Thm:Geo_horosphere_of_small_domain_is_subset_of_geo_horosphere_of_big_domain}, there exists a constant $C>0$ such that for all $R>0$
\[
\clos{H_D(\wt\gamma, e^{2T}R)}\cap\bdy\OM\subseteq \clos{H_\OM(\wh\gamma, e^{2T}CR)}\cap\bdy\OM.
\]
Therefore, using Lemma~\ref{Lem:Relation_between_origin_geo_horo_and_truncated_geo_horo} and Proposition~\ref{Prop:Geo_horo_is_nsingleton_for_expo_comb_dom}, we get
\[
[-1,1]\subseteq \clos{H_\OM(\gamma, CR)}\cap\bdy\OM.
\]
Hence, $\OM$ does not have Denjoy-Wolff property. This completes the proof of proposition.
\end{proof}
\end{example}

The third example is an example of infinitely connected bounded domain that satisfies the Denjoy-Wolff property. The main idea behind the construction is to remove infinitely many copies of the compact set from the unit disc such that
\begin{itemize}
\item[(i)] the diameter of compact sets tends to zero as we approach to the boundary of the unit disc, and 
\item[(ii)] after a suitable modification, it become a subdomain of the unit disc and so that subdomain satisfies the Denjoy–Wolff property.
\end{itemize}
We record the following result about the required compact set and associated subdomain.
\begin{result}[Proposition~4.1 in \cite{Benini-Bracci}]\label{Res:Double_sided_comb_dom_is_DW_dom}
Let $r>0$ and let $D\subset \C$ be a domain which is interior of $U:= r(\Sigma\cup G)$ where $G:=\{z\in \C:|\rprt z|\leq 1, \iprt z= 1\}$, and
\[
\Sigma:=\bigcup_{j =2}^\infty\Big\{z\in\C:\rprt z  = \pm\frac{1}{j},-1 \leq \iprt z\leq 0\Big\}\cup\Big(\bigcup_{i=1}^4L_i\Big)
\]
where
\begin{align*}
L_{1} &\defeq \{z\in\C:\rprt z= -1, -1\leq \iprt z\leq 1\},\\
L_2 &\defeq \{z\in\C: |\rprt z|\leq 1, \iprt z = -1\}\\
L_{3} &\defeq \{z\in\C:\rprt z= 1, -1\leq \iprt z\leq 1\} \text{ and }\\
L_{4}&\defeq \{z\in\C:\rprt z= 0, -1\leq \iprt z\leq 0\}
\end{align*}
Then for every geodesic ray $\wt\gamma:[0,\infty) \to D$ of $(D,K_D)$, there exists $x\in \bdy D$ such that for all $R>0$
\[
\clos{H_D(\wt\gamma, R)}\cap\bdy D=\{x\}.
\]
\end{result}

\begin{example}\label{Ex:The_inf_con_domain_with_the_DWP}
Let $\nu \in \N$ be a fixed natural number. Let $S:=\{e^{\frac{2\pi ik}{\nu}}: k = 1,\ldots,\nu\}\subset\bdy\D$. Let $L_k$ be the radial line that joins $0$ and point $\xi_k:=e^{\frac{2\pi ik}{\nu}}$. Let $(\eta^k_n)_{n \geq 1}\subset L_k$ be a sequence with $\eta_1^k = 0$ such that $(\eta^k_n)_{n \geq 1}$ converges to $\xi_k$ and such that $|\eta^k_n|< |\eta^k_{n+1}|$ for all $n \geq 1$. Moreover, we assume that if $(\eta^{k_{1}}_n)_{n \geq 1}\subset L_{k_1}$ and $(\eta^{k_{2}}_n)_{n \geq 1} \subset L_{k_2}$ are two sequences then $|\eta^{k_{1}}_n|= |\eta^{k_{1}}_n|$ for all $n \geq 1$. We can choose $(r_n)_{n \geq 1}$ with $r_n>0$ and $r_n \to 0$ as $ \to \infty$ such that if $k \in \{1,\ldots, \nu\}$, then $\clos{B_n(\eta^k_n,r_n)}\cap\clos{B_{n+1}(\eta^k_{n+1},r_{n+1})} = \emptyset$ for all $n \in \N$.
Consider a domain defined as:
\[
\OM:=\D\setminus\bigcup_{k =1}^{\nu}\Sigma_k
\]
where 
\[
\Sigma_k:= \bigcup_{n =1}^\infty(\eta^k_n + r_n\Sigma).
\]
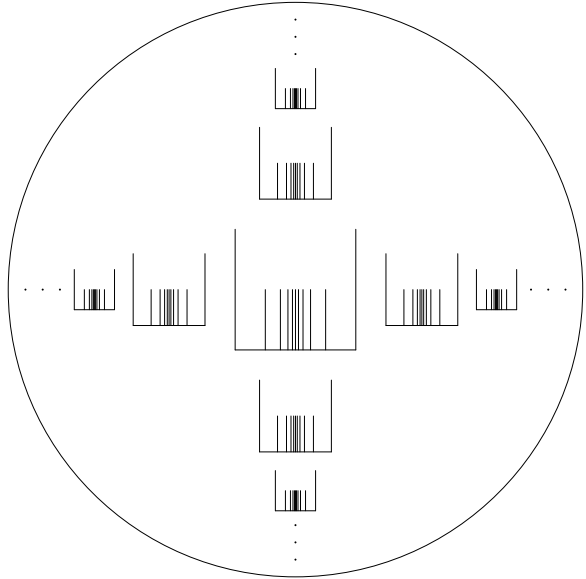
\begin{figure}[htbp]
  \centering
  \begin{tikzpicture}[
    scale=3.8,
    line cap=round,
    line join=round
  ]
    \draw[line width=0.10pt] (0,0) circle[radius=1cm];

    \foreach \cx/\cy/\kscale in {
      0/0/0.21,
      0.44/0/0.125, 0.70/0/0.07,
      -0.44/0/0.125, -0.70/0/0.07,
      0/0.44/0.125, 0/0.70/0.07,
      0/-0.44/0.125, 0/-0.70/0.07
    } {
      \pgfmathsetmacro{\CombOutlineWidth}{0.35*\kscale}
      \pgfmathsetmacro{\CombToothWidth}{0.3*\kscale}
      \pgfmathsetmacro{\CombSpineWidth}{0.45*\kscale}
      \begin{scope}[shift={(\cx,\cy)},scale=\kscale]
        \draw[line width=\CombOutlineWidth pt] (-1,-1) -- (-1,1);
        \draw[line width=\CombOutlineWidth pt] (1,-1) -- (1,1);
        \draw[line width=\CombOutlineWidth pt] (-1,-1) -- (1,-1);

        \foreach \n in {2,4,8,20} {
          \pgfmathsetmacro{\CombTooth}{1/\n}
          \draw[line width=\CombToothWidth pt]
            (\CombTooth,-1) -- (\CombTooth,0)
            (-\CombTooth,-1) -- (-\CombTooth,0);
        }
        \draw[line width=\CombSpineWidth pt] (0,-1) -- (0,0);
      \end{scope}
    }

    \foreach \dotx/\doty in {
      0.82/0, 0.88/0, 0.94/0,
      -0.82/0, -0.88/0, -0.94/0,
      0/0.82, 0/0.88, 0/0.94,
      0/-0.82, 0/-0.88, 0/-0.94
    } {
      \fill (\dotx,\doty) circle[radius=0.10pt];
    }
  \end{tikzpicture}
  \caption{An illustration of domain $\OM$ when $\nu =4$.}
  \label{fig:shrinking-copies-axis}
\end{figure}

We want to show that $\OM$ have Denjoy-Wolff property. Note $\OM$ satisfies Condition~(M$\star$) at every point $\xi\in\bdy\D\setminus S$. Therefore, using Theorem~\ref{Thm:Main_thm_outside_totally_disc_set}, it suffices to prove the following:
\begin{proposition}
Let $\OM\subset\C$ be a domain as in Example~\ref{Ex:The_inf_con_domain_with_the_DWP}. Then for every geodesic ray $\gamma:[0,\infty) \to \OM$ of $(\OM,K_\OM)$, there exists $x \in \bdy\OM$ such that $\clos{H_\OM(\gamma,R)}\cap\bdy\OM=\{x\}$ for all $R>0$.  
\end{proposition}
\begin{proof}
Given $k \in \{1,\ldots,\nu\}$ and $n \in \N$, consider a line 
\[
G^k_n:=\{z\in \D:-r_n \leq \rprt (z-\eta^k_n)\leq r_n , \iprt (z-\eta^k_n)=r_n \}.
\]
Let $D^k_n$ be a subdomain of $\D$ which is an interior of $(\eta^k_n + r_n\Sigma)\cup G^k_n$. Let $0<c^k_n<d^k_n<r_n$. Consider
\begin{align*}
W^k_n:&=\{z\in D^k_n: -r_n \leq \rprt (z-\eta^k_n)\leq r_n, -r_n<\iprt(z-\eta^k_n)<c^k_n\} \text{ and }\\ 
\wt W^k_n:&= \{z\in D^k_n:-r_n \leq \rprt (z-\eta^k_n)\leq r_n,-r_n<\iprt(z-\eta^k_n)<d^k_n\}.    
\end{align*}
By following the same argument of the proof of Claim~1 in Proposition~\ref{Prop:Lattice_of_expo_comb_domain_is_not_DW_dom}, we may show that $W^k_n$ and $\OM\setminus \wt W^k_n$, and $\wt W^k_n$ and $\OM\setminus D^k_n$ have visibility property with respect to $(\OM,K_\OM)$. Let $\gamma:[0,\infty) \to \OM$ be a geodesic of $(\OM,K_\OM)$. We may have following two cases:
\smallskip

\noindent{\bf Case~1.} There exists $T>0$ such that $\wh\gamma:=\left.\gamma\right|_{[T,\infty)} \subset D^k_n$ for some $k \in \{1,\ldots, \nu\}$ and $n \in \N$.
\smallskip

In this case, we shall use Result~\ref{Res:Double_sided_comb_dom_is_DW_dom} for domain $D^k_n$. Without loss of generality, we may assume that $\wh\gamma\subset W^k_n$. By reparameterization, we may assume that $\wh\gamma:[0,\infty) \to \OM$. Let $p\in W^k_n$ with $p\neq\wh\gamma(0)$ and let $(t_n)_{n\geq 1}\subset[0,\infty)$ be a sequence such that $t_n\to \infty$ as $n \to \infty$. Let $\wt\gamma_n:[0,a_n]\to D^k_n$ be a geodesic of $(D^k_n,K_{D^k_n})$ such that $\wt\gamma_n(0)=p$ and $\wt\gamma_n(a_n)=\wh\gamma(t_n)$. Then, up to a subsequence, $(\wt\gamma_n)_{n \geq 1}$ converges uniformly on compact subsets of $[0,\infty)$ to a geodesic $\wt\gamma:[0,\infty) \to \OM$ of $(D^k_n,K_{D^k_n})$. Then there exists $M>0$ such that $\left.\wt\gamma\right|_{[M,\infty)}\subset W^k_n$. By the abuse of notation, we denote $\left.\wt\gamma\right|_{[M,\infty)}$ by $\wt\gamma$. By reparameterization, we may assume that $\wt\gamma:[0,\infty) \to D^k_n$. By Theorem~\ref{Thm:Global_Loc_Res}, there exists a constant $C_0>0$ such that for all $s,t\in [0,\infty)$
\[
K_{D^k_n,}(\wt\gamma(s),\wt\gamma(t))\leq K_\OM(\wt\gamma(s),\wt\gamma(t))+C_0.
\]
Therefore, $\wt\gamma$ is $(1,\k)$-almost-geodesic of $(\OM,K_\OM)$ where $\kappa:=C_0$. Now, by following the same argument of the proof of Claim~2 of Proposition~\ref{Prop:Lattice_of_expo_comb_domain_is_not_DW_dom}, we get
\[
K^{H}_D(\wt\gamma,\wh\gamma ) <+\infty.
\]
By Lemma~\ref{Lem:Aysm_&_Haus_dis_finite_r_equi}, there exists a constant $C^\prime>0$ such that for all $t \geq 1$ 
\[
K_\OM(\wt\gamma(t),\wh\gamma(t)) <C^\prime.
\]
By Theorem~\ref{Thm:Geo_horosphere_of_big_domain_is_subset_of_geod_of_small_domain_upto_some_compact_set_and-constant} and Lemma~\ref{Lem:Relation_between_origin_geo_horo_and_truncated_geo_horo}, there exists a constant $C_2>0$ such that
\[
H_\OM(\gamma,R)\subseteq H_{D^k_n}(\wt\gamma,e^{2T}C_2R)
\]
Now, by Result~\ref{Res:Double_sided_comb_dom_is_DW_dom}, there exists $x\in\bdy\OM$ such that 
\[
\clos{H_{D^k_n}(\wt\gamma,e^{2T}C_2R)}\cap\bdy\OM=\{x\}.
\]
Therefore, we get that $\clos{H_\OM(\gamma,R)}\cap\bdy\OM=\{x\}$.
\smallskip

\noindent{\bf Case~2.} There do not exist $T>0$, $k \in \{1,\ldots, \nu\}$ and $n \in \N$ such that $\left.\gamma\right|_{[T,\infty)} \subset D^k_n$.
\smallskip

In this case, we choose a boundary point $\xi\in\bdy\OM, \epsilon>0$ and $T^\prime>0$ such that $\wh\gamma\defeq \left.\gamma\right|_{[T^\prime,\infty)}\subset B(\xi,\epsilon)\cap\OM$. It is clear that $B(\xi,\epsilon)\cap\OM$ is a bounded simply connected domain such that its boundary is Jordan curve. Hence, $U\defeq B(\xi,\epsilon)\cap\OM$ has the Denjoy-Wolff property. Therefore, using Result~\ref{Thm:Main_Thm_of_BB}, for every geodesic ray $\wt\gamma:[0,\infty)\to U$ of $(U,K_U)$, there exists $x^\prime\in\bdy U$ such that for all $R>0$
\begin{equation}\label{Eq:Small_dom_geo_horo_int_at_sing}
\clos{H_U(\wt\gamma,R)}\cap\bdy U =\{x^\prime\}.    
\end{equation}
Let $W\defeq B(\xi,\epsilon/3)\cap\OM$ and $W_1\defeq B(\xi,\epsilon/2)\cap\OM$. It is easy to see that $W$ and $\OM\setminus W_1$, and $W_1$ and $\OM\setminus U$ have visibility property with respect to $(\OM,K_\OM)$. Now, by arguing same as in Case~1, we construct a geodesic ray $\wt\gamma:[0,\infty)\to U$ fo $(U,K_U)$ such that $\wt\gamma\subset W$. Therefore, by Theorem~\ref{Thm:Global_to_loc_vis}, we get that $\wt\gamma$ is a $(1,\k)$-almost-geodesic of $(\OM,K_\OM)$ for some $\k>0$. Also, by following same arguments as in the Case~1, we get that there exists a constant $C^{\prime\prime}>0$ such that for all $t\geq 0$
\[
K_\OM(\wt\gamma(t),\wh\gamma(t)) <C^{\prime\prime}.
\]
By Theorem~\ref{Thm:Geo_horosphere_of_big_domain_is_subset_of_geod_of_small_domain_upto_some_compact_set_and-constant} and Lemma~\ref{Lem:Relation_between_origin_geo_horo_and_truncated_geo_horo}, there exists a constant $C_2>0$ such that
\[
H_\OM(\gamma,R)\subseteq H_U(\wt\gamma,e^{2T^\prime}C_2R).
\]
Now, by~\eqref{Eq:Small_dom_geo_horo_int_at_sing}, we get
\[
\clos{H_\OM(\gamma,R)}\cap\bdy \OM \subseteq \clos{H_U(\wt\gamma,e^{2T^\prime}C_2R)}\cap\bdy \OM =\{x^\prime\}.
\]

Therefore, for every geodesic ray $\gamma:[0,\infty)\to\OM$ of $(\OM,K_\OM)$, there exists $x\in\bdy\OM$ such that $\clos{H_{\OM}(\gamma,R)}\cap\bdy\OM=\{x\}$ for all $R>0$. This completes the proof of proposition.
\end{proof}
\end{example}


{\bf Declaration:} No artificial intelligence tools were used to obtain any of the statements and their proofs presented in this paper.  All results are proved by the authors.
\vspace{0.8em}

{\bf Conflict of interest:} The authors declare that there is no conflict of interest regarding
the publication of this paper.
\vspace{0.8em}

{\bf Acknowledgments:} The first named author was supported by DST-INSPIRE Fellowship (IF230372), New Delhi, India. The second named author was supported by the institute seed grant IIT/SRIC/BS/RP551/2025-2026 IIT Bhubaneswar.


\begin{thebibliography}{99}
\bibitem{Abate:Horo_and_itera}
M.\,Abate, {\em Horospheres and iterates of holomorphic maps}, Math. Z. 198 (2) (1988) 225--238.

\bibitem{AbateBook}
M.\,Abate, {\em Iteration theory of holomorphic maps on taut manifolds}, Mediterranean Press, Cosenza, 1989.

\bibitem{Abate:Ite_theory_and_comp_div}
M.\,Abate, {\em Iteration theory, 
compactly divergent sequences and commuting holomorphic maps}, Ann.\,Scuola Norm.\,Sup.\,Pisa Cl.\,Sci., 18, 167-191, (1991).

\bibitem{AR2014}
M.\,Abate and J.\,Raissy, {\em Wolff-{D}enjoy theorems in nonsmooth convex domains}, Ann. Mat. Pura Appl. (4), 193 (5) (2014) 1503--1518.

\bibitem{AbateNewbook}
M.\,Abate, {\em Holomorphic dynamics on hyperbolic Riemann surfaces}, De Gruyter Stud. Math., 89 de Gruyter, Berlin, 356 pp., 2023.

\bibitem{Beardon1}
A.\,F.\,Beardon, {\em Iteration of contractions and analytic maps}, J.\,London\,Math.\,Soc. \textbf{41} (1990), 141--150.

\bibitem{Beardon2}
A.\,F.\,Beardon, {\em The dynamics of contractions}, Ergod.\,Th.\,$\&$\,Dynam.\,Sys. \textbf{17} (1997)  1257–-1266.

\bibitem{Benini-Bracci}
A.\,M.\,Benini and F.\,Bracci, {\em The Denjoy-Wolff Theorem in simply connected domains}, Trans.\,Amer.\,Math.\,Soc.\, (2026).

\bibitem{Bharali_Maitra} G.\,Bharali and A.\,Maitra, {\em A weak 
notion of visibility, a family of examples, and Wolff--Denjoy theorems}, 
Ann. Sc. Norm. Super. Pisa Cl. Sci. (5) 22 (2021), 195--240.

\bibitem{BZ2017}
G.\,Bharali and A.\,Zimmer, {\em Goldilocks domains, a weak notion of 
visibility, and applications}, Adv. Math. 310 (2017), 377--425.

\bibitem{BZ2023} G.\,Bharali and A.\,Zimmer,
{\em Unbounded visibility domains, the end compactification, and
applications}, Trans. Amer. Math. Soc. 376 (2023), no.~8, 5949--5988.

\bibitem{BCM}
F.\,Bracci, M.\,D.\,Contreras and S.\,D{\'i}az-Madrigal, {\em Continuous 
semigroups of holomorphic self-maps of the unit disc}, Springer 2020.

\bibitem{Bracci-Gaussier}
F.\,Bracci, H.\,Gaussier, {\em Horosphere Topology}, Ann.~Scuola~Norm.~Sup.~di~Pisa, Cl. Sci. 20,\,1,\, 239--289,\, (2020).

\bibitem{BGNT2024}
F.\,Bracci, H.\,Gaussier, N.\,Nikolov, and P.\,J.\,Thomas, 
{\em Local and global visibility and Gromov hyperbolicity of domains with 
respect to the Kobayashi distance}, Trans.~Amer.~Math.~Soc.~{\bf 377} (2024), no.~1, 471--493.

\bibitem{Bracci_Gaussier_Zimmer}
F.\,Bracci, H.\,Gaussier and A.\,Zimmer, {\em Homeomorphic extension  of quasi-isometries for convex domains in $\C^d$ and iteration theory}, Math.\,Ann. 379, 691–718 (2021).

\bibitem{BNT2022}
F.\,Bracci, N.\,Nikolov and P.\,J.\,Thomas, {\em Visibility of Kobayashi geodesics 
in convex domains and related properties}, Math.\,Z.,\,\textbf{301}(2), (2022), 2011--2035.

\bibitem{BB1}
F.\,Bracci and A.\,Y.\,{\"O}kten, {\em Some open questions and conjectures about visibility and iteration in bounded convex domains in $\C^{N}$}, Complex Analysis and its Synergies, (2025).

\bibitem{BO_Failure}
F.\,Bracci and A.\,Y.\, \"{O}kten, {\em On the failure of the Denjoy-Wolff Theorem in convex domains}, Proc.\,Amer.\,Math.\,Soc., (2026).

\bibitem{Bridson_Haefliger}
M.\,R.\,Bridson and A.\,Haefliger, {\em Metric Spaces of Non-positive Curvature}, Grundlehren 
der mathematischen Wissenschaften \textbf{319}, Springer-Verlag Berlin, 1999.

\bibitem{BKR}
M.\,Budzy\`{n}ska, T.\,Kuczumov and S.\,Reich {\em Theorems of Dennjoy-Wolff type}, Ann.\,Mat.\,Pura\,Appl. (4) 192 (2013), 621--648.

\bibitem{CGMS2024}
V.\,S.\,Chandel, S.\,Gorai, A.\,Maitra and A.\,D.\,Sarkar, {\em Visibility property in one and several variables and its applications}, (2024), \url{https://arxiv.org/abs/2406.15298}.

\bibitem{CMS2024}
V.\,S.\,Chandel, A.\,Maitra and A.\,D.\,Sarkar, {\em Notions of visibility with respect to the
Kobayashi distance: comparison and applications}, Ann. Mat. Pura Appl. 203, 475--498, 2024.

\bibitem{Collingwood-Lowater_Book}
E.\,F.\,Collingwood and A.\,J\,Lohwater, {\em The theory of cluster sets}, Cambridge Tracts in Mathematics and Mathematical Physics, No. 56. Cambridge University Press, Cambridge, 1966.

\bibitem{Denjoy} 
A.\,Denjoy, {\em Sur l'it\'eration des fonctions analytiques}, C. R. Acad. Sci. Paris 182 (1926), 255--257.

\bibitem{EON}
P.\,Eberlein and B.\,O'Neill, {\em Visibility manifolds}, Pacific J. Math. \textbf{46} (1973), 45--109.

\bibitem{Herve}
M.\,Herv\'e, {\em Quelques propri\'et\'es des applications analytiques d'une boule \`a m dimensions dan ellem\^eme}, J. Math. Pures Appl. 9 (42) (1963), 117--147.

\bibitem{JarPfl} 
M.\,Jarnicki and P.\,Pflug, {\em Invariant Distances and Metrics in Complex Analysis}, 
2nd extended edition, De Gruyter Expositions in Mathematics {\bf 9}, 2013.

\bibitem{Karlsson}
A.\,Karlsson, {\em Non-expanding maps and busemann functions}, Ergodic Theory Dyn.\,Syst. 21, 1447–1457 (2001)

\bibitem{Masanta}
R.\,Masanta, {\em Visibility domains relative to the Kobayashi distance in complex manifolds},  Trans. Amer. Math. Soc., (2025).

\bibitem{NOT2024}
N.\,Nikolov, A.\,Y.\,{\"O}kten and P.\,J.\,Thomas, {\em Local and 
global notions of visibility with respect to Kobayashi distance: a comparison}, 
Ann. Polon. Math.~{\bf 132} (2024), no.~2, 169--185.

\bibitem{PapaBook}
A.\,Papadopoulos, {\em Metric spaces, convexity and nonpositive curvature}, IRMA Lectures in Mathematics and Theoretical Physics.~{\bf 6} (2), European Mathematical Society, 2014.

\bibitem{RoyKobDist} 
H.\,L.\,Royden, {\em Remarks on the Kobayashi metric}, Several Complex Variables, College Park, MD, USA,
1970, Part 2, Lecture Notes in Mathematics {\bf 185}.

\bibitem{ADS2023}
A.\,D.\,Sarkar, {\em Localization of the Kobayashi distance for any visibility domain}, J.~Geom.~Anal.~{\bf 33} (2023), no.~5, Paper No.~144, 16~pp.

\bibitem{Whyburn:topology}
G.\,T.\,Whyburn, {\em Analytic Topology}, American Mathematical Society Colloquium Publications, {\bf 28}, American Mathematical Society, New York 1942.

\bibitem{Wolff1}
J.\,Wolff, {\em Sur l'it\'eration des fonctions holomorphes dans une r\'egion, et dont les valeurs appartiennent \`a cette r\'egion}, C. R. Acad. Sci. Paris 182 (1926), 42--43.

\bibitem{Wolff2}
J.\,Wolff, {\em Sur l'it\'eration des fonctions born\'ees}, C. R. Acad. Sci. Paris 182 (1926), 200-201.
\end{thebibliography}
\end{document}